\documentclass[11pt]{amsart}
\usepackage{enumerate,amssymb}
\usepackage{amsopn,amsfonts,amsmath}
\usepackage[dvips]{graphicx}
\usepackage{float}
\usepackage[american]{babel}
\usepackage{enumerate}
\usepackage{hyperref}
\usepackage{fancyhdr}
\usepackage{color, colortbl}
\usepackage{float}
\allowdisplaybreaks
\numberwithin{equation}{section}
\def\E{{\mathbb E}}
\def\indic{{\rm {\large 1}\hspace{-2.3pt}{\large
l}}}

\def\R{{\mathbb R}}

\def\xP{{\mathbb P}}

\newtheorem{theorem}{Theorem}[section]
\newtheorem{proposition}[theorem]{Proposition}
\newtheorem{lemma}[theorem]{Lemma}
\newtheorem{corollary}[theorem]{Corollary}
\newtheorem{definition}[theorem]{Definition}

\newtheorem{assumption}[theorem]{Assumption}

\begin{document}
\title[]{High-dimensional instrumental
variables regression and confidence sets - v2/2012}

\author[Gautier]{Eric Gautier}
\address{CREST, ENSAE ParisTech,
3 avenue Pierre Larousse, 92 245 Malakoff Cedex, France.}
\email{\href{mailto:eric.gautier@ensae-paristech.fr}{eric.gautier@ensae-paristech.fr}}

\author[Tsybakov]{Alexandre B. Tsybakov}
\address{CREST, ENSAE ParisTech,
3 avenue Pierre Larousse, 92 245 Malakoff Cedex, France.}
\email{\href{mailto:alexandre.tsybakov@ensae-paristech.fr}
{alexandre.tsybakov@ensae-paristech.fr}}

\date{First version December 2009, this version: November 2012.}
\thanks{\emph{Background facts}:  This was a revision of the arXiv preprint \href{http://arxiv.org/pdf/1105.2454v1.pdf}{arXiv:1105.2454v1} which was online 
\href{https://sites.google.com/view/eric-gautier/papers?authuser=1}{here}.  A link is given in footnote 2 of \href{http://arxiv.org/abs/1105.2454v5}{arXiv:1105.2454v5}. % but for more exposure we make this file more widely available. 
(3.5) is a different {\em STIV} estimator from the one in \href{http://arxiv.org/pdf/1105.2454v1.pdf}{arXiv:1105.2454v1} where, instead of one conic constraint, there are as many conic constraints as moments (instruments) allowing to use more directly moderate deviations for self-normalized sums. The idea first appeared in formula (6.5) in \href{http://arxiv.org/pdf/1105.2454v1.pdf}{arXiv:1105.2454v1} when some instruments can be endogenous. This was called the {\em STIV} until the \href{https://www.dropbox.com/s/otioun8wzz7htzq/GTrevision3.pdf?dl=0}{2014 version}. Eventually, the original {\em STIV} was put back in the v2 on arXiv because we did not agree on the usefulness of the additional complication. %It was not clear this estimator was superior and it is much more difficult to compute. 
For reference and to avoid confusion with the {\em STIV} estimator, this estimator should be called {\em C-STIV} 
%Much later, the idea was used in \href{https://arxiv.org/abs/1708.08353v1}{arXiv:1708.08353v1}  in the context of a model with error in variables. Unexpectedly, Section 3.3 of \href{https://arxiv.org/abs/1708.08353v2}{arXiv:1708.08353v2} had a whole section on the {\em C-STIV} in the particular case without endogenous regressors, with various mentions that it can be generalized to the instrumental variables context. % without the first author of this paper responsible for this alternative STIV estimator, claiming this is a new idea. 
%After clarification that this was the  {\em STIV} from 2012 to 2014, \href{https://www.dropbox.com/s/otioun8wzz7htzq/GTrevision3.pdf?dl=0}{2014 version} being shared to offer to write a separate paper on the {\em C-STIV} that includes both authors of the original idea, and that the idea was presented during the one-to-one meetings % to the second author of \href{https://arxiv.org/abs/1712.08102v1.pdf}{arXiv:1712.08102v1} 
%after the presentation at Harvard-MIT in 2011 (see below), \href{https://arxiv.org/abs/1708.08353v2}{arXiv:1708.08353v2}  was removed from arXiv. %from arXiv 
%Still, % by Alex Belloni. 
%in  \href{https://arxiv.org/abs/1708.08353v3}{arXiv:1708.08353v3}, this paper is not cited. % and (9.22) is cited instead of (6.5) from \href{http://arxiv.org/pdf/1105.2454v1.pdf}{arXiv:1105.2454}. 
%Later, in \href{https://arxiv.org/abs/1712.08102v1.pdf}{arXiv:1712.08102v1}, %was put online and 
%the {\em C-STIV} is reintroduced again (see formula (2.4)) with two authors of \href{https://arxiv.org/abs/1708.08353v3}{arXiv:1708.08353v3} but without any of this one or citation.
}
\thanks{\emph{Keywords}:  Instrumental variables, sparsity, {\em STIV} estimator,
endogeneity, high-dimensional regression, conic programming,
heteroscedasticity, confidence regions, non-Gaussian errors,
variable selection, unknown variance, sign consistency.}
\thanks{We thank James Stock and three anonymous
referees for comments that greatly improved this paper.  We also
thank Don Andrews, Ivan Canay, Victor Chernozhukov, Azeem Shaikh and
the seminar participants at Brown, Bocconi, CEMFI, CREST, Compi\`egne,
Harvard-MIT, Institut Henri Poincar\'e, LSE, Paris 6 and 7,
Princeton, Queen Mary, Toulouse, Valpara\'iso, Wharton, Yale, as
well as participants of SPA, Saint-Flour, ERCIM 2011, the 2012
CIREQ conference on High Dimensional Problems in Econometrics and 4th
French Econometrics Conference for
helpful comments.}

\begin{abstract}
We propose an instrumental variables method for inference in
structural linear models with endogenous regressors in the
high-dimensional setting where the number of possible regressors $K$
can be much larger than $n$.  Our new procedure, called {\em STIV}
(Self Tuning Instrumental Variables) estimator, is realized as a
solution of a conic optimization program.  We allow for partial
identification and very weak distributional assumptions including
heteroscedasticity.  We do not need prior knowledge of the variances
of the errors.  A key ingredient is sparsity, {\it i.e.}, the vector
of coefficients has many zeros. The main result of the paper is
nested confidence sets, sometimes non-asymptotic, around the
identified region under various level of sparsity.  In the presence
of very many instruments, a number exponential in the sample size,
too many non-zero coefficients, or when there is an endogenous
regressor which is only instrumented by instruments which are too weak, our
confidence sets have infinite volume.  We show that a variation on the {\em STIV}
estimator is a new robust method to estimate low dimensional models
without sparsity and possibly many weak instruments. We obtain
rates of estimation and show that, under appropriate assumptions, a
thresholded {\em STIV} estimator correctly selects the non-zero
coefficients with probability close to 1. In a non sparse setting,
we obtain a sparse oracle inequality that shows how well the
regression coefficients are estimated when the model can be well
approximated by a sparse model. In our {\em IV} regression setting,
the standard tools from the literature on sparsity, such as the
restricted eigenvalue assumption are inapplicable.  Therefore, for
our analysis we develop a new approach based on data-driven
sensitivity characteristics.  We also obtain confidence sets, when
the endogenous regressors in the high-dimensional structural equation
have a sparse reduced form and we use a
two-stage procedure, akin to two-stage least squares.
We finally study the properties of a
two-stage procedure to deal with the detection of non-valid
(endogenous) instruments.
\end{abstract}
\maketitle

\section{Introduction}\label{s1}
In this article we consider a structural model of the form
\begin{equation}\label{estruct}
y_i=x_i^T\beta+u_i, \quad i=1,\dots,n,
\end{equation}
where $x_i$ are vectors of explanatory variables of dimension
$K\times1$, $u_i$ is a zero-mean random error possibly correlated
with some or all regressors.
We consider the problem of inference on the structural parameter $\beta$
from $n$ independent, not necessarily identically distributed,
realizations $(y_i,x_i^T,z_i^T)$, $i=1,\dots,n$.  This allows for
heteroscedasticity.  We denote by $x_{ki}$, $k=1,\dots, K$, the
components of $x_i$.  The regressors $x_{ki}$ are called endogenous
if they are correlated with $u_i$ and they are called exogenous
otherwise.  It is well known that
endogeneity occurs when a regressor correlated both
with $y_i$ and regressors in the model is unobserved; in the
errors-in-variables model when the measurement error is independent
of the underlying variable; when a regressor is determined
simultaneously\footnote{This is the case where $\beta$ is actually structural.}
with the response variable $y_i$.
A random vector $z_i$ of dimension
$L\times 1$ is called a vector of instrumental
variables (or instruments) if it satisfies
\begin{equation}\label{einstr}
\forall i=1,\hdots,n,\ \E[z_iu_i]=0,
\end{equation}
where $\E[\,\cdot\,]$ denotes the expectation.
%When the data is a realization of identically distributed random
%vectors\footnote{For ease of discussion on identification.},
Having access to instrumental variables makes it possible to identify the
vector $\beta$ when $\E[z_ix_i^T]$ has full column rank.
The case where $K\ge L$ is typically a case where identification
fails (see Example 6 below).
%This can occur when a discrete endogenous regressor is
%instrumented by a discrete variable that has less points in its
%support than the endogenous regressor or when we want to allow the
%instruments to have a direct effect on the outcome (see, {\em e.g.},
%Koles\'ar, Chetty, Friedman, et al. (2011)).

Though the method of this paper can be applied in standard setups\footnote{This
is also very attractive because it is a feasible, non combinatoric, method, unlike a
variable selection method like BIC.},%, or recent proposals that rely on nonconvex optimization.},
we are mainly interested in the more challenging high-dimensional
setting where $K$ can be much larger than $n$
and one of the following two assumptions is
satisfied:
\begin{itemize}
\item[(i)] only few coefficients $\beta_k$ are non-zero ($\beta$ is {\it
sparse}),
\item[(ii)] $\beta$ can be well approximated by a sparse vector
($\beta$ is {\it approximately sparse}).
\end{itemize}
The first assumption means that $K$ is not the actual number of regressors
but the number of potential regressors.  All of the regressors do not necessarily
appear in the true model underlying \eqref{estruct}.
The second assumption is likely to be satisfied in a wage
equation because in some data sets there are many variables that
have a nonzero effect on wage but their marginal effect
is too small to matter.
The large $K$ relative
to $n$ problem allows to deal with various important
situations encountered in empirical applications.  Here are some examples of
such applications.  In examples 3, 4 and 5 the most likely assumption is the approximate sparsity
assumption (ii).

\noindent{\bf Example 1. Economic theory is not explicit enough about which variables belong to the true model.}
Sala-i-Martin (1997)
and Belloni and Chernozhukov (2001b) give examples from development economics
where model selection is important.
In macroeconomics, development economics or
international finance, it is common to consider cross-country
regressions. Because $n$ is of the order of a few dozens,
it is important to allow for $K$ to be much larger than $n$.

\noindent{\bf Example 2. Rich heterogeneity.}
When there is a rich heterogeneity one usually wants to control
for many variables and possibly interactions, or to carry out a stratified analysis where
models are estimated in small population
sub-groups ({\it e.g.}, groups defined by
the value of an exogenous discrete variable). In both cases
$K$ can be large relative to $n$.

\noindent {\bf Example 3. Many endogenous regressors due to a type of non-parametric specification.}
This occurs when one considers a
structural equation of the form
\begin{align}
y_i&=f(x_{{\rm end},i})+u_i\label{eq:NP}\\
&=\sum_{k=1}^{K}\alpha_kf_k(x_{{\rm end},i})+u_i\notag\\
&=x_i^T\beta+u_i\label{eq:NP2}
\end{align}
and the stronger notion of exogeneity $\mathbb{E}[u_i|z_i]=0$ (zero
conditional mean assumption), $x_{{\rm end},i}$ is a low dimensional
vector of endogenous regressor, $x_i=(f_1(x_{{\rm
end},i}),\hdots,f_K(x_{{\rm end},i}))^T$,
$\beta=(\alpha_1,\hdots,\alpha_K)^T$ and $(f_k)_{k=1}^K$ are
functions from a dictionary\footnote{We use the word dictionary,
common in the machine learning community, because high-dimensional
methods allow to consider very many series terms and to mix bases.
This is not exactly the nonparametric {\em IV} setup because in a
truly nonparametric model there would be an extra approximation
error term in \eqref{eq:NP2}, of order at most $n^{-1/2}$ when $K$
is very large and we make minimal smoothness assumptions.  High
dimensional methods allow in several cases to obtain an adaptive
estimation method (that does not need to know the smoothness of the
unknown function) but we will not discuss this aspect in this
article.}. Exogenous regressors could also be included in the right
hand side of \eqref{eq:NP} similar to a partial linear model (see
next example).  Belloni and Chernozhukov (2011b) gives the example
of a wage equation with many transformation of education to properly
account for nonlinearities. Another typical application is the
estimation of Engle curves where it is important to include
nonlinearities in the total budget (see, {\it e.g.}, Blundell, Chen
and Kristensen (2007)). When one estimates Engle curves using
aggregate data, $n$ is again usually of the order of a few dozens.
It is well known that education in a wage equation and total budget
in Engle curves are endogenous variables.

\noindent {\bf Example 4. Many exogenous regressors due to a type of
semi-parametric specification.} Consider a partially linear model of
the form
\begin{equation}\label{eq:SP}
y_i=x_{{\rm end},i}^T\beta_{\rm end}+f(x_{{\rm exo},i})+u_i\quad \mathbb{E}[u_i|x_{{\rm exo},i}]=0.
\end{equation}
$x_{{\rm exo},i}$ is a low dimensional vector of exogenous regressors.
If $f$ can be properly decomposed as a linear combination of functions
from a large enough
dictionary $(f_k)_{k=1}^{K_c}$, with coefficients $\alpha_k$, we obtain
\begin{align}
y_i&=x_{{\rm end},i}^T\beta_{\rm end}+\sum_{k=1}^{K_c}\alpha_kf_k(x_{{\rm exo},i})+u_i\\
&=x_i^T\beta+u_i\label{eq:SP2}
\end{align}
where $\beta=\left(\beta_{\rm end}^T,\alpha_1,\hdots,\alpha_{K_c}\right)^T$
and $x_i=\left(x_{{\rm end},i}^T,f_1(x_{{\rm exo},i}),\hdots,f_{K_c}(x_{{\rm exo},i})\right)^T$.
In that case, one is usually interested in the marginal effects $\beta_{\rm end}$
of  $x_{{\rm end},i}$ holding fixed $x_{{\rm exo},i}$, rather than the whole vector
$\beta$.

\noindent {\bf Example 5. Many control variables to justify the use of an instrument.}
Suppose that we are interested in the parameter $\overline{\beta}$ in
\begin{equation}\label{eq:control}
y_i=\overline{x}_i^T\overline{\beta}+v_i,
\end{equation}
where some of the variables in $\overline{x}_i$ are endogenous, but that we have at our disposal
a variable $z_i$ that we want to use as an instrument but which does not satisfy
$\mathbb{E}[z_iv_i]=0$.  Suppose that we also have observations of vectors of controls
$w_i$ such that $\mathbb{E}[v_i|w_i,z_i]=\mathbb{E}[v_i|w_i]$ (conditional mean independence).
Then we can rewrite \eqref{eq:control} as
\begin{equation}\label{eq:control2}
y_i=\overline{x}_i^T\overline{\beta}+f(w_i)+u_i
\end{equation}
where $f(w_i)=\mathbb{E}[v_i|w_i]$ and $u_i=v_i-\mathbb{E}[v_i|w_i,z_i]$ is such that $\mathbb{E}\left[z_iu_i\right]=0$\footnote{By the
law of iterated conditional expectations.}.
It yields
\begin{equation}\label{eq:control3}
y_i=\overline{x}_i^T\overline{\beta}+\sum_{k=1}^{K_c}\alpha_kf_k(w_i)+u_i
\end{equation}
when $f$ can be decomposed on $(f_k)_{k=1}^{K_c}$.
This model can be rewritten in the form \eqref{estruct} with
$\beta=\left(\overline{\beta}^T,\alpha_1,\hdots,\alpha_{K_c}\right)^T$
and $x_i=\left(\overline{x}^T,f_1(w_i),\hdots,f_{K_c}(w_i)\right)^T$.
%This is a linear model with endogenous regressors and $K$ large where $K$ is the sum
%of $K_c$ and the dimension of $\overline{\beta}$.
Again one is usually interested in the subvector $\overline{\beta}$.

\noindent {\bf Example 6. The instruments can have a direct effect on the outcome.}
%Assume here that the data is a realization of identically distributed random vectors.
Koles\'ar, Chetty, Friedman, et al. (2011) considers the case where
one wants to allow the instruments to have a direct effect on the
outcome, {\it i.e.}, to be potentially on the right hand side of
\eqref{estruct}.  We are thus in a setting where $K>L$ which implies
that $\beta$ is not identified from the moment conditions
\eqref{einstr}.  One typically needs exclusion restrictions ({\it
i.e.}, that some exogenous variables are instruments and do not
appear on the right hand side of \eqref{estruct}) which corresponds
to some coefficients being equal to zero in the specification where
all instruments can potentially appear on the right hand side of
\eqref{estruct}.  Exclusion restrictions is therefore a reason for
sparsity. This specification is more flexible as we assume that
there exists some exclusion restrictions without telling which one
in advance. In a setup where there are many (weak) instruments\footnote{As
exemplified in Angrist and Krueger (1991), under the stronger notion
of exogeneity based on the zero conditional mean assumption,
considering interactions or functionals of instruments can lead to a
large number of instruments.}, allowing the instruments to
potentially have a direct effect, implies that $K$ is large.

Statistical inference under the sparsity scenario when the dimension
is larger than the sample size is now an active and challenging
field.  The most studied techniques are the Lasso, the Dantzig
selector (see, {\em e.g.}, Cand\`es and Tao~(2007), Bickel, Ritov
and Tsybakov~(2009); more references can be found in the recent book
by B\"uhlmann and van de Geer~(2011), as well as in the lecture
notes by Koltchinskii~(2011), Belloni and Chernozhukov~(2011b)), and
agregation methods (see Dalalyan and Tsybakov~(2008), Rigollet and
Tsybakov~(2011) and the papers cited therein).  A central concern in
this literature is to propose methods that are computationally
feasible.  The Lasso for example is a convex relaxation of $l^0$
penalized least squares methods like BIC.  The last are $NP$-hard
and it is impossible in practice to consider the case where $K$ is
larger than a few dozens. The Dantzig selector is solution of a
simple linear program. In recent years, these techniques became a
reference in several areas, such as biostatistcs and imaging. Some
important extensions to model from econometrics have been obtained
by Belloni and Chernozhukov~(2011a) who study the $\ell_1$-penalized
quantile regression and give an application to cross-country growth
analysis and by Belloni, Chernozhukov and Hansen~(2010) who use the
Lasso to estimate the optimal instruments with an application to the
impact of eminent domain on economic outcomes. St\"adler, Buhlmann
and van de Geer (2010) studies the estimation of mixtures of
high-dimensional regressions, this is important in econometrics to
handle group heterogeneity. Caner~(2009) studies a Lasso-type GMM
estimator. Rosenbaum and Tsybakov~(2010) deal with the
high-dimensional errors-in-variables problem and discuss an
application to hedge fund portfolio replication. Belloni and
Chernozhukov (2011b) also presents several applications of
high-dimensional methods to economics. The high-dimensional setting
in a structural model with endogenous regressors that we are
considering here has not yet been analyzed.  Note that the direct
implementation of the Lasso or Dantzig selector fails in the
presence of a single endogenous regressor simply because the zero
coefficients in the structural equation \eqref{estruct} do not
correspond to the zero coefficients in a linear projection type
model.  We also obtain confidence sets in a high-dimensional
framework.

The main message of this paper is that, in model \eqref{estruct}
containing endogenous regressors, under point identification, the
high-dimensional vector of coefficients or a subset of important
coefficients can be estimated together with proper confidence sets
using instrumental variables.  In  partially identified settings, we
obtain confidence sets around the identified region under a sequence
of sparsity scenarios varying the upper bound $s$ on the number of
non-zero coefficients of $\beta$ (that we call {\em sparsity
certificate}).   Thus there is no restriction on the size of $K$.
When $K$ is too large, one looses identification. Our confidence
sets can have infinite volume when $s$ and/or $K$ is too large.
There is no restriction either on the strength of the instruments.
If a very large number of instruments is used (exponential in the
sample size), or when all the instruments are weak, the method
yields infinite volume confidence sets.  The price to pay for
including an irrelevant instrument in a preliminary set of $L$
instruments is just a factor of $\sqrt{\log(L+1)}/\sqrt{\log(L)}$ in
the size of the confidence sets.  The {\em STIV} estimator is thus a
method that is robust to weak instruments and can handle very many
({\it i.e.} exponential in the sample size\footnote{The terminology
very many is borrowed from Belloni, Chen, Chernozhukov et
al.~(2010).}) instruments. Also it is not required in principle to
know in advance which regressor is endogenous and which is not, one
needs a set of valid (exogenous) instruments.  This is achieved by
the {\em STIV} estimator (Self Tuning Instrumental Variables
estimator) that we introduce below. Based on it, we can also perform
variable selection.  Under various assumptions on the data
generating process we obtain non-asymptotic or asymptotic results.
We also provide meaningful bounds when either (i) or (ii) above
holds and $\log(L)$ is small compared to $n$.  We believe that a
non-asymptotic framework is the most natural setting to consider in
high dimensions.  We also restrict mostly our attention to a one
stage method.  This is even more justified in {\em IV} estimation of
structural equations with endogenous regressors and weak instruments
where inference usually relies on non-standard asymptotics (see,
{\it e.g.}, Stock, Wright and Yogo (2002) and Andrews and Stock
(2007) for a review and the references cited therein). Nelson and
Startz (1990a, b) contains a simulation study where, due to weak
instruments, the finite sample distribution of the two-stage least
squares estimator can be non normal and even bimodal.

The {\em STIV} estimator is an extension of the Dantzig selector of
Cand{\`e}s and Tao (2007). The results of this paper extend those on
the Dantzig selector (see Cand{\`e}s and Tao~(2007), Bickel, Ritov
and Tsybakov~(2009) and further references in B\"uhlmann and van de
Geer~(2011)) in several ways: By allowing for endogenous regressors
when instruments are available, by working under weaker sensitivity
assumptions than the restricted eigenvalue assumption of Bickel,
Ritov and Tsybakov~(2009), which in turns yields tighter bounds, by
imposing weak distributional assumptions, by introducing a procedure
independent of the noise level and by providing confidence
sets.

The {\em STIV} estimator is also very much inspired by the
Square-root Lasso of Belloni, Chernozhukov and Wang (2010) which
proposes a pivotal method independent of the variance of the errors
in the Gaussian linear model with fixed regressors.  The Square-root
Lasso is a very important contribution to the literature on the
Lasso where almost all articles require to know the variance of the
errors which is related to the degree of penalization required.  The
most common practice is to adjust the degree of penalization is to
use cross validation or BIC.
%St\"adler, Buhlmann and van de Geer (2010) is another method based
%on $l^1$ penalized log-likelihood, Antoniadis (2010) proposes a
%$l^1$-penalized Huber's loss estimator which is equivalent to the
%Square-root Lasso.
The {\em STIV} estimator adds extra linear constraints coming from
the restrictions \eqref{einstr} to the Square-root Lasso which
allows one to deal with endogenous regressors. The implementation of
the {\em STIV} estimator also correspond to solving a simple conic
optimization program. Confidence sets require as well solving at
most $2K^2$ linear programs for lower bounding the sensitivities
that we introduce below.  So our method is easy and fast to
implement in practice. Our confidence sets rely on moderate
deviations for self-normalized sums and either on the sparsity
certificates or on perfect model selection.  Perfect model selection
requires a separation from zero of the non-zero
coefficients\footnote{We believe that this is unavoidable for model
selection but presumably strong for some econometrics applications.
However this result is a useful addition to the current state of the
art in the theory of high-dimensional regression even with exogenous
regressors.}, while the sparsity certificate approach
does not. %Moderate deviations are by essence
%finite sample results. The confidence sets can sometimes be refined
%in two stage.
We are indebted to Belloni, Chernozhukov and Wang (2010) who are the
first to use moderate deviations  for self-normalized sums results
for high-dimensional regression (see also Belloni, Chen,
Chernozhukov and Hansen (2011)).  We will see that under some
distributional assumptions these can lead to confidence sets that
have finite sample coverage properties.

\section{Basic Definitions and Notation}\label{s2}
We set $\bold{Y}=(y_1,\dots,y_n)^T$, $\bold{U}= (u_1,\dots,u_n)^T$,
and we denote by $\bold{X}$ and $\bold{Z}$ the matrices of dimension
$n\times K$ and $n\times L$ respectively with rows $x_i^T$ and
$z_i^T$, $i=1,\hdots,n$.

The sample mean is denoted by $\E_n[\,\cdot\,]$. We use the notation
$$
\E_n[X_k^aU^b]\triangleq\frac1{n}\sum_{i=1}^n x_{ki}^au_i^b, \quad \E_n[Z_l^a
U^b]\triangleq\frac1{n}\sum_{i=1}^n z_{li}^a u_i^b,
$$
where $x_{ki}$ is the $k$th component of vector $x_i$, and $z_{li}$
is the $l$th component of $z_i$ for some $k\in\{1,\dots,K\}$,
$l\in\{1,\dots,L\}$, $a\ge 0, b\ge0$. Similarly, we define the
sample mean for vectors; for example, $\E_n[UX]$ is a row vector
with components $\E_n[UX_k]$. We also define the corresponding
population means:
$$
\E[X_k^aU^b]\triangleq\frac1{n}\sum_{i=1}^n \E[x_{ki}^au_i^b], \quad \E[Z_l^a
U^b]\triangleq\frac1{n}\sum_{i=1}^n \E[z_{li}^a u_i^b],
$$
and set, for $l=1,\hdots,L$,
$$
z_{l*}\triangleq\max_{i}|z_{li}|,\quad (x\bullet z)_{l}=\left\{\max_{k=1,\hdots,K}\E_n\left[\left(\frac{X_k}
{\E_n[X_k^2]^{1/2}}Z_l\right)^2\right]\right\}^{1/2}.
$$
We denote by ${\bf D}_{{\bf X}}$ the diagonal $K\times K$ matrices
with diagonal entries $\E_n[X_k^2]^{-1/2}$ for $k=1,\dots,K$. We use
the notation ${\bf D_Z^{(I)}}$ for the $L\times L$ matrix of entries
$(x\bullet z)_{l}^{-1}$ for $l\in I$ and $z_{l*}^{-1}$ otherwise.
The set $I$ is a subset of $\{1,\hdots,L\}$.  It always contains the
index of the instrument which is identically equal to 1.  The
corresponding coefficient $\beta_k$ is the usual constant in
\eqref{estruct}.  The matrices ${\bf D}_{{\bf X}}$ and ${\bf
D_Z^{(I)}}$ are normalization matrices.  The user specifies
the set $I$ based on the observed distribution.
For fast implementation of
our algorithm, $z_{l*}^{-1}$ is the preferred normalization.  For
heavy tail distributions, though, one may want to use the second
normalization.  The more heavy tail instruments we have, the more
computationally intensive our method will be.  Indeed, we will add a
conic constraint for each heavy tail instrument.  For a
vector $\beta\in\R^K$, let $J(\beta)=\{k\in\{1,\hdots,K\}:\
\beta_k\ne0\}$ be its support, {\em i.e.}, the set of indices
corresponding to its non-zero components $\beta_k$. We denote by
%$|\beta|_0$ the cardinality of the set $J(\beta)$ {\bf [Eric: non utilis\'e?]}
%and by
$|J|$ the cardinality of a set $J\subseteq \{1,\hdots,K\}$ and by
$J^c$ its complement: $J^c=\{1,\hdots,K\}\setminus J$.  The subset
of indices $\{1,\hdots,K\}$ corresponding to variables
in \eqref{estruct} that are known in advance to be exogenous and
serve as their own instruments is denoted by $J_{{\rm exo}}$.  There
might be more true exogenous variables in \eqref{estruct}, they are
excluded from this list of indices if they have not been included as
instruments.  The $\ell_p$ norm of a vector $\Delta$ is denoted by
$|\Delta|_p$, $1\le p\le\infty$. For
$\Delta=(\Delta_1,\dots\Delta_K)^T\in\R^K$ and a set of indices
$J\subseteq \{1,\ldots,K\}$, we consider $\Delta_J\triangleq
(\Delta_1\indic_{\{1\in J\}}, \ldots,\Delta_K\indic_{\{K\in
J\}})^T$, where $\indic_{\{\cdot\}}$ is the indicator function. For
a vector $\beta\in\R^K$, we set $\overrightarrow{{\rm
sign}(\beta)}\triangleq({\rm sign}(\beta_1),\hdots,{\rm
sign}(\beta_K))$ where
$${\rm sign}(t)\triangleq\left\{\begin{array}{ll}
                  1 & {\rm if\ }t>0 \\
                  0 & {\rm if\ }t=0 \\
                  -1 & {\rm if\ }t<0
                \end{array}\right.
$$
For $a\in \R$, we set $a_+\triangleq\max(0,a)$, $a_+^{-1}\triangleq
(a_+)^{-1}$, and $a/0\triangleq \infty$ for $a>0$. We adopt the
convention $0/0\triangleq0$ and $1/\infty \triangleq0$.

We denote by
$$\mathcal{I}dent=\left\{\beta:\ \mathbb{E}[z_i(y_i-x_i^T\beta)]=0\right\}$$
the identified region.  It is an affine space which is reduced to a
point when the model \eqref{estruct}-\eqref{einstr} is point
identified. It is possible to impose some restrictions like a known
sign or a prior upper bound on the size of the coefficients.
However, for simplicity, we will just consider the case when we know
an a priori upper bound  $s$ on the sparsity of $\beta$, {\em i.e.},
we know that $|J(\beta)|\le s$ for some integer $s$. We call this a
{\em sparsity certificate}. We introduce
$$\mathcal{B}_s=\mathcal{I}dent
\bigcap\left\{\beta:\ |J(\beta)|\le s\right\}.$$
Note that $\mathcal{B}_K=\mathcal{I}dent$.  Thus, considering
confidence sets around $\mathcal{B}_s$ allows to deal with both the
case when we rely on sparsity certificates and the case when we do
not.

\section{The {\em STIV} Estimator}\label{s3}
The sample counterpart of the moment conditions
(\ref{einstr}) can be written in the form
\begin{equation}\label{eq:struct}
\frac1n \bold{Z}^T(\bold{Y}-\bold{X}\beta)=0.
\end{equation}
This is a system of $L$ equations with $K$ unknown parameters. If
$L>K$, it is overdetermined.  In general ${\rm
rank}(\bold{Z}^T\bold{X})\le \min(K,L,n)$, thus when $L\le K$ or
when $n<K$ the matrix does not have full column rank.  Furthermore,
replacing the population equations (\ref{einstr}) by
(\ref{eq:struct}) induces a huge error when $L$, $K$ or both are
larger that $n$.  So, looking for the exact solution of
(\ref{eq:struct}) in the high-dimensional setting makes no sense.
However, we can stabilize the problem by restricting our attention
to a suitable ``small" candidate set of vectors $\beta$, for
example, to those satisfying the constraint
\begin{equation}\label{eq:constr}
\left|\frac1n \bold{Z}^T(\bold{Y}-\bold{X}\beta)\right|_{\infty}\le \tau,
\end{equation}
where $\tau >0$ is chosen such that (\ref{eq:constr}) holds for
$\beta$ in $\mathcal{B}_s$ with high probability.  We can then refine the
search of the estimator in this ``small" random set of vectors
$\beta$ by minimizing an appropriate criterion.  It is possible to
consider different small sets in \eqref{eq:constr}, however the use
of the sup-norm makes the inference robust in the presence of weak
instruments.  This will be clarified later.

In what follows, we use this idea with suitable modifications.
First, notice that it makes sense to normalize the matrix ${\bf Z}$.
This is quite intuitive because, otherwise, the larger the
instrumental variable, the more influential it is on the estimation
of the vector of coefficients. For technical reasons, we choose
normalization where we multiply ${\bf Z}$ by ${\bf D_Z^{(I)}}$.  The
constraint (\ref{eq:constr}) is modified as follows:
\begin{equation}\label{eq:constr1}
\left|\frac1n {\bf D_Z^{(I)}}\bold{Z}^T(\bold{Y}-\bold{X}\beta)\right|_{\infty}\le \tau.
\end{equation}

Along with the constraint of the form (\ref{eq:constr1}), we include
more constraints to account for the unknown (average in the case of
heteroscedasticity) level $\sigma$ of the ``effective noise"
$z_{il}u_i$.
%  in particular, if the errors $u_i$
% are i.i.d., $\sigma^2$ corresponds to their unknown variance.
Specifically, we say that a pair $(\beta,\sigma)\in\R^K\times\R^+$
satisfies the {\em IV}-{\em constraint} if it belongs to the set
\begin{equation}\label{IVC}
\widehat{\mathcal{I}}^{(I)}\triangleq\left\{(\beta,\sigma):\
\beta\in\R^K,\ \sigma>0,\ \left|\frac1n {\bf D_Z^{(I)}}\bold{Z}^T(\bold{Y}
-\bold{X}\beta)\right|_{\infty}\le \sigma
r,\ ({\bf D_Z^{(I)}})_{ll}^2\widehat{Q}_l(\beta)\le \sigma^2, \forall l\in I\right\}
\end{equation}
for some $r>0$ (specified below), and
$$\widehat{Q}_l(\beta)\triangleq\frac1n\sum_{i=1}^nz_{li}^2(y_i-x_i^T\beta)^2.$$
Note that the instrument $z_{li}=1$ for all $i=1,\hdots,n$ belongs to $I$ and
the corresponding value of
$\widehat{Q}_l(\beta)$ is
$$\widehat{Q}_l(\beta)=\frac1n\sum_{i=1}^n(y_i-x_i^T\beta)^2.$$
%which for $\beta$ in $\mathcal{B}_s$ corresponds to the level of the noise.

\begin{definition}
We call the {\em STIV} estimator any solution
$(\widehat{\beta}^{(c,I)},\widehat{\sigma}^{(c,I)})$ of the following minimization
problem:
\begin{equation}\label{IVS}
\min_{(\beta,\sigma)\in \widehat{\mathcal{I}}^{(I)}} \left(\,\left|{\bf
D}_{{\bf X}}^{-1}\beta\right|_1+c\sigma \right),
\end{equation}
where $0<c<1$. \end{definition} In this formulation the instruments
show up both in the {\it IV}-constraint and in the penalization
through $\sigma$.  As discussed above, the {\it IV}-constraint
includes the constraints coming from the moment conditions
\eqref{einstr} accounting for instrument exogeneity.  We use
$\widehat{\beta}^{(c,I)}$ as an estimator of $\beta$ in
$\mathcal{B}_s$. Finding the {\em STIV} estimator is a conic
program; it can be efficiently solved, see Section \ref{s900}. Note
that the {\em STIV} estimator is not necessarily unique. Minimizing
the $\ell_1$ criterion $\left|{\bf D}_{{\bf X}}^{-1}\beta\right|_1$
is a convex relaxation of minimizing the $\ell_0$ norm, {\it i.e.},
the number of non-zero coordinates of $\beta$. This usually ensures
that the resulting solution is sparse. The term $c\sigma$ is
included in the criterion to prevent from choosing $\sigma$
arbitrarily large; indeed, the {\em IV}-constraint does not prevent
from this.  The matrix ${\bf D}_{{\bf X}}^{-1}$ arises from
re-scaling of ${\bf X}$. It is natural to normalize the regressors
by their size in a procedure that can do variable selection.  This
way changing units does not change which coefficients are found to
be zeros or non-zeros.  For the particular case where ${\bf Z}={\bf
X}$ and ${\bf D_Z^{(I)}}={\bf D}_{{\bf X}}$ is the diagonal matrix
with entries $\left(\mathbb{E}_n[X_k^2]^{-1/2}\right)_{k=1}^K$, the
{\em STIV} estimator provides an extension of the Dantzig selector
to the setting with unknown variance of the noise.

In this particular case, the {\em STIV} estimator can also be
related to the Square-root Lasso of Belloni, Chernozhukov and
Wang~(2010), which solves the problem of unknown variance in
high-dimensional regression with deterministic regressors and i.i.d.
errors.  The definition of {\em STIV} estimator contains the
additional constraint \eqref{eq:constr1}, which is not present in
the conic program for the Square-root Lasso.  This is due to the
fact that we have to handle the endogeneity.

\section{Summary of the Main Results}\label{s40}
The main result of this paper are nested confidence sets around the
regions $(\mathcal{B}_s)_{s\le K}$.  We allow for
set-identification, arbitrarily weak instruments, for a situation
where $K>L$ (for example when all the instruments can have a direct
effect on the outcome or for various combinations of discrete
regressors and instruments).  The dimension $K$ could be as large as
we want, even larger than an exponential in the sample size, in that
case the model \eqref{estruct} together with \eqref{einstr} becomes
set-identified.  Also, in finite samples and when $L$ and/or $K$ are
much larger than $n$, we can obtain finite volume confidence sets
for low values of $s$ and infinite volume confidence sets for larger
values of $s$, even in point-identified settings. These results are
presented in Section \ref{s5}.  We first present various assumptions
that one can make on the data generating process.  These sometimes limit
the number of instruments.  We give the corresponding way to adjust the
constant $r$ in the definition of the set
$\widehat{\mathcal{I}}^{(I)}$ appearing in the constrained
optimization for the {\em STIV} estimator.

In all but the first scenario on the data generating process, 
our confidence sets are of the form:\\
For every $\beta$ in $\mathcal{B}_s$, with probability at least
$1-\alpha$, for any $c$ in $(0,1)$ and any set of indices $I$ containing the index of the instrument which is unity,
any solution $(\widehat\beta^{(c,I)},\widehat\sigma^{(c,I)})$ of the minimization problem (\ref{IVS}) satisfies
$$
\left|\left({\bf D_X}^{-1} (\widehat{\beta}^{(c,I)}-\beta)\right)_{J_0}\right|_p \le
\frac{2\widehat{\sigma}^{(c,I)}r}{\overline{\kappa}_{p,J_0}^{(c,I)}(s)}\left(1-\frac{r}{\kappa_{1}^{(c,I)}(s)}
\right)_+^{-1}\, \quad \forall \ p\in[1,\infty],\ \forall J_0\subset\{1,\hdots,K\},
$$
where the set $J_0$ is specified by the econometrician, depending on
which subset of the regressor he is interested in, and for all
$k=1,\dots,K$,
$$|\widehat{\beta}_k^{(c,I)}-\beta_k| \le
\frac{2\widehat{\sigma}^{(c,I)}r}{\mathbb{E}_n[X_k^2]^{1/2}\,
\kappa^{(c,I)*}_{k}(s) }\left(1-\frac{r}{\kappa_{1}^{(c,I)}(s)}\right)_+^{-1}.
$$
Under the first scenario on the data generating process the $1-\alpha$ probability
event depends on $I$ so that the statement is not uniform in $I$ but for fixed $I$.
The constants $\overline{\kappa}_{p,J_0}^{(c,I)}(s)$, $\kappa_{1}^{(c,I)}(s)$ or
$\kappa^{(c,I)*}_{k}(s)$ are easy to calculate, data-dependent, lower
bounds on the sensitivities that we introduce in Section \ref{s4}.
The sensitivities provide a generalization of the restricted
eigenvalues to non-symmetric and non-square matrices (see Section
\ref{s900} in the appendix for a comparison).  The most difficult
constant to calculate is $\kappa_{1}^{(c,I)}(s)$, it requires to solve
$2K^2$ linear programs. The confidence sets could be infinite
when $\kappa_{1}^{(c,I)}(s)\le r$.  We will see that this occurs either when
we have very many instruments or when there exists a regressor for which
all instruments are too weak instruments or $s$ is too large.
We will also discuss how the
constants $\kappa^{(c,I)*}_{k}(s)$ and $\kappa_{1}^{(c,I)}(s)$ compare
to the usual concentration parameter defined for low-dimensional
structural equations with homoscedastic Gaussian errors.  Recall that we
relax both heteroscedasticity and Gaussian errors.  When we are only interested in
a subset of the coefficients, our confidence sets are smaller than those that one would
obtain by projection of joint confidence sets.  They account for the a priori
upper bound on the sparsity $s$.  The confidence sets do not require that the non-zero
coefficients are large enough but require to make a stand and give a prior
upper bound $s$.  As we will see in Section~\ref{s54},
the {\em STIV} estimator usually estimates too many non-zero
coefficients when the underlying vector is sparse so that 
$|J(\widehat{\beta}^{(c,I)})|$ can be considered as
a conservative upper bound on the true sparsity.
For increasing values of $s$ one obtain nested confidence sets.
These become infinite when
$s$ approaches $K$ and $K$ is larger than $n$.
If one is interested in a few coefficients, it is possible to 
choose the value of $c$ in $(0,1)$ that yields the smaller 
confidence sets.  This is easily obtained by taking grid values for $c$.
Similarly, for all but the first scenario on the data generating process, 
it is possible to try several sets $I$ and choose according to the size of the 
confidence sets around the vector of the coefficients of interest.

In Section \ref{s5}, we present a sparse oracle inequality in the
case (ii) above where the underlying model is not sparse.  It shows 
that the {\em STIV} estimator estimates the coefficients as well as 
a method that would know in advance which is the best
sparse approximation in terms of bias/variance trade-off.

In Section~\ref{s54}, we present rates of convergence and model
selection results.  Indeed, the right hand-side in the above two
upper bounds are random and these inequalities do not tell, for
example, how well is the vector of coefficient estimated depending
on $L$ and $n$.  We also obtain that, if the absolute values of the
coefficients $|\beta_k|$ are large enough on the support of $\beta$,
then with probability close to one $J(\beta) \subseteq
J(\widehat\beta^{(c,I)})$.  This is used to obtain a second type of
confidence sets where we plug-in an estimator of $J(\beta)$ in the
definition of the sensitivities to obtain a lower bound, and thus, a
proper data-driven upper bound on the estimation error.  Recall that
the first type of confidence sets uses sparsity certificates and
does not rely on such a separation from zero of the non-zero
coefficients assumption.  The fact that, with probability close to
one $J(\beta) \subseteq J(\widehat\beta^{(c,I)})$ is confirmed
in the simulation study where the {\em STIV}
estimator usually selects too many regressors.  This result
allows to recover exactly, with probability close to 1, the true
support of the vector of coefficients
as well as the sign of the coefficients by a thresholding
rule.

Section \ref{s7} discusses some special cases and extensions.
Because the {\em STIV} estimator is a one stage method that is
robust to weak-instruments and allows for very many instruments, we
present a variation on the {\em STIV} estimator to estimate
low-dimensional ($K<n$) and non sparse structural models.  This is a
new method to construct confidence sets that is robust to weak
instruments, heteroscedasticity and non-Gaussian errors, and allows
for very-many-instruments.  We also present the refined properties
of our estimator when none of the regressors and instruments have
heavy tail distribution.  This simplifies its calculation and gives
finer results, especially in the very many or weak-instruments case.
Finally, we present the properties of a two-stage {\em STIV}
procedure with estimated sparse linear projection type instruments,
akin to two-stage least squares. Here the two stages are
high-dimensional regressions and the second stage has endogenous
regressors.  This is related to the literature on selection of
instruments and optimal instruments.  Here we do not touch on
optimality because it is an open question to define optimality with
a high-dimensional structural equation with endogenous regressors.
This method only works if the endogenous regressors have a sparse
reduced form.  In the case of an approximately sparse reduced form
(case (ii) above) we would obtain rates of convergence but our
method do not yield confidence sets.

In Section \ref{s6}, we consider the following model
\begin{align*}
&y_i=x_i^T\beta+u_i,\\
&\E\left[z_iu_i\right]=0,\\
&\E\left[\overline{z}_iu_i\right]=\theta.
\end{align*}
Here, $z_{li}$ for $i=1,\hdots,n$ and $l=1,\hdots,L$, is a set of
instruments that are known in advance to be valid (exogenous) and
${\overline{z}}_i$ for $i=1,\hdots,n$ and $l=1,\hdots,L$, is a
second set of instruments.  One wants to decide which instrument
from this second list is valid and which is invalid (endogenous). In
this setup, we again allow for the dimensions of $x_i$, $z_i$ an
$\overline{z}_i$ to be much larger than $n$.

Section \ref{s8} discusses way to implement our algorithm and to
calculate the constants that drive the size of the confidence sets.
This is based either on conic or linear programs.  These are readily
available in many classical softwares.  In this section we also
present a simulation study.  All the proofs are given in the
appendix.

\section{Sensitivity Characteristics}\label{s4}
In the usual linear regression in low dimension, when
$\bold{Z}=\bold{X}$ and the Gram matrix $\bold{X}^T\bold{X}/n$ is
positive definite, the sensitivity is given by the minimal
eigenvalue of this matrix.  In high-dimensional regression, the
theory of the Lasso and the Dantzig selector comes up with a more
sophisticated sensitivity analysis; there the Gram matrix cannot be
positive definite and the eigenvalue conditions are imposed on its
sufficiently small submatrices. This is typically expressed via the
restricted isometry property of Cand{\`e}s and Tao~(2007) or the
more general restricted eigenvalue condition of Bickel, Ritov and
Tsybakov~(2009).  In our structural model with endogenous
regressors, these sensitivity characteristics cannot be used, since
instead of a symmetric Gram matrix we have a rectangular matrix
$\bold{Z}^T\bold{X}/n$ involving the instruments.  More precisely,
we will deal with its normalized version
$$\Psi_n\triangleq\frac1n {\bf D_Z^{(I)}}\bold{Z}^T\bold{X}{\bf D}_{{\bf X}}.
$$
In general, $\Psi_n^{(I)}$ is not a square matrix.  For $L=K$, it is a
square matrix but, in the presence of at least one endogenous
regressor, $\Psi_n^{(I)}$ is not symmetric.

We now introduce some scalar sensitivity characteristics related to
the action of the matrix $\Psi_n^{(I)}$ on vectors in the cone
$$
C_J^{(c)}\triangleq\left\{\Delta\in\R^K:\
|\Delta_{J^c}|_1\le\frac{1+c}{1-c}|\Delta_{J}|_1 \right\},
$$
where $0<c<1$ is the constant in the definition of the {\em STIV}
estimator, $J$ is any subset of $\{1,\hdots,K\}$.  When the
cardinality of $J$ is small, the vectors $\Delta$ in the cone
$C_J^{(c)}$ have a substantial part of their mass concentrated on a set
of small cardinality.  We call $C_J^{(c)}$ the {\em cone of dominant
coordinates}.  The set $J$ that will be used later is the set
$J(\beta)$ which is small if $\beta$ is sparse.  The use of similar
cones to define sensitivity characteristics is standard in the
literature on the Lasso and the Dantzig selector (see, Bickel, Ritov
and Tsybakov~(2009)); the particular choice of the constant
$\frac{1+c}{1-c}$ will become clear from the proofs. It follows from
the definition of $C_J^{(c)}$ that
\begin{equation}\label{eq:technical}
|\Delta|_1\le \frac{2}{1-c}|\Delta_{J}|_1 \le
\frac{2}{1-c}|J|^{1-1/p}|\Delta_J|_{p}, \quad \forall \ \Delta\in C_J^{(c)}, \ 1\le p \le \infty.
\end{equation}
For $p\in[1,\infty]$, we define the {\em $\ell_p$ sensitivity} as
the following random variable:
$$\kappa_{p,J}^{(c,I)}\triangleq\inf_{\Delta\in C_{J}^{(c)}:\ |\Delta|_p=1}\left|\Psi_n^{(I)}\Delta \right|_{\infty}.$$
Similar, but different, quantities named cone invertibility factors
have been introduced in Ye and Zhang (2010).\\
Given a subset $J_0\subset\{1,\dots,K\}$ and $p\in[1,\infty]$,
we define the $l_p$-$J_0$-{\em block sensitivity} as
\begin{equation}\label{eq:technical2}
\kappa_{p,J_0,J}^{(c,I)}\triangleq\inf_{\Delta\in C_J^{(c)}:\
|\Delta_{J_0}|_p=1}\left|\Psi_n^{(I)}\Delta \right|_{\infty}.
\end{equation}
By convention, we set $\kappa_{p,{\varnothing} ,J(\beta^*)}^{(c,I)}=\infty$.
We use the notation $\kappa_{k,J}^{(c,I)*}$ for {\em coordinate-wise
sensitivities}, {\it i.e.}, for block sensitivities when $J_0=\{k\}$
is a singleton\footnote{They coincide for all values of $p\in[1,\infty]$.}:
\begin{equation}\label{eq:coordsens}
\kappa_{k,J}^{(c,I)*}\triangleq \inf_{\  \Delta\in C_J^{(c)}: \ \Delta_k= 1}
\left|\Psi_n^{(I)}\Delta \right|_{\infty}.
\end{equation}
Note that here we restrict the minimization to vectors $\Delta$ with
positive $k$th coordinate, $\Delta_k=1$, since replacing $\Delta$ by
$-\Delta$ yields the same value of $|\Psi_n^{(I)}\Delta|_{\infty}$.

The heuristic behind the sensitivity characteristics is the
following.  As we will see in the Appendix, for a fixed value of
$\beta$ in $\mathcal{B}_s$, we adjust $r$ in the definition of
$\widehat{\mathcal{I}}^{(I)}$ such that $\left(\beta,\max_{l\in
I}({\bf D_Z^{(I)}})_{ll}\sqrt{\widehat{Q}_l(\beta)}\right)$ belongs
to $\widehat{\mathcal{I}}^{(I)}$ on an event $E_{\alpha}$ of
probability $1-\alpha$. This  yields that for a proper $\tau$, on
$E_{\alpha}$,
\begin{equation}\label{eqCI}
|\Psi_n^{(I)}\Delta|_{\infty}\le\tau
\end{equation}
where $\Delta={\bf D}_{{\bf X}}^{-1}(\hat{\beta}^{(c,I)}-\beta)$. When
$\tau$ can be calculated directly from the data, this yields a
confidence region for $\beta$ in $\mathcal{B}_s$. Moreover, because
we minimize the objective function \eqref{IVS}, we will see that, on
the same event $E_{\alpha}$, $\Delta$ is constrained to belong to
the subset $C_{J(\beta)}^{(c,I)}$ of $\mathbb{R}^K$. The sensitivities allow
to deduce from \eqref{eqCI} and the cone condition, confidence
regions for various losses.  Suppose, for example, that one is
interested in a confidence region for the subvector of $\beta$ in
$\mathcal{B}_s$, which corresponds to the coefficients of indices in
$J_0\subset\{1,\hdots,K\}$, and considers the $l_p$-loss, then one
gets
$$|\Delta_{J_0}|_p\le\frac{\tau}{\kappa_{p,J_0,J(\beta)}^{(c,I)}}.$$
Indeed, this is trivial when $|\Delta_{J_0}|_p=0$ and, when it is not,
it immediately follows from
$$\frac{|\Psi_n^{(I)}\Delta|_{\infty}}{|\Delta_{J_0}|_p}
\ge\inf_{\widetilde{\Delta}:\ \widetilde{\Delta}\ne0,\
\widetilde{\Delta}\in C_{J(\beta)}^{(c,I)}}\frac{|\Psi_n^{(I)}
\widetilde{\Delta}|_{\infty}}{|\widetilde{\Delta}_{J_0}|_p}.$$  Remark
that, working with the coordinate-wise sensitivities
$\kappa^{(c,I)*}_{k,J(\beta)}$ is much better than projecting the region
\eqref{eqCI} onto the axes because it also takes into account the
cone condition and, thus, the sparsity of the underlying $\beta$. We
see that the sensitivities are quantities that are intrinsic to the
estimation properties of the {\em STIV} estimator\footnote{They are
also intrinsic to the Dantzig selector and the Lasso.}.  We show in
Section \ref{s91} that the assumption that the {sensitivities}
$\kappa_{p,J}^{(c,I)}$ are positive is weaker and more flexible than the
restricted eigenvalue (RE) assumption of Bickel, Ritov and
Tsybakov~(2009).  Unlike the RE assumption, it is applicable to
non-square non-symmetric matrices.  Another nice feature of the
sensitivities that we introduce is that lower bounds on these
sensitivities can be efficiently calculated (see below),
this opens the path to confidence statements in high-dimensional
regression.

The coordinate-wise sensitivities are measures of the strength of the
instruments.  The coordinate-wise sensitivity for index $k$ could be interpreted as {\em
restricted maximal partial empirical correlation between the
instruments and the regressor} $(x_{ki})_{i=1}^n$.  Indeed,
\eqref{eq:coordsens} can be written as
\begin{equation}\label{eRestMaxEmpirCorr}
\kappa_{k,J}^{(c,I)*}=\inf_{\lambda\in\mathbb{R}^{K-1}:\ \sum_{l\in J^c}|
\mathbb{E}_n[X_l^2]^{1/2}\lambda_l|\le \frac{1+c}{1-c}
\left(1+\sum_{l\in J\setminus\{k\}}|\mathbb{E}_n[X_l^2]^{1/2}\lambda_l|\right)}
\max_{l=1,\hdots,L}({\bf D_Z^{(I)}})_{ll}\left|\frac{1}{n}
\sum_{i=1}^n z_{li}\left(x_{ki}-x_{\{k\}^ci}^T\lambda\right)\right|.
\end{equation}
It is easy to check that \eqref{eRestMaxEmpirCorr} yields that
when $L<|J|$, $\kappa_{k,J}^{(c,I)*}=0$.
For confidence sets or rates of estimation statements, the set $J$ will
correspond to $J(\beta)$ for $\beta$ in $\mathcal{B}_s$ so that the
restriction in the infimum corresponds to vectors having most of
their mass on the support of the underlying vector of
coefficients $\beta$ from the set $\mathcal{B}_s$ (which includes the constant).
When an exogenous
variable serves as its own instrument, the coordinate-wise sensitivity
for that regressor should be large.  This is because there exists,
among the list of instruments, $(x_{ki})_{i=1}^n$ themselves.
Otherwise, the coordinate-wise sensitivity of a regressor is
small if all instruments have small restricted partial empirical
correlation with that regressor.  Because of the maximum, one good instrument is enough
to have a large coordinate-wise sensitivity.  It is small only if
all instruments are weak.\\
{\bf Example 7.} For comparison, consider a
structural equation with only one endogenous regressor
$(x_{k_{\rm end}i})_{i=1}^n$.  Assume that $(z_{li})_{i=1}^n$ for $l=1,\hdots,L$
are fixed and write the reduced form equation
\begin{equation}\label{eq:reducedform0}
x_{k_{\rm end}i}=\sum_{l\in J^c}\tilde{z}_{li}\zeta_{l}+\sum_{l\in J}\tilde{z}_{li}
\zeta_{l}+v_{i},\quad i=1,\dots,n,
\end{equation}
where $(\zeta_{l})_{l=1}^L$ are unknown coefficients, $J$
is the set of indices of the exogenous regressors that have non-zero
coefficients in \eqref{estruct} (this is possible when $|J|\le n$),
$\tilde{z}_{li}=z_{li}$ for $l\in J$ and $i=1,\hdots,n$, while
$(\tilde{z}_{li})_{i=1}^n$ are residuals from the regression of the
original instruments $(z_{li})_{i=1}^n$ for $l\in J^c$ on
$(z_{li})_{i=1}^n$ for $l\in J$.  Assume finally that
$(u_i,v_i)$ are i.i.d. and have a mean zero
bivariate normal distribution.  Denote by $\sigma_v^2$
the variance of $v_i$.
It is easy to check that $\lambda$ in $\mathbb{R}^{K-1}$ such that
$x_{\{k_{\rm end}\}^ci}^T\lambda=\sum_{l\in J}\tilde{z}_{li}
\zeta_{l}$ satisfies the constraint in \eqref{eRestMaxEmpirCorr}
which yields using \eqref{eq:reducedform0}
\begin{equation}\label{eInstrStrength}
\kappa_{k_{\rm end},J}^{(c,I)*}\ge\max\left(\max_{l\in J}
({\bf D_Z^{(I)}})_{ll}\left|\frac{1}{n}\sum_{i=1}^nz_{li}v_i\right|, \max_{l\in J^c}
({\bf D_Z^{(I)}})_{ll}\left|\frac{1}{n}\sum_{i=1}^nz_{li}v_i+\frac{1}{n}e_l^T{\bf \widetilde{Z}}_{J}^T
{\bf \widetilde{Z}}_{J^c}\zeta_{J^c}\right|\right)
\end{equation}
where $(e_l)_{l=1}^L$ is the canonical basis of $\mathbb{R}^L$.
Note that $\left|\frac{1}{n}\sum_{i=1}^nz_{li}v_i\right|$ is small for large $n$ because
the instruments are fixed and $v_i$ is mean zero.  This shares similarity
with the {\em concentration parameter} (see, {\it e.g}, Andrews and Stock (2007)),
defined as
$$\mu^2=\frac{\zeta_{J^c}^T{\bf \widetilde{Z}}_{J}^T
{\bf \widetilde{Z}}_{J^c}\zeta_{J^c}}{\sigma_v^2},$$
which is a classical measure of the strength of the set of instruments
for low dimensional structural equations, under homoscedasticity.
The quantity in \eqref{eInstrStrength} isolates the strength of each instrument.
Thus we see that, unlike the concentration parameter which considers the set of all instruments,
the coordinate-wise sensitivities depend on the strength of the best instrument.
Due to our rescaling of the instruments, which allows to relax a lot the usual distributional
assumptions and to handle heteroscedasticity, the coordinate-wise sensitivities and concentration
parameter are not comparable in terms of the dependence on $\sigma_v^2$.

A first reason to come up with lower bounds on the sensitivities is
that they depend on $J(\beta)$ and thus on the unknown $\beta$. A
second reason is to obtain an easy to calculate lower bound.  For
every sensitivity, we propose an easy to calculate lower bound.  The
following proposition will be useful.
\begin{proposition}\label{p4}
\begin{enumerate}[\textup{(}i\textup{)}]
\item\label{p410} Let $J,\widehat{J}$ be two subsets of
$\{1,\dots,K\}$ such that $J\subseteq \widehat{J}$. Then, for all
$J_0\subset\{1,\dots,K\}$, all $p\in[1,\infty]$, all $c$ in $(0,1)$ and
all set of indices $I$ containing the index of unity,
$\kappa_{p,J_0,J}^{(c,I)}\ge\kappa_{p,J_0,\widehat{J}}^{(c,I)}$;
\item\label{p411} For all $J_0\subset\{1,\dots,K\}$ all $p\in[1,\infty]$,
all $c$ in $(0,1)$ and
all set of indices $I$ containing the index of unity,
$\kappa_{p,J_0,J}^{(c,I)}\ge\kappa_{p,J}^{(c,I)}$.
\item\label{p42i} For all $p\in[1,\infty]$, all $c$ in $(0,1)$ and
all set of indices $I$ containing the index of unity,
\begin{equation}\label{k2}
\left(\frac{2|J|}{1-c}\right)^{-1/p}\kappa_{\infty,J}^{(c,I)}
\le\kappa_{p,J}^{(c,I)}\le \frac{2}{1-c}|J|^{1-1/p} \kappa_{1,J}^{(c,I)},
\end{equation}
and for all $J_0\subset\{1,\dots,K\}$, all $p\in[1,\infty]$, all $c$ in $(0,1)$ and
all set of indices $I$ containing the index of unity,
\begin{equation}\label{k3}
|J_0|^{-1/p}\kappa_{\infty,J_0,J}^{(c,I)}
\le\kappa_{p,J_0,J}^{(c,I)}\le |J_0|^{1-1/p}\kappa_{1,J_0,J}^{(c,I)};
\end{equation}
\item\label{p44i} For all $J_0\subset\{1,\dots,K\}$, $\kappa_{\infty,J_0,J}^{(c,I)}
=\min_{k\in J_0}\kappa^{(c,I)*}_{k,J}$.
\end{enumerate}
\end{proposition}
Results for the $l_p$-sensitivities are easily deduced from \eqref{p410}
and \eqref{p44i} by taking $J_0=\{1,\hdots,K\}$.
\eqref{p411} and \eqref{k3} allow to minorize the bloc sensitivities
when $J_0$ can be large and a direct calculation would be too difficult
from a numerical point of view.
%The proof of Proposition \ref{p4} is given in Section~\ref{s93}.

We can control $\kappa_{p,J_0,J(\beta)}^{(c,I)}$ without knowing $J(\beta)$
by means of {\em sparsity certificate}.  Assume that we have an
upper bound $s$ on the sparsity of $\beta$, {\em i.e.}, we know that
$|J(\beta)|\le s$ for some integer $s$ and that we use it as well to
define $\mathcal{B}_s$.  This does not require that non-zero
coefficients are large enough but just that a maximum of $s$ are
non-zero.  In view of (\ref{eq:technical}), if $|J|\le s$, then for
any $\Delta$ in the cone $C_J^{(c)}$ we have $|\Delta|_1\le
\frac{2s}{1-c}|\Delta|_{\infty}.$  Thus, for all $J$ such that
$|J|\le s$, we can bound the {coordinate-wise sensitivities} as
follows:
\begin{eqnarray}\label{eq:certif1}
\kappa_{k,J}^{(c,I)*}&\ge& \inf_{\Delta_k= 1,\,|\Delta|_1\le
a|\Delta|_{\infty}}\left|\Psi_n^{(I)}\Delta \right|_{\infty}
\\
\nonumber &\ge& \min_{j=1,\dots,K} \left\{ \min_{\ \ \Delta_k= 1,
|\Delta|_1\le a|\Delta_j|}\left|\Psi_n^{(I)}\Delta
\right|_{\infty}\right\} \triangleq \kappa_{k}^{(c,I)*}(s),
\end{eqnarray}
where $a=\frac{2s}{1-c}$\,.  For given $s$, this bound is data-driven
since the minimum in curly brackets can be computed by solving $2K$
linear programs (see Section \ref{s900}).  Then, using \eqref{p44i},
we can deduce a lower bound on $\kappa_{\infty,J_0,J}^{(c,I)}$
\begin{equation}\label{eq:kappa}
\kappa_{\infty,J_0,J}^{(c,I)}\ge \min_{k\in J_0}\kappa_{k}^{(c,I)*}(s).
\end{equation}
Using \eqref{k2} and \eqref{eq:kappa} we get computable lower bounds
for all $\kappa_{p,J}^{(c,I)}$, $p\in[1,\infty]$, which depend only on $s$
and on the data. In particular, for $|J|\le s$,
\begin{equation}\label{eq:kappa1}
\kappa_{1,J}^{(c,I)}\ge
\frac{1-c}{2s}\min_{k=1,\hdots,K}\kappa_{k}^{(c,I)*}(s)\triangleq
\kappa_1^{(c,I)}(s).
\end{equation}
This can thus be obtained by solving $2K^2$ linear programs.
Analogously to \eqref{eq:certif1}, the sparsity certificate approach
yields a bound for block sensitivities
\begin{eqnarray}\label{eq:certif2}
\kappa_{1,J_0,J}^{(c,I)}&\ge& \inf_{\ |\Delta_{J_0}|_1= 1,\,|\Delta|_1\le
a|\Delta|_{\infty}}\left|\Psi_n^{(I)}\Delta \right|_{\infty}
\\
\nonumber &\ge& \min_{j=1,\dots,K} \left\{ \min_{\ \
|\Delta_{J_0}|_1= 1, \ |\Delta|_1\le a|\Delta_j|}\left|\Psi_n^{(I)}\Delta
\right|_{\infty}\right\} \triangleq \kappa_{1,J_0}^{(c,I)}(s).
\end{eqnarray}
In Section \ref{s900} we show that the expression in curly brackets
in \eqref{eq:certif2} can be computed by solving $2^{|J_0|}$ linear
programs. Thus, the values $\kappa_{1,J_0}^{(c,I)}(s)$ can be readily
obtained for sets $J_0$ of small cardinality.  Otherwise we do not
advertise this lower bound but simply to use the following lower
bound for $p=1$.  For $|J|\le s$, \eqref{p411}, \eqref{k2} and
\eqref{eq:kappa} yield
\begin{equation}\label{eq:certifgen}
\kappa_{p,J_0,J}^{(c,I)}\ge\max\left(|J_0|^{-1/p}\min_{k\in J_0}\kappa^{(c,I)*}_k(s),
\left(\frac{1-c}{2s}\min_{k=1,\hdots,K}
\kappa^{(c,I)*}_k(s)\right)\right).
\end{equation}
This bound can be calculated by solving $2K^2$ linear programs even
for large sets $J_0$.

For the sake of completeness, we present an alternative to the sparsity
certificate approach.  It corresponds to computing
$\kappa_{1,J}^{(c,I)}$ and $\kappa_{k,J}^{(c,I)*}$ directly.
This is numerically feasible only for
$J$ of small cardinality.  Indeed, we show in Section \ref{s900}
that obtaining the
coordinate-wise sensitivities corresponds to solving $2^{|J|}$
linear programs.  Using (\ref{k2}) and (\ref{eq:kappa}), we obtain
computable lower bounds for all $\kappa_{p,J}^{(c,I)}$, $p\in[1,\infty]$.
The lower bounds are valid for any given index set $J$.
However, we
will need to compute the characteristics for the inaccessible set
$J=J(\beta)$, where $\beta$ is the unknown parameter from $\mathcal{B}_s$
on which we are making inference.
To circumvent this problem, we can plug in an estimator ${\widehat J}$
of $J(\beta)$.  For example, we can take ${\widehat
J}=J(\widehat\beta^{(c,I)})$. The confidence bounds remain valid whenever
$J(\beta) \subseteq {\widehat J}$, since then
$\kappa_{p,J(\beta)}^{(c,I)}\ge \kappa_{p,{\widehat J}}^{(c,I)},$ by Proposition
\ref{p4}~(i). Theoretical guarantees for the inclusion $J(\beta)
\subseteq J(\widehat\beta^{(c,I)})$ to hold with probability close to 1
require a separation from zero: $|\beta_k|$ is not too small on the
support of $\beta$ (see Theorem~\ref{t1b}~(iv)). On the other hand,
one typically observes in simulations that the relevant set
$J(\beta)$ is either estimated exactly or overestimated
by its
empirical counterpart ${\widehat J}=J({\widehat \beta}^{(c,I)})$, so that
the required inclusion is satisfied for such a simple choice of
${\widehat J}$.  Belloni and Chernozhukov (2010) study the property of
post model selection least squares estimation in the case of prediction loss,
they obtain theoretical bounds on the error made when we do not have
$J(\beta)\subset J({\widehat \beta}^{(c,I)})$.  We do not touch upon the very difficult
question of obtaining confidence sets with coverage at most $1-\alpha$ that
allow to deal with the possibility that $J(\beta)\subset J({\widehat \beta}^{(c,I)})$
is not satisfied.  We state that this type of confidence sets has
coverage $1-\gamma$ where $\gamma>\alpha$.  They are approximate  $1-\alpha$
confidence sets.  In our simulation study we obtain
that the confidence regions obtained via the sparsity certificate using
for $s=|J(\widehat{\beta}^{(c,I)})|$ are almost identical to the
second ones which also require the separation from zero.
Because the first method does not require such an assumption
and is feasible (non combinatoric) even when $|J(\beta)|$
is large, we strongly advertise the first approach based
on the sparsity certificate and to use a conservative sparsity
certificate $s$.  We also suggest to draw nested confidence sets
varying the degree of sparsity that we assume.

In the next proposition, we present a simple lower bound on
$\kappa_{p,J}^{(c,I)}$ for general $L\times K$ rectangular matrices
$\Psi_n^{(I)}$.  Its proof, as well as other lower bounds on
$\kappa_{1,J}^{(c,I)}$, can be found in Section \ref{s9}.  It is important
to note that adding rows to matrix $\Psi_n^{(I)}$ ({\it i.e.}, adding
instruments) increases the sup-norm $|\Psi_n^{(I)} \Delta|_\infty$, and
thus potentially increases the sensitivities $\kappa_{p,J(\beta)}^{(c,I)}$
and their computable lower bounds. This has a positive effect since
the inverse of the computable lower bounds of the sensitivities
drive the width of the confidence set for $\beta$ in
$\mathcal{B}_s$, see Theorem \ref{t1}.  Thus, adding instruments
potentially improves the confidence set, which is quite intuitive.
On the other hand, the price for adding instruments in terms of the
rate of convergence is only logarithmic in the number of
instruments, as we will see in the next section.

\begin{proposition}\label{p5}
Fix $J\subseteq\{1,\dots,K\}$.
Assume that there exist $\eta_1>0$ and $0<\eta_2<1$ such that
\begin{equation}\label{k2_0}
\forall
k\in J,\ \exists l(k):\ \left\{\begin{array}{l}
         |(\Psi_n^{(I)})_{l(k)k}|\ge\frac{\eta_1}{1-c}\,,\\
         \frac{\max_{k'\ne k}|(\Psi_n^{(I)})_{l(k)k'}|}
         {|(\Psi_n^{(I)})_{l(k)k}|}\le\frac{(1-\eta_2)(1-c)}{2|J|}\,.
       \end{array}\right.
\end{equation}
Then
$$\kappa_{p,J}^{(c,I)}\ge(2|J|)^{-1/p}(1-c)^{-1+1/p}\eta_1\eta_2.$$
\end{proposition}

%The proof of Proposition \ref{p5} is given in Section~\ref{s93}.

Assumption (\ref{k2_0}) is similar in spirit to the coherence
condition introduced by Donoho, Elad and Temlyakov~(2006) for
symmetric matrices, but it is more general because it deals with
rectangular matrices. Since the regressors and instruments are
random, the values $\eta_1$ and $\eta_2$ can, in general, be random.
Remarkably, for estimation of the coefficients of the endogenous
variables, it suffices to have a ``good" row of the matrix $\Psi_n^{(I)}$.
This means that it is enough to have, among all instruments, one
good instrument. The way the instruments are ordered is not
important. Good instruments correspond to the rows $l(k)$, for which
the value $|(\Psi_n^{(I)})_{l(k)k}|$ measuring the relevance of the
instrument for the $k$th variable is high. On the other hand, the
value $\max_{k'\ne k}|(\Psi_n^{(I)})_{l(k)k'}|$ accounting for the
relation between the instrument and the other variables should be
small. An instrument which is well ``correlated" with two variables
of the model is not satisfactory for this assumption.

%\newpage

\section{Confidence Sets and Sparse Oracle Inequality}\label{s5}
\subsection{Distributional Assumptions and Control of The Stochastic Part}\label{s51}
For every $\beta$ in $\mathcal{B}_s$, $u_i=y_i-x_i^T\beta$ is such
that, for every $l=1,\hdots,L$, $z_{li}u_i$ are independent and, for
every $i=1,\hdots,n$, $\mathbb{E}[z_{li}u_i]=0$. The size of the
constant $r$ in the definition of $\widehat{\mathcal{I}}^{(I)}$ is
directly related to the coverage probability of the confidence sets.
It should be adjusted so that with probability $1-\alpha$ (the
coverage probability),
$$\max_{l=1,\hdots,L}\frac{\left|\frac{1}{n}\sum_{i=1}^nz_{li}u_i\right|}
{\sqrt{\frac{1}{n}\sum_{i=1}^nz_{li}^2u_i^2}}\le r.$$
This is a sup-norm of so called {\em self-normalized sums}.
We propose different possible choices of  $r$ based on
different distributional assumptions (and sometimes an upper bound
on the number of instruments).

\noindent {\bf Scenario 1.} Suppose that the errors $u_i$ are
identically distributed, independent from the $z_i$'s, and of
distribution known up to the variance, then the quantiles of
$$\max\left(\max_{l\in I^c}\frac{\left({\bf D_Z^{(I)}}\right)_{ll}
\left|\frac{1}{n}\sum_{i=1}^nz_{li}u_i\right|}
{\sqrt{\frac{1}{n}\sum_{i=1}^nu_i^2}},\max_{l\in I}\frac{
\left|\frac{1}{n}\sum_{i=1}^nz_{li}u_i\right|}
{\sqrt{\frac{1}{n}\sum_{i=1}^nz_{li}^2u_i^2}}\right),$$ conditional on $z_i$
for $i=1,\hdots,n$ can be obtained numerically, for example using a
Monte-Carlo method, and $r$ is adjusted accordingly.  This is an
approach proposed by Belloni and Chernozhukov (2010) for the linear
regression model without endogeneity.  This approach is worth
mentioning because it has the advantage to yield a smaller $r$, and thus
smaller confidence sets as we will see later.  Indeed, it does
not rely on a crude union bound but account for the
correlation between the regressors.

\noindent{\bf Scenario 2.} Under the assumption
\begin{assumption}\label{ass:distr1}
For every $i=1,\hdots,n$ and $l=1,\hdots,L$, $z_{li}u_i$ are
symmetric and neither of $z_{li}u_i$ is almost surely equal to 0.
\end{assumption}
we choose
\begin{equation}\label{er1}
r=\sqrt{\frac{2\log(L/2\alpha)}{n}}.
\end{equation}
Symmetry is a very likely assumption if \eqref{estruct} is a first
difference between two time periods in a linear panel data model.

\noindent{\bf Scenario 3.} Under the assumption
\begin{assumption}\label{ass:distr2}
For every $l=1,\hdots,L$, $z_{li}u_i$ are i.i.d. and symmetric,
neither of $z_{li}u_i$ is almost surely equal to 0 and $L$ is such
that
$$L<\frac{9\alpha}{4e^3\Phi(-\sqrt{n})}.$$
\end{assumption}
we choose
\begin{equation}\label{er2}
r=-\frac{1}{\sqrt{n}}\Phi^{-1}\left(\frac{9\alpha}{4Le^3}\right).
\end{equation}
This is a slightly tighter constant that the one given for Scenario
2 but Assumption \eqref{ass:distr2} does not allow for
heteroscedastic errors.  The upper bound on the number of
instruments is of the order of an exponential in $n$.

\noindent{\bf Scenario 4.} We denote by
$d_{n,\delta}=\min_{l=1,\hdots,L}d_{n,\delta,l}$ where
$$d_{n,\delta,l}\triangleq\frac{(\sum_{i=1}^n
\E[z_{li}^2u_i^2])^{1/2}}
{\left(\sum_{i=1}^n\E[|z_{li}u_i|^{2+\delta}]\right)^{1/(2+\delta)}}
$$
when the numerator and denominator are well defined.   Under the
assumption
\begin{assumption}\label{ass1} There exists $\delta$ positive such that,
for all $i=1,\hdots,n,\ l=1,\hdots,L$,
$
\E[\left|z_{li} u_i\right|^{2+\delta}]<\infty$,
neither of $z_{li}u_i$ is almost surely equal to 0 and $L$ is
such that
$$L\le\frac{\alpha}{2\Phi(-d_{n,\delta})\left(1+A_0
\left(1+d_{n,\delta}^{-1}\right)^{2+\delta}\right)}$$
where $A_0>0$ is the absolute constant of Theorem \ref{th:jing}
in Section \ref{s90}.
\end{assumption}
we choose
\begin{equation}\label{er3}
r=-\frac{1}{\sqrt{n}}\Phi^{-1}\left(\frac{\alpha}{2L}\right).
\end{equation}
Note that if, for any fixed  $l$ in $\{1,\hdots,L\}$, the variables $z_{li}u_i$
are i.i.d., then for $l=1,\hdots,L$, $d_{n,\delta,l}=n^{\frac{\delta}{4+2\delta}}
\frac{(\E[z_{l1}^2u_1^2])^{1/2}}{(\E[|z_{l1}u_1|^{2+\delta}])^{1/(2+\delta)}}$ which
tends to infinity with $n$.
%Thus, in the i.i.d. setting, we can still consider a number of instruments that
%is exponential in $n$.
With this choice of $r$ the confidence sets will only be
asymptotically valid\footnote{We only consider asymptotically valid
confidence sets because the moderate deviations result for the
self-normalized sums that we use depends on the constant $A_0$ which is
universal but not explicit.}, in the asymptotic where $L$ and $n$
increase in such a way that
$\Phi^{-1}\left(\frac{\alpha}{2L}\right)d_{n,\delta}^{-1}\to 0$.
This scenario allows for heteroscedasticity and non-symmetric
errors.  This is very desirable in an {\em IV} setup.

\noindent{\bf Scenario 5.} We denote by
$\gamma_4\triangleq\max_{l=1,\hdots,L}\gamma_{4l}$ where
$\gamma_{4l}=\mathbb{E}[(z_{li}u_i)^4]/(\mathbb{E}[(z_{li}u_i)^2])^2$.
Under the assumption
\begin{assumption}\label{ass:distr3}
For every $l=1,\hdots,L$, $z_{li}u_i$ are i.i.d.,
$\mathbb{E}[(z_{li}u_i)^4]<\infty$ and $z_{li}u_i$ is not 0 almost
surely, for some $c_4$ positive $\gamma_4\le c_4$ and
\begin{equation}\label{econdas}
L<\frac{\alpha}{2e+1}\exp\left(\frac{n}{c_4}\right).
\end{equation}
\end{assumption}
we choose
\begin{equation}\label{er4}
r=\sqrt{\frac{2\log(L(2e+1)/\alpha)}{n-c_4\log(L(2e+1)/\alpha)}}.
\end{equation}
It is reasonable to assume that $n-c_4\log(L(2e+1)/\alpha)\ge n/2$
as soon as $n$ is relatively large relative compared to $\log(L/\alpha)$. In
that case, we can take
\begin{equation}\label{er4b}
r=2\sqrt{\frac{\log(L(2e+1)/\alpha)}{n}}
\end{equation}
and have confidence sets with valid finite sample properties.
It is also possible to proceed in two stages.  We first start by
choosing $r$ as in \eqref{er4} with a very rough upper bound on
$\gamma_4$, for example using \eqref{er4b}.  As we will see, it
yields a point estimate and a confidence region for $\beta$, which
in turns yields a point estimate and a confidence set for
$\gamma_4$.  In a second stage, it is possible to plug-in either the
point estimate for $\gamma_4$ or the upper bond from the confidence
set for $\gamma_4$.  This two-stage approach yields
approximately valid confidence sets.  A similar approach has
been used in the simulations in Bertail, Gauth\'erat and
Harari-Kermadec (2005) to obtain finite sample confidence sets for
inference based on empirical $\varphi^*-$discrepencies and
quasi-empirical likelihood methods\footnote{Self-normalized sums
also appear naturally in empirical likelihood contexts.}.

Scenarios 2-5 rely on moderate deviations for self-normalized sums
that are recalled in Section \ref{s90} and are respectively from
Efron (1969), Pinelis (1994), Jing, Shao and Wang (2003) and
Bertail, Gauth\'erat and Harari-Kermadec (2009)\footnote{Bertail,
Gauth\'erat and Harari-Kermadec (2009) provides an upper bound of
tail probabilities of self-normalized sums with an explicit constant
unlike Jing, Shao and Wang (2003).  The result of Jing, Shao and
Wang (2003) is mostly useful to study large deviations, that is why
under Scenario 4 we only consider asymptotically valid confidence
sets.}.  Each moderate deviations result relies on a different set
of assumptions and our confidence sets with finite sample
validity will rely on a restricted class of distributions for the
data generating process.  This is related to the Bahadur and Savage
(1956) impossibility result, see also Romano and Wolf
(2000)\footnote{As explained in Romano and Wolf (2000), the Bahadur
and Savage result implies that we need to make some restriction on
the set of distributions $\mathcal{P}$ on the line in order to
construct conservative confidence intervals for the mean of
$\mathbb{P}$ in $\mathcal{P}$ that are bounded.  A conservative
interval of coverage level $1-\alpha$ is such that the parameter is
contained in the interval with probability at least $1-\alpha$ for
all $\mathbb{P}$ in $\mathcal{P}$ and sample size $n$.}.

\subsection{Confidence Sets}\label{s52}

\begin{theorem}\label{t1}
For every $\beta$ in $\mathcal{B}_s$, under one of the scenarios 2-5
of Section \ref{s51}, together with its respective 
choice of $r$, with probability at least $1-\alpha$ (approximately
at least $1-\alpha$ for scenario 4 and the two-stage procedure with
scenario 5), for any $c$ in $(0,1)$  and any set of indices 
$I$ containing the index of the instrument which is unity,
for any solution $(\widehat\beta^{(c,I)},\widehat\sigma^{(c,I)})$ of the
minimization problem (\ref{IVS}) we have
\begin{equation}\label{eq:t1:1}
\left|\left({\bf D_X}^{-1} (\widehat{\beta}^{(c,I)}-\beta)\right)_{J_0}\right|_p \le
\frac{2\widehat{\sigma}^{(c,I)}r}{\kappa_{p,J_0,J(\beta)}^{(c,I)}}\left(1-\frac{r}{\kappa_{1,J(\beta)}^{(c,I)}}
\right)_+^{-1}\, \quad \forall \ p\in[1,\infty],\ \forall J_0\subset\{1,\hdots,K\},
\end{equation}
\begin{equation}\label{eq:t1:3}
|\widehat{\beta}_k^{(c,I)}-\beta_k| \le
\frac{2\widehat{\sigma}^{(c,I)}r}{\mathbb{E}_n[X_k^2]^{1/2}\, \kappa^{(c,I)*}_{k,J(\beta)}
}\left(1-\frac{r}{\kappa_{1,J(\beta)}^{(c,I)}}\right)_+^{-1}\, \quad \forall k=1,\dots,K,
\end{equation}
and
\begin{equation}\label{eq:t1:4}
\widehat{\sigma}^{(c,I)}\le\max_{l\in I}({\bf D_Z^{(I)}})_{ll}
\sqrt{\widehat{Q}_l(\beta)}\left(1+\frac{r}{c\kappa_{1,J(\beta),J(\beta)}^{(c,I)}}\right)
\left(1-\frac{r}{c\kappa_{1,J(\beta),J(\beta)}^{(c,I)}}\right)_+^{-1}.
\end{equation}
Under scenario 1 of Section \ref{s51} with the corresponding
choice of $r$, for any set of indices
$I$ containing the index of the instrument which is unity, 
for every $\beta$ in $\mathcal{B}_s$, 
with probability at least $1-\alpha$, 
for any $c$ in $(0,1)$,
any solution $(\widehat\beta^{(c,I)},\widehat\sigma^{(c,I)})$  of the
minimization problem (\ref{IVS}) satisfies 
\eqref{eq:t1:1}, \eqref{eq:t1:3} and \eqref{eq:t1:4}.
\end{theorem}
%The proof of Theorem \ref{t1} is given in Section \ref{s93}.

Consider the most likely case where $L\ge K$ and the parameter $\beta$ is point identified.
One can then consider Theorem \ref{t1} with $s=K$ and $\mathcal{B}_s=\mathcal{I}dent=\{\beta^*\}$.
Still, in the high-dimensional
framework, we have that ${\rm rank}\ \Psi_n^{(I)}\le \min(n,K)=n<K$,
thus ${\rm dim}({\rm Ker}\ \Psi_n^{(I)})\ge K-n>0$.  Hence, in the high-dimensional
framework,
$\Psi_n^{(I)}\Delta=0$ can have non-zero solutions.
It is because of the cone constraint
that the sensitivities $\kappa_{p,J_0,J(\beta^*)}^{(c,I)}$
and $\kappa^{(c,I)*}_{k,J(\beta^*)}$ on the right hand-side of
\eqref{eq:t1:1} and \eqref{eq:t1:3} can be different from 0
and thus yield non-trivial upper bounds.

Consider now the case of Example 6 where the instruments can have a direct
effect on the outcome.  In that case $L<K$ and identification fails without
further assumptions.  This is simply because
${\rm dim\ Ker}\left(\mathbb{E}[z_ix_i^T]\right)\ge K-L>0$.
If we knew the number $s$ of non-zero coefficients, the underlying
vector $\beta$ is the solution of one system
\begin{equation}\label{eqidentEx6}
\left\{\begin{array}{c}
            \beta_{J^c}=0\\
           \mathbb{E}[z_i(y_i-x_{Ji}^T\beta_{J})]=0
         \end{array}
\right.
\end{equation}
for each subset $J$ of $\{1,\hdots,K\}$ of size $s$
(including the index of the constant regressor).
For every set $J$ such that $\mathbb{E}[z_ix_{Ji}^T]$
has full column rank, \eqref{eqidentEx6} has at most one solution.
When $s<L$, $\mathbb{E}[z_ix_{Ji}^T]\beta_{J}=\mathbb{E}[z_iy_i]$
does not necessarily have a solution because the system is overdetermined.
Note as well that, because $s$ is the number
of non-zero coefficients, in addition to \eqref{eqidentEx6}, we know
that for every $k$ in $\{1,\hdots,K\}$, $\beta_{k}\ne 0$.  This again
restricts the set of solutions.
Note that there are also cases where the set of solutions of
\eqref{eqidentEx6} is an affine space.
For example, if $s>L$, $\mathbb{E}[z_ix_{Ji}^T]$
does not have full column rank.
Geometrically $\mathcal{B}_s$ is the union of the intersections of $\mathcal{I}dent$
with the spaces $\{\forall k\in J,\ \beta_k=0\}$ if for every $k$ in $\{1,\hdots,K\}$, $\beta_{k}\ne 0$.
This is a union of affine spaces.
We have observed that, adding the prior information
that there are at most $s<L$ non-zero coefficients, is very likely to
yield a set $\mathcal{B}_s$ which is a union of points.  Recall that,
in the context of Example 6, $s<L$ means that there are in reality instruments
that do not have a direct effect on the outcome.
There are thus unknown exclusion restrictions which have an identification power.
Again, because of the cone constraints, the sensitivities $\kappa_{p,J_0,J(\beta^*)}^{(c,I)}$
and $\kappa^{(c,I)*}_{k,J(\beta^*)}$ on the right hand-side of
\eqref{eq:t1:1} and \eqref{eq:t1:3} can be different from 0 even in the high-dimensional framework
where $n$ is small relative $K$.

Consider now the term $\left(1-r/\kappa_{1,J(\beta)}^{(c,I)}\right)_+^{-1}$
in the upper bounds of Theorem \ref{t1}.
Observe that the confidence sets of level at least $1-\alpha$
for $\beta$ in $\mathcal{B}_s$ can have infinite volume on the
random event $\frac{r}{\kappa_{1,J(\beta)}^{(c,I)}}>1$.  Indeed,
the upper bounds (\ref{eq:t1:1}) and (\ref{eq:t1:3}) are infinite.
This occurs either when $r$ is
large or when $\kappa_{1,J(\beta)}^{(c,I)}$ is too small.  Recall that $r$
is of the order of $\sqrt{\log(L)/n}$ in all scenarios.  Because
increasing $L$ increases the sensitivities, once
$\kappa_{1,J(\beta)}^{(c,I)}$ is sufficiently bounded away from zero, the
first condition corresponds to a situation with very many
instruments relative to the sample size.  This occurs when $L$ is as
large as an exponential in $n$.  To interpret the second condition,
note that Proposition \ref{p5} (ii) yields that for any $k$ in
$J(\beta)$,
$$\kappa_{1,J(\beta)}^{(c,I)}\le\kappa_{k,J(\beta)}^{(c,I)*}.$$
This can occur on the event where a regressor
$(x_{ki})_{i=1}^n$ has a restricted maximal partial empirical
correlation with the instruments less than $r$.
In the setting of Example 7, this becomes increasingly likely when
the parameter $\zeta_{J^c}$ in the data generating
approaches zero (see \eqref{eInstrStrength}).
Confidence sets of infinite volume with positive probability
is a desired feature of a procedure that is
robust to weak instruments (see Dufour (1997) building on Gleser and Hwang (1987)).
We will show in Section
\ref{sLDWI} that for low dimensional models ($K$ is small and there
is no sparsity) a slight modification of the {\em STIV} estimator is
a new procedure that is robust to weak instruments
(see Andrews and Stock (2007) for a review of existing methods).
Note also that the term $\left(1-r/\kappa_{1,J(\beta)}^{(c,I)}\right)_+^{-1}$
appears because we have a procedure that relaxes distributional assumptions
and, in the homoscedastic case, does not require to know the variance of the errors
or an upper bound.  When the regressors are deterministic and the errors are Gaussian,
an upper bound on the variance is sufficient to adjust the penalization and obtain a
consistent Lasso or Dantzig estimator.
In the linear regression with exogenous regressors,
$\mathbb{E}[y_i^2]$ provides such an upper bound on the variance.
This no longer works in the presence of endogenous regressors.

When
$\tau_1\triangleq 1-\frac{r}{\kappa_{1,J(\beta)}^{(c,I)}}$ is close to 1 and
the sensitivities $\kappa_{p,J_0,J(\beta)}^{(c,I)}$ and
$\kappa^{(c,I)*}_{k,J(\beta)}$ are bounded away from zero, the upper bounds
in \eqref{eq:t1:3} is of the order $O(r)=O(\sqrt{\log(L)/n})$. Thus,
we have an extra $\sqrt{\log(L)}$ factor as compared to the usual
root-$n$ rate.  It is a modest price for using a large number $L$ of
instruments.  Under the premise of Proposition \ref{p5}, for
$\tau_1\approx 1$ it is sufficient to have $|J(\beta)|\le
Cr^{-1}=O(\sqrt{n/\log(L)})$ where $C>0$ is a proper constant.  This
is quite a reasonable condition on the sparsity $|J(\beta)|$.

%Note that a condition somewhat more stringent than positivity of
%$1-\frac{r^2}{\kappa_{1,J(\beta^*)}}$ is required in Belloni,
%Chernozhukov and Wang (2010) to obtain a finite sample bound on the
%prediction loss of the Square-root Lasso. This agrees with our value
%of $\tau_1$ when $J_{{\rm end}}=\varnothing$ and thus
%$\frac{r}{\kappa_{J_{{\rm end}},J(\beta^*)}^*}=0$ (indeed, in
%Belloni, Chernozhukov and Wang (2010) all the variables are
%exogenous).

The only unknown ingredient of the inequalities \eqref{eq:t1:1} and
(\ref{eq:t1:3}) is the set $J(\beta)$ that determines the
sensitivities.  To turn these inequalities into valid confidence
bounds, it suffices to provide data-driven lower estimates on the
sensitivities.  As discussed in Section \ref{s4}, there are two ways
to do it.  The first one is based on the sparsity certificate, {\it
i.e.}, assuming some known upper bound $s$ on $|J(\beta)|$; then we
get bounds depending only on $s$ and on the data.
\begin{corollary}\label{cor1}
For every $\beta$ in $\mathcal{B}_s$, under scenarios 2-5
of Section \ref{s51} together with its respective assumption and
choice of $r$, with probability at least $1-\alpha$ (approximately
at least $1-\alpha$ for scenario 4 and the two-stage procedure with
scenario 5), for any $c$ in $(0,1)$  and any set of indices $I$ containing the index of the instrument which is unity,
for any solution $(\widehat\beta^{(c,I)},\widehat\sigma^{(c,I)})$ of the
minimization problem (\ref{IVS}) we have
\begin{equation}\label{eq:CIlp}
\left|\left({\bf D_X}^{-1} (\widehat{\beta}^{(c,I)}-\beta)\right)_{J_0}\right|_p \le
\frac{2\widehat{\sigma}^{(c,I)}r}{\overline{\kappa}_{p,J_0}^{(c,I)}(s)}\left(1-\frac{r}{\kappa_{1}^{(c,I)}(s)}
\right)_+^{-1}\, \quad \forall \ p\in[1,\infty],\ \forall J_0\subset\{1,\hdots,K\},
\end{equation}
where $\overline{\kappa}_{p,J_0}^{(c,I)}(s)$ is any lower bound on
$\kappa_{p,J_0, J(\beta)}^{(c,I)}$ based on the sparsity certificates that
is convenient to calculate (see, {\it e.g.}, \eqref{eq:certif2} and
\eqref{eq:certifgen})
\begin{equation}\label{eq:CI}
|\widehat{\beta}_k^{(c,I)}-\beta_k| \le
\frac{2\widehat{\sigma}^{(c,I)}r}{\mathbb{E}_n[X_k^2]^{1/2}\, \kappa^{(c,I)*}_{k}(s)
}\left(1-\frac{r}{\kappa_{1}^{(c,I)}(s)}\right)_+^{-1}\, \quad \forall k=1,\dots,K,.
\end{equation}
Under scenario 1 of Section \ref{s51} with the corresponding
choice of $r$, for any set of indices
$I$ containing the index of the instrument which is unity,
for every $\beta$ in $\mathcal{B}_s$,
with probability at least $1-\alpha$,
for any $c$ in $(0,1)$,
any solution $(\widehat\beta^{(c,I)},\widehat\sigma^{(c,I)})$ of the
minimization problem (\ref{IVS}) satisfies
\eqref{eq:CIlp} and \eqref{eq:CI}.
\end{corollary}
In simulations the {\em STIV} estimator has always more non-zeros
than the truth so it is reasonable to take $s=|J(\widehat{\beta}^{(c,I)})|$
or some larger value.  We also advertise the possibility of drawing
nested confidence sets for increasing values of $s$.

The second way is to plug in, instead of $J(\beta)$, some
data-driven upper estimate $\widehat J$, {\it i.e.}, a set
satisfying $J(\beta)\subseteq \widehat J$ on the intersection of the
event of probability at least $1-\alpha$ (or approximately
$1-\alpha$ for scenario 4 and the two-stage procedure with scenario
5) of Theorem \ref{t1} and an event of probability close to~1.
Theorem \ref{t1b} provides examples of such estimators $\widehat J$.
In particular, assertion (iv) of Theorem \ref{t1b} guarantees that,
under some assumptions, the estimator $\widehat J= J(\widehat
\beta^{(c,I)})$ has the required property.  The statement of the
corresponding result is postponed to Section \ref{s54}.

We do not touch upon optimality of the confidence regions in this
article.  Because of the popularity of two-stage least squares, we
present the properties of a high dimensional two-stage method akin
to two-stage least squares in Section \ref{sLPI}.  The {\em STIV}
estimator depends on the tuning parameter $c$.
A smaller $c$ implies a smaller cone
$C_{J(\beta)}^{(c,I)
}$ and thus a larger sensitivity.  On the other hand,
because we penalize less $\sigma$ in \eqref{IVS}, $\widehat{\sigma}^{(c,I)}$
is larger.  Because the dependence is nonlinear and the
sensitivities are random and depend on the distribution of the data
generating process, there does not exist a universally good value
for $c$\footnote{The same is true for the Square-root Lasso.}.
It is important to note that our results hold on the $1-\alpha$ probability event,
for any value of $c$ and any set of indices $I$ containing the index of the instrument
which is unity (for scenarios 2-5 only for the set $I$).  
It is therefore possible to consider values of $c$ or $I$ which
are random and depend on the data.
If one is interested in a specific coefficient,
it is possible to pick the values of $c$ and $I$ that yields the smaller confidence interval
for that specific regressor.  Because
the procedure is fast to implement, varying $c$ on a
grid on $(0,1)$ and picking the value that yields the smaller confidence
interval is an easy thing to do.  Varying $I$ is too complicated but,
if one hesitates on a few specifications of $I$, it is possible to build confidence sets
for these various possibilities and choose the set $I$ that yields the smaller confidence sets.
%For example there could be (rare) situations
%where changing $c$ would allow to obtain confidence regions with
%finite volume.

\subsection{Sparse Oracle Inequality}\label{s53}
We now consider the approximately sparse setting.  The sparsity
assumption is quite natural in empirical economics since usually
only a moderate number of covariates is included in the model.
However, one might be also interested in the case when $\beta$ is
only approximately sparse.  This means that most of the coefficients
$\beta$ are not exactly zero but too small to matter, whereas the
remaining ones are relatively large. This setting received some
attention in the statistical literature. For example, the
performance of Dantzig selector and $MU$-selector under such
assumptions is studied by Cand{\`e}s and Tao (2007) and Rosenbaum
and Tsybakov (2010) respectively. We will derive a similar result
for the {\em STIV} estimator.

Consider the enlarged cone
$$\widetilde{C}_{J}^{(c)}\triangleq\left\{\Delta\in\R^K:\ |\Delta_{J^c}|_1
\le\frac{2+c}{1-c}|\Delta_{J}|_1\right\}$$
and define, for $p\in[1,\infty]$ and $J_0\subset\{1,\hdots,K\}$
$$\widetilde{\kappa}_{p,J_0,J}^{(c,I)
}\triangleq
        \inf_{\Delta\in\R^K:\ |\Delta_{J_0}|_p=1,\ \Delta\in\widetilde
        {C}_{J}^{(c)}}\left|\Psi_n^{(I)}\Delta
\right|_{\infty}$$ and $\widetilde{\kappa}_{1,J}^{(c,I)
}$ corresponds to
$\widetilde{\kappa}_{p,J_0,J}^{(c,I)}$ with $p=1$ and $J_0=\{1,\hdots,K\}$.

The following theorem is an analog of the above results for the
approximately sparse case.

%%%%%%%%%%%%
\begin{theorem}\label{tapproxsparse}
For every $\beta$ in $\mathcal{B}_K$, under scenarios 2-5
of Section \ref{s51} together with its respective assumption and
choice of $r$,  with probability at least $1-\alpha$ (approximately
at least $1-\alpha$ for scenario 4 and the two-stage procedure with
scenario 5), for any $c$ in $(0,1)$  and any set of indices $I$ containing
the index of the instrument which is unity, for any solution
$(\widehat\beta^{(c,I)},\widehat\sigma^{(c,I)})$ of the
minimization problem (\ref{IVS}) we have, for every
$J_0\subset\{1,\hdots,K\}$,
\begin{equation}\label{eq:tapproxsparse}
\left|\left({\bf D_X}^{-1} \left(\widehat{\beta}^{(c,I)
}-\beta\right)\right)_{J_0}\right|_p
\le \min_{J\subset\{1,\hdots,K\}} \left\{\max\left(
\frac{2\widehat{\sigma}^{(c,I)} r }{\widetilde{\kappa}_{p,J_0,J}^{(c,I)}}
\left(1-\frac{r}{\tilde{\kappa}_{1,J}^{(c,I)}}\right)_+^{-1},
%\frac{2\sqrt{2}\sigma
%r}{(\widetilde{\kappa}_{p,J}-3rc^{-1}(1-c)^{-1} |J|^{1-1/p} )_+}, \
\ \frac{6 \left|\left({\bf
D_X}^{-1}\beta\right)_{J^c}\right|_1}{1-c}\right)\right\} \,.
\end{equation}
Under scenario 1 of Section \ref{s51} with the corresponding
choice of $r$, for any set of indices
$I$ containing the index of the instrument which is unity,
for every $\beta$ in $\mathcal{B}_s$,
with probability at least $1-\alpha$,
for any $c$ in $(0,1)$,
any solution $(\widehat\beta^{(c,I)},\widehat\sigma^{(c,I)})$ of the
minimization problem (\ref{IVS}) satisfies 
\eqref{eq:tapproxsparse} for every $J_0\subset\{1,\hdots,K\}$.
\end{theorem}

We can interpret Theorem \ref{tapproxsparse} as the fact that the
{\em STIV} estimator automatically realizes a ``bias/variance"
trade-off related to a non-linear approximation.  Inequality
(\ref{eq:tapproxsparse}) means that this estimator performs as well
as if the optimal subset $J$ were known.

\section{Rates of Convergence and Selection of the Variables}\label{s54}
Let us consider rates of convergence of our estimator.  We need to
replace the random right hand-side of \eqref{eq:t1:1},
\eqref{eq:t1:3} and \eqref{eq:t1:4} by deterministic upper bounds.
In this section we consider that $c$ in $(0,1)$ and the set $I$ are fixed.
\begin{assumption}\label{ass1b}
For every $\beta\in\mathcal{B}_s$ and $\gamma_1\in(0,1)$, there
exists a constant $\sigma_*>0$ such that
$$\mathbb{P}\left(\max_{l\in I}({\bf D_Z^{(I)}})_{ll}^2\E_n[Z_l^2(Y-X^T\beta)^2]\le \sigma_*^2\right)\ge 1-\gamma_1.$$
\end{assumption}

The second assumption concerns the population counterparts of the
sensitivities.
%It is stated in terms of subsets $J_0$ of $\{1,\dots,
%K\}$ and constants $p\ge 1,$ $k\in \{1,\dots, K\} $ that can differ
%from case to case and will be specified later.

\begin{assumption}\label{ass1c}
For every $\beta\in\mathcal{B}_s$ and $\gamma_2\in(0,1)$, there
exists constants $c_{p}^{(c,I)}>0$, $c_{1,J_0}^{(c,I)}>0$ for $J_0=J(\beta)$ and
$J_0=\{k\}$ for $k=1,\hdots,K$, such that, with probability at least
$1-\gamma_2$,
\begin{eqnarray}
\label{Assa1}
\kappa_{p,J(\beta)}^{(c,I)} &\ge& c_{p}^{(c,I)
}|J(\beta)|^{-1/p},\\
\label{Assa2} \kappa_{1,J_0,J(\beta)}^{(c,I)} &\ge& c_{1,J_0}^{(c,I)}.
\end{eqnarray}
\end{assumption}
If $J_0=\{k\}$ is a singleton we write for brevity
$c_{1,J_0}^{(c,I)}=c_{k}^{(c,I)*}$.

The dependence on $|J(\beta)|$ of the right hand-side of
(\ref{Assa1}) is
motivated by Proposition \ref{p5}. In (\ref{Assa2}), %and
%(\ref{Assa3})
we do not indicate the dependence of the bounds on $|J(\beta)|$
explicitly because it can be different for different
sets $J_0$. %{\bf ??? If $J_0=\{k\}$ the constant $c_{J_0}^*$, which
%we denote in this case by $c_{k}^*$, can be taken independent on
%$|J(\beta^*)|$.}
For general $J_0$, combining Proposition \ref{p4}~(ii) and
Proposition \ref{p5} suggests that the value $c_{1,J_0}^{(c,I)}$ can be
bounded from below by a quantity of the order $|J(\beta)|^{-1}$.
Note, however, that this is a coarse bound valid for any set~$J_0$.

The last assumption defines a population counterpart of
$\mathbb{E}_n[X_k^2]^{1/2}$.

\noindent\begin{assumption}\label{ass1X} For every $\gamma_3\in(0,1)$ and
$k\in\{1,\hdots,K\}$, there exist constants $v_k>0$ such that
$$\mathbb{P}\left(\mathbb{E}_n[X_k^2]^{1/2}\ge v_k, \ \forall
k\in \{1,\hdots,K\}\right)\ge 1-\gamma_3.$$
\end{assumption}

Assumptions \ref{ass1b}, \ref{ass1c} and \ref{ass1X} are very weak.
They are required to obtain rates of convergence, {\it i.e.}, a
deterministic bound on the estimation error on a high probability
event (see, {\it e.g.}, \eqref{eq:th1b:3} below).  We set
$\gamma=\alpha+\sum_{j=1}^3\gamma_j$, and
$$\tau^{(c,I)*}\triangleq \left(1+\frac{r}{c c_{1,J(\beta)}^{(c,I)}}\right)
\left(1-\frac{r}{c
c_{1,J(\beta)}^{(c,I)}}\right)_+^{-1}\left(1-\frac{r|J(\beta)|}{c_{1}^{(c,I)}}
\right)_+^{-1}\,.$$

\begin{theorem}\label{t1b}
For every $\beta$ in $\mathcal{B}_s$, under the assumptions of
Theorem~\ref{t1} and Assumption \ref{ass1b}, the following holds.
\begin{itemize}
\item[(i)]  Let part (\ref{Assa2}) of Assumption \ref{ass1c} be satisfied.
Then, with probability at least $1-\alpha-\gamma_1-\gamma_2$,
for any solution $\widehat\sigma^{(c,I)}$ of (\ref{IVS}) we have
$$
\widehat{\sigma}^{(c,I)}\le\sigma_* \left(1+\frac{r}{c
c_{1,J(\beta)}^{(c,I)}}\right) \left(1-\frac{r}{c
c_{1,J(\beta)}^{(c,I)}}\right)_+^{-1}.
$$
\item[(ii)]  Fix $p\in[1,\infty]$.  Let
Assumption \ref{ass1c} be satisfied.  Then, with probability at
least $1-\alpha-\gamma_1-\gamma_2$, for any solution
$\widehat\beta^{(c,I)}$ of (\ref{IVS}) we have
\begin{equation}\label{eq:th1b:1}
\left|{\bf D_X}^{-1} \left(\widehat{\beta}^{(c,I)}-\beta\right)\right|_p
\le \frac{2\sigma_* r |J(\beta)|^{1/p}\tau^{(c,I)}*}{c_p^{(c,I)}}\,.
\end{equation}
\item[(iii)] Let Assumptions \ref{ass1c} and \ref{ass1X}
be satisfied.  Then with probability at least $1-\gamma$, for
any solution $\widehat\beta^{(c,I)}$ of (\ref{IVS}) we have
\begin{equation}\label{eq:th1b:3}
|\widehat{\beta}_k^{(c,I)}-\beta_k| \le \frac{2\sigma_* r
\tau^{(c,I)*}}{c_{k}^{(c,I)*}v_k}\, , \quad  \ k=1,\dots,K.
\end{equation}
\item[(iv)]  Let the assumptions of {\rm (iii)} hold, and
$|\beta_k|>\frac{2\sigma_* r\tau^{(c,I)*}}{c_{k}^{(c,I)*}v_k}$ for all $k\in
J(\beta)$.  Then, with probability at least $1-\gamma$, for any
solution $\widehat\beta^{(c,I)}$ of (\ref{IVS}) we have
$$J(\beta) \subseteq J(\widehat\beta^{(c,I)}).$$
\end{itemize}
\end{theorem}

%The proof of Theorem \ref{t1b} is given in Section~\ref{s93}.

For reasonably large sample size ($n\gg \log(L)$), the value $r$ is
small, and $\tau^{(c,I)*}$ is approaching 1 as $r\to0$. Thus, the bounds
(\ref{eq:th1b:1}) and (\ref{eq:th1b:3}) are of the order of
magnitude $O(r |J(\beta)|^{1/p})$ and $O(r)$ respectively. These are
the same rates, in terms of the sparsity $|J(\beta)|$, the dimension
$L$, and the sample size $n$, that were proved for the Lasso and
Dantzig selector in high-dimensional regression with Gaussian errors
and without endogenous variables (here $L=K$) in Cand{\`e}s and
Tao~(2007), Bickel, Ritov and Tsybakov~(2009), Lounici~(2008) (see
also B\"uhlmann and van de Geer~(2011) for references to more recent
work).

The lower bound for estimation of a high dimensional linear
regression model, without endogeneity, under fixed and Gaussian
design, with $K$ potential regressor and at most $s$ non-zero
coefficients, is known to be
$\sqrt{\frac{|J(\beta)|\log(K/|J(\beta)|)}{n}}$ (see Verzelen
(2010), Ye and Zhang (2010) and Raskutti, Wainwright and Yu (2011),
rates for prediction are also obtained in Rigollet and Tsybakov
(2011)).  Proving minimax lower bounds for the model
\eqref{estruct}-\eqref{einstr} with $L\ne K$ is subject of future
investigation.

Note that Theorem \ref{t1b} assumes that we work on some event such that the lower bounds $c_p^{(c,I)}$, $c_k^{(c,I)*}$
and $c_{1,J(\beta)}^{(c,I)}$ are positive. When for a particular value of $\beta$ there are 0 on a large probability event
then on that event we do not have convergence rates or model selection results.
Also the extra condition in (iv) is not satisfied for every $\beta$ in $\mathcal{B}_s$.
Therefore when $\mathcal{B}_s$ is not a point then the conclusion of (iv) only holds for certain values of the parameter
in the set $\mathcal{B}_s$.

From (\ref{eq:t1:3}), Theorem~\ref{t1b} (iv) and Proposition \ref{p4} (i), we obtain the
following confidence sets of level $1-\gamma$ for~$\beta_k$ in
$\mathcal{B}_s$.
\begin{corollary}\label{cor2}
For every $\beta$ in $\mathcal{B}_s$, under one of the 5 scenarios
of Section \ref{s51}, together with its respective assumption and
choice of $r$, under the assumptions of Theorem \ref{t1b} (iv), with
probability at least $1-\gamma$ (approximately at least $1-\gamma$
for scenario 4 and the two-stage procedure with scenario 5) for any
solution $(\widehat\beta,\widehat\sigma)$ of the minimization
problem (\ref{IVS}) we have
\begin{equation}\label{eq:CIlpPI}
\left|\left({\bf D_X}^{-1} (\widehat{\beta}^{(c,I)}-\beta)\right)_{J_0}\right|_p \le
\frac{2\widehat{\sigma}^{(c,I)}r}{\overline{\kappa}_{p,J_0,J(\widehat{\beta}^{(c,I)})}^{(c,I)}}
\left(1-\frac{r}{\kappa_{1,J(\widehat{\beta}^{(c,I)})}^{(c,I)}}
\right)_+^{-1}, \quad \forall \ p\in[1,\infty],\ \forall J_0\subset\{1,\hdots,K\},
\end{equation}
where $\overline{\kappa}_{p,J_0,J(\widehat{\beta}^{(c,I)})}^{(c,I)}$ is any lower bound on
$\kappa_{p,J_0,J(\widehat{\beta}^{(c,I)})}^{(c,I)}$ that is convenient to calculate (see, {\it e.g.},
\eqref{eq:certif2}and \eqref{eq:certifgen}), and for all $k=1,\dots,K$,
\begin{equation}\label{eq:conf_int}
|\widehat{\beta}_k^{(c,I)}-\beta_k| \le
\frac{2\widehat{\sigma}^{(c,I)}r}{\mathbb{E}_n[X_k^2]^{1/2}\, \kappa^{(c,I)*}_{k,J(\widehat{\beta}^{(c,I)})}
}\left(1-\frac{r}{\kappa_{1,J(\widehat{\beta}^{(c,I)})}^{(c,I)}}\right)_+^{-1}.
\end{equation}
\end{corollary}
Here we do not control exactly the coverage probability because its
probability is now $1-\gamma$.  Having $\gamma_j$ very small for
$j=1,2,3$ and thus coverage probability of approximately at least
$1-\alpha$ requires taking a large $\sigma^*$ and lower bounds on
the sensitivities.  Such a choice is compatible with (iv) from
Theorem \ref{t1b} only when $\beta_k$ are large enough for $k\in
J(\beta)$.

Theorem~\ref{t1b} (iv) provides an upper estimate on the set of
non-zero components of $\beta$.  We now consider the problem of the
exact selection of variables.  For this purpose, we use the
thresholded {\em STIV} estimator whose coordinates are defined by
\begin{equation}\label{eq:thresh}
\widetilde{\beta}_k^{(c,I)}(\omega_k^{(c,I)})\triangleq\left\{\begin{array}{ll}
                  \widehat{\beta}_k^{(c,I)} & {\rm if\ }|\widehat{\beta}_k^{(c,I)}|>
                  \omega_k^{(c,I)}, \\
                  0 & {\rm otherwise,}
                \end{array}\right.
\end{equation}
where $\widehat{\beta}_k^{(c,I)}$ are the coordinates of the {\em STIV}
estimator $\widehat{\beta}^{(c,I)}$, and $\omega_k^{(c,I)}>0, \ k=1,\hdots,K$, are
thresholds that will be specified below. We will use the sparsity
certificate approach, so that the thresholds will depend on the
upper bound $s$ on the number of non-zero components of $\beta$.
%They will be of the general
%form
%$$
%\omega_{k,J}\triangleq \frac{2\widehat\sigma r }{\kappa_{k,
%J}^*x_{k*}}\left(1-\frac{r}{\kappa^*_{J_{\rm end}, J}
%}-\frac{r^2}{\kappa^*_{J_{\rm end}^c, J}}\right)_+^{-1}\,
%$$
%for appropriate sets $J$.
Because, for practical use, we want to a data driven thresholding
rule we need to strengthen Assumption~\ref{ass1c} as follows.
\begin{assumption}\label{ass_select}
Fix an integer $s\ge 1$. For every $\gamma_2\in(0,1)$, for $\beta$
in $\mathcal{B}_s$, there exist constants $c_{1,J(\beta)}^{(c,I)}>0$,
$c_{1,J_0}^{(c,I)}(s)>0$ for $J_0=\{k\}$, $\forall k$ and
$J_0=\{1,\hdots,K\}$, such that, with probability at least
$1-\gamma_2$,
\begin{eqnarray}
\label{eq:ass_select} \kappa_{1,J(\beta),J(\beta)}^{(c,I)} &\ge&
c_{1,J(\beta)}^{(c,I)}  \quad {\rm and} \quad \kappa_{1,J_0}^{(c,I)}(s) \ge
c_{1,J_0}^{(c,I)}(s)
\end{eqnarray}
\end{assumption}
If $J_0=\{k\}$ is a singleton we write for brevity
$c_{1,J_0}^{(c,I)}(s)=c_{k}^{(c,I)*}(s)$ and $c_{1,\{1,\hdots,K\}}^{(c,I)}(s)=c_{1}^{(c,I)}(s)$.
Set
$$
\tau^{(c,I)*}(s)\triangleq \left(1+\frac{r}{c c_{1,J(\beta)}^{(c,I)}}\right)
\left(1-\frac{r}{c
c_{1,J(\beta)}^{(c,I)}}\right)_+^{-1}\left(1-\frac{r}{c_{1}^{(c,I)}(s)}\right)_+^{-1}\,.
$$

The following theorem shows that, based on thresholding of the {\em
STIV} estimator, we can reconstruct exactly the set of non-zero
coefficients $J(\beta)$ with probability close to~1. Even more, we
achieve the sign consistency, {\em i.e.,} we reconstruct exactly the
vector of signs of the coefficients of $\beta$ with probability
close to~1.

\begin{theorem}\label{th:threshold}
For every $\beta$ in $\mathcal{B}_s$, let the assumptions of
Theorem~\ref{t1} and
Assumptions~\ref{ass1b},~\ref{ass1X},~\ref{ass_select} be satisfied.
Assume that $|J(\beta)|\le s$, and $|\beta_k|>\frac{4\sigma_*
r\tau^{(c,I)*}(s)}{c_{k}^{(c,I)*}(s)v_k}$ for all $k\in J(\beta)$.  Take the
thresholds
$$
\omega_{k}^{(c,I)}(s)\triangleq \frac{2\widehat\sigma^{(c,I)} r
}{\kappa_{k}^{(c,I)*}(s)\mathbb{E}_n[X_k^2]^{1/2}}\left(1-\frac{r}{\kappa_{1}^{(c,I)}(s)}\right)_+^{-1},
$$
and consider the estimator $\widetilde{\beta}^{(c,I)}$ with coordinates
$\widetilde{\beta}_k^{(c,I)}(\omega_{k}^{(c,I)}(s)), k=1,\dots, K$. Then, with
probability at least $1-\gamma$, we have
\begin{equation}\label{eSelb}
\overrightarrow{{\rm sign} (\widetilde{\beta}^{(c,I)})}=\overrightarrow{{\rm
sign}(\beta)}.
\end{equation}
As a consequence, $ J(\widetilde{\beta}^{(c,I)})=J(\beta)$.
\end{theorem}
Sign consistency implies perfect model selection but it is also of
independent interest.  Indeed, the above result gives conditions
under which one can make valid statements on the sign of the
coefficients.  Unfortunately this is only a theoretical result
because it is not possible to verify the condition that
$|\beta_k|>\frac{4\sigma_* r\tau^{(c,I)*}(s)}{c_{k}^{(c,I)*}(s)v_k} $ for all
$k\in J(\beta)$.

%The proof of Theorem~\ref{th:threshold} is given in
%Section~\ref{s93}.
%%%%%%%%%%%%%%%%%%%%%%%%%%%%%%%

\section{Further Results on the STIV estimator}\label{s7}
\subsection{Low Dimensional Models and Very-Many-Weak-Instruments}\label{sLDWI}
In this section we consider a slight modification of the {\em STIV}
estimator to deal with the case where $K$ is small relative to $n$
and we know that all $K$ coefficients are nonzero. It is a new
procedure that is robust to a situation where each individual
instrument is weak. It can handle very-many instruments, {\it i.e.}
$L$ can be much larger than $n$. It does not rely on normal errors and the confidence
sets under Scenarios 2 and 4 (asymptotic) allow for
heteroscedasticity. Again, we do not rely on non-standard
asymptotics.

Because this is a very active and challenging area of research in
econometrics we cannot make an exhaustive literature review and
refer to Andrews and Stock (2007) and the references therein.
%Related work on weak-instruments and many-weak-instruments include, for example,
%Staiger and Stock (1997), Kleibergen (2002, 2005),

Because $K$ is small and $L$ is large \eqref{estruct}-\eqref{einstr}
is usually point identified, but in the presence of weak instruments
one is at the verge of identification. Thus for the sake of
completeness we do not exclude partial identification from our
result.

\begin{definition}
We call the {\em STIV-R} estimator any solution
$(\widehat{\beta}^{(I)},\widehat{\sigma}^{(I)})$ of the following minimization
problem:
\begin{equation}\label{IVSRWI}
\min_{(\beta,\sigma)\in \widehat{\mathcal{I}}^{(I)}} \sigma .
\end{equation}
\end{definition}
In that case there are no cone constraints in the definition of
the sensitivities and we simply drop the index $J(\beta)$ in the
notation of the sensitivities.  We also drop the exponent $c$ everywhere because
the minimization problem \eqref{IVSRWI} no longer involves the constant
$c$ because the model is low dimensional.
Unlike the high-dimensional setup
of the previous sections, the sensitivities can be directly
calculated from the data and we do not have to rely on lower bounds. 
Obtaining confidence sets in a non high-dimensional structural equation
is much more direct and easy. We obtain the following result.
\begin{theorem}\label{tRWI}
For every $\beta$ in $\mathcal{I}dent$, under one of the 5 scenarios
of Section \ref{s51} together with its respective assumption and
choice of $r$, with probability at least $1-\alpha$ (approximately
at least $1-\alpha$ for scenario 4 and the two-stage procedure with
scenario 5) for any solution $(\widehat\beta^{(I)},\widehat\sigma^{(I)})$ of the
minimization problem (\ref{IVSRWI})
\begin{equation}\label{eq:t1RWI:1}
\left|\left({\bf D_X}^{-1} (\widehat{\beta}^{(I)}-\beta)\right)_{J_0}\right|_p \le
\frac{2\widehat{\sigma}^{(I)}r}{\kappa_{p,J_0}^{(I)}}\left(1-\frac{r}{\kappa_{1}^{(I)}}\right)_+^{-1}\,,
\quad \forall \ p\in[1,\infty],\ \forall J_0\subset \{1,\hdots,K\},
\end{equation}
for all $k=1,\dots,K,$
\begin{equation}\label{eq:t1RWI:3}
|\widehat{\beta}^{(I)}_k-\beta_k| \le
\frac{2\widehat{\sigma}^{(I)}r}{\mathbb{E}_n[X_k^2]^{1/2} \kappa^{(I)*}_{k}
}\left(1-\frac{r}{\kappa_{1}^{(I)}}\right)_+^{-1}\,
\end{equation}
and
\begin{equation}\label{eq:t1RWI:4}
\widehat{\sigma}^{(I)}\le\max_{l=1,\hdots,L}({\bf D_Z^{(I)}})_{ll}\sqrt{\widehat{Q}_l(\beta)}.
\end{equation}
\end{theorem}
The same discussion that we had in the case of a high-dimensional
structural equation with endogenous regressors holds. The term
$$\left(1-\frac{r}{\kappa_{1}^{(I)}}\right)_+^{-1}$$ in the upper bounds
of Theorem \ref{tRWI} can yield infinite volume confidence regions.
This occurs either when $r$, of the order of $\sqrt{\log(L)/n}$ in
all scenarios, is large or when $\kappa_{1}^{(I)}$ is too small. The first
condition occurs when $L$ is as large as an exponential in $n$.  To
interpret the second condition, again Proposition \ref{p5} (ii)
yields that for any $k$ in $\{1,\hdots,K\}$,
$$\kappa_{1}^{(I)}\le\kappa_{k}^{(I)*}.$$
Therefore a necessary condition for $\kappa_{1}^{(I)}$ to be small is when
there exists in model \eqref{estruct} a regressor $(x_{ki})_{i=1}^n$
with low coordinate-wise sensitivity. Without the cone constraint
the coordinate wise sensitivities are defined through
\begin{equation*}
\kappa_{k}^{(I)*}=\inf_{\lambda\in\mathbb{R}^{K-1}}
\max_{l=1,\hdots,L}({\bf D_Z^{(I)}})_{ll}\left|\frac{1}{n}
\sum_{i=1}^n z_{li}\left(x_{ki}-x_{\{k\}^ci}^T\lambda\right)\right|
\end{equation*}
and can be interpreted as {\em maximal partial empirical
correlations with the instruments}.  This means that there is a
regressor for which all\footnote{Not as a set, see the discussion in
Section \ref{s4} about the comparison with the usual concentration
parameter.} instruments are weak in the sense that it is such that
$$\inf_{\lambda\in\mathbb{R}^{K-1}}
\max_{l=1,\hdots,L}({\bf D_Z^{(I)}})_{ll}\left|\frac{1}{n}
\sum_{i=1}^n
z_{li}\left(x_{ki}-x_{\{k\}^ci}^T\lambda\right)\right|\le r.$$
Again, in the setting of Example 7, this becomes increasingly likely
as $\zeta_{J^c}$ approaches zero (see \eqref{eInstrStrength}).
Rates of convergence could be obtained  as well as in Section
\ref{s54}.

\subsection{Potentially Sharper Results When None of the Instruments Have Heavy Tail}\label{sNHT}
In the subsequent sections we restrict the subset of indices $I$ in
\eqref{IVC} to a singleton (the index of the regressor which is
constant and equal to 1).  When $I$ is a singleton, the
{\em STIV} estimator is straightforward to obtain because
\eqref{IVC} only involves one conic constraint.  This procedure
should be used when all regressors and instruments are bounded or
have relatively thin tails.
In this section we consider the properties of the {\em STIV}
estimator with a different matrix ${\bf D}_{{\bf X}}$ in \eqref{IVS}.
We use for the matrix ${\bf D}_{{\bf X}}$ the diagonal $K\times
K$ matrix with diagonal entries $x_{k*}^{-1}$, $k=1,\dots,K$ where
$x_{k*}\triangleq\max_{i}|x_{ki}|$.

The following theorem which is in the same spirit as Theorem
\ref{t1}.
\begin{theorem}\label{t1HT}
For every $\beta$ in $\mathcal{B}_s$, under one of the 5 scenarios
of Section \ref{s51} together with its respective assumption and
choice of $r$, with probability at least $1-\alpha$ (approximately
at least $1-\alpha$ for scenario 4 and the two-stage procedure with
scenario 5), for any $c$ in $(0,1)$,
for any solution $(\widehat\beta^{(c)},\widehat\sigma^{(c)})$ of the
minimization problem (\ref{IVS}) we have for all $p\in[1,\infty]$
and $J_0\subset\{1,\hdots,K\}$,
\begin{equation}\label{eq:t1HT:1}
\left|\left({\bf D_X}^{-1} (\widehat{\beta}^{(c)}-\beta)\right)_{J_0}\right|_p \le
\frac{2\widehat{\sigma}^{(c)}r}{\kappa_{p,J_0,J(\beta)}^{(c)}}\left(1-\frac{r}{\kappa_{1,J_{{\rm exo}}^c,J(\beta)}^{(c)}}
-\frac{r^2}{\kappa_{1,J_{{\rm exo}},J(\beta)}^{(c)}}\right)_+^{-1}\,
\end{equation}
for all $k=1,\dots,K,$
\begin{equation}\label{eq:t1HT:3}
|\widehat{\beta}_k^{(c)}-\beta_k| \le
\frac{2\widehat{\sigma}^{(c)}r}{x_{k*}\,\, \kappa_{k,J(\beta)}^{(c)*}
}\left(1-\frac{r}{\kappa_{1,J_{{\rm exo}}^c,J(\beta)}^{(c)}}-\frac{r^2}
{\kappa_{1,J_{{\rm exo}},J(\beta)}^{(c)}}\right)_+^{-1}\
\end{equation}
and
\begin{equation}\label{eq:t1HT:4}
\widehat{\sigma}^{(c)}\le\sqrt{\widehat{Q}(\beta)}\left(1+\frac{r}{c\kappa_{1,J(\beta),J(\beta)}^{(c)}}\right)
\left(1-\frac{r}{c\kappa_{1,J(\beta),J(\beta)}^{(c)}}\right)_+^{-1}.
\end{equation}
\end{theorem}

Note that the upper bound \eqref{eq:t1:1} involves the term
$$\tau_1=1-\frac{r}{\kappa_{1,J(\beta)}^{(c)}}$$
while \eqref{eq:t1HT:1} in Theorem \ref{t1HT} involves
$$\tau_2\triangleq1-\frac{r}{\kappa_{1,J_{{\rm
exo}}^c,J(\beta)}^{(c)}} -\frac{r^2}{\kappa_{1,J_{{\rm exo}},J(\beta)}^{(c)}}.$$
Though the sensitivities in the two cases are not entirely
comparable because we use different matrices ${\bf D}_{{\bf X}}$, we
know, from Proposition \ref{p4} \eqref{p411} that restricting a
sensitivity to a sub-block increases the sensitivity.  When the set
$J_{{\rm exo}}^c$ of regressors that have not been identified as
exogenous and therefore used as instruments is small, the
sensitivity $\kappa_{1,J_{{\rm exo}}^c,J(\beta)}^{(c)}$ is large, whereas
the small sensitivity of its complement $\kappa_{1,J_{{\rm
exo}},J(\beta)}^{(c)}$ is compensated by the small value $r^2$ in the
numerator.  In the extreme case where $J_{{\rm exo}}^c=\varnothing$,
we have $\frac{r}{\kappa_{1,J_{{\rm exo}}^c,J(\beta)}^{(c)}}=0$, so that
$\tau_2\le 1-\frac{r^2}{\kappa_{1,J(\beta)}^{(c)}}$.  Under the premise of
Proposition \ref{p5}, for it to be positive it is sufficient to have
$|J(\beta)|\le Cr^{-2}=O(n/\log(L))$ where $C>0$ is a proper
constant.  Recall that in the case of Theorem \ref{t1} we need that
$|J(\beta)|\le Cr^{-1}=O(\sqrt{n/\log(L)})$.

The following results gives proper confidence sets.
\begin{corollary}\label{cor4}
For every $\beta$ in $\mathcal{B}_s$, under one of the 5 scenarios
of Section \ref{s51} together with its respective assumption and
choice of $r$, with probability at least $1-\alpha$ (approximately
at least $1-\alpha$ for scenario 4 and the two-stage procedure with
scenario 5), for any $c$ in $(0,1)$, for any solution $(\widehat\beta^{(c)},\widehat\sigma^{(c)})$ of the
minimization problem (\ref{IVS}) we have
\begin{equation}\label{eq:CIlpHT}
\left|\left({\bf D_X}^{-1} (\widehat{\beta}^{(c)}-\beta)\right)_{J_0}\right|_p \le
\frac{2\widehat{\sigma}^{(c)}r}{\overline{\kappa}_{p,J_0}^{(c)}(s)}\left(1-\frac{r}{\overline{\kappa}_{1,J_{{\rm exo}}^c}^{(c)}(s)}
-\frac{r^2}{\overline{\kappa}_{1,J_{{\rm exo}}}^{(c)}(s)}\right)_+^{-1}, \quad \forall \ p\in[1,\infty],\ \forall J_0\subset\{1,\hdots,K\},
\end{equation}
where $\overline{\kappa}_{p,J_0}^{(c)}(s)$, $\overline{\kappa}_{1,J_{{\rm
exo}}^c}^{(c)}(s)$, $\overline{\kappa}_{1,J_{{\rm exo}}}^{(c)}(s)$ are any lower
bounds on $\kappa_{p,J_0, J(\beta)}^{(c)}$, $\kappa_{1,J_{{\rm
exo}}^c,J(\beta)}^{(c)}$ and $\kappa_{1,J_{{\rm exo}},J(\beta)}^{(c)}$ based on
the sparsity certificates that is convenient to calculate (see, {\it
e.g.}, \eqref{eq:certif2} and \eqref{eq:certifgen}), and for all
$k=1,\dots,K$,
\begin{equation}\label{eq:CIHT}
|\widehat{\beta}_k^{(c)}-\beta_k| \le
\frac{2\widehat{\sigma}^{(c)}r}{x_{k*} \kappa^{(c)*}_{k}(s)
}\left(1-\frac{r}{\overline{\kappa}_{1,J_{{\rm exo}}^c}^{(c)}(s)}
-\frac{r^2}{\overline{\kappa}_{1,J_{{\rm exo}}}^{(c)}(s)}\right)_+^{-1}.
\end{equation}
\end{corollary}
As we have seen in Section \ref{s4}, equation \eqref{eq:certif2}
yields a lower bound for the block sensitivities.  It is only
computationally feasible for small blocks.  So this lower bound
could only be used at most for one of the two block sensitivities
appearing in $\tau_2$, the other one (possibly the two) should be
lower bounded using \eqref{eq:certifgen}.  When $J_{{\rm exo}}^c$ is
small there are cases when the method of this section yields finite
volume confidence regions while due to $\tau_1$ we can only deduce
from Theorem \ref{t1} infinite volume confidence regions.

It is easy to check that Theorem \ref{t1b} and Corollary \ref{cor2}
also hold replacing Assumption \ref{ass1X} by
\noindent\begin{assumption}\label{ass1XHT} For every $0<\gamma_3<1$, there
exist constants $v_k>0$ such that
$$
\mathbb{P}\left(x_{k*}\ge v_k, \ \forall \ k\in J(\beta)\right)\ge 1-\gamma_3.
$$
\end{assumption}
by assuming Assumption \ref{ass_select} allowing for $J_0$ to also
be $J_{{\rm exo}}$ and $J_{{\rm exo}}^c$ and using the following
expression for $\tau^{(c)*}$.
$$
\tau^{(c)*}\triangleq \left(1+\frac{r}{c c_{1,J(\beta)}^{(c)}}\right)
\left(1-\frac{r}{c
c_{1,J(\beta)}^{(c)}}\right)_+^{-1}\left(1-\frac{r}{c_{1,J_{\rm exo}^c}^{(c)}}
-\frac{r^2}{c_{1,J_{\rm exo}}^{(c)}}\right)_+^{-1}\,.
$$
Theorem \ref{th:threshold} also holds using
$$
\omega_{k}^{(c)}(s)\triangleq \frac{2\widehat\sigma^{(c)} r
}{\kappa_{k}^{(c)*}(s)x_{k*}}\left(1-\frac{r}{\kappa_{1,J_{\rm exo}^c}^{(c)}(s)
}-\frac{r^2}{\kappa_{1,J_{\rm exo}}^{(c)}(s)}\right)_+^{-1}\,,
$$
with
$$
\tau^{(c)*}(s)\triangleq \left(1+\frac{r}{c c_{1,J(\beta)}^{(c)}}\right)
\left(1-\frac{r}{c
c_{1,J(\beta)}^{(c)}}\right)_+^{-1}\left(1-\frac{r}{c_{1,J_{\rm
exo}^c}^{(c)}(s)} -\frac{r^2}{c_{1,J_{\rm exo}}^{(c)}(s)}\right)_+^{-1}\,.
$$

Note that due to the new scaling, the constant $v_k$ can be much
larger than in Section \ref{s54}.  This does not change the rates
when the regressors are bounded, but it does if the regressors are
unbounded.  If a regressor has a normal distribution, this implies
an increased $\log(n)$ in the rate.  When $\tau^{(c,I)*}$ would be infinite
under the scaling of Section \ref{s3} this is still a mild loss. We
present in Section \ref{sSim} the implementation of this method in
the difficult setting where the regressors are Gaussian.  This
method still performs well in this case.

\subsection{The {\em STIV} Estimator with Linear Projection
Instruments}\label{sLPI} The results of the previous sections show
that the {\em STIV} estimator can handle a very large number of
instruments, up to an exponential in the sample size. Moreover,
adding instruments always improves the sensitivities. In this
section, we consider the case where $L>K$ but we look for a smaller
set of instruments, namely, of size $K$.

This is a classical problem when the structural equation is
low-dimensional. The two-stage least square is a leading and very
popular example. Under the stronger zero conditional mean
assumption, optimal instruments attain the semi-parametric
efficiency bound (see Amemiya (1974), Chamberlain (1987), Newey
(1990) and Hahn (2002)). In the homoscedastic case, the optimal
instruments correspond to the projection of the endogenous variables
on the space of variables measurable with respect to all the
instruments\footnote{This is extremely hard in the presence of many
instruments due to the curse of dimensionality but it can be done
when the instruments are functionals of a low dimensional original
instrument.}.  These optimal instruments are expressed in terms of
conditional expectations that are not available in practice and
should be estimated. Related work includes, for example, Donald and
Newey (2001), Chamberlain and Imbens (2004), Bai and Ng (2009), Okui
(2010), Carrasco and Florens (2000) and Carrasco (2012). Belloni,
Chen, Chernozhukov et al.~(2010) is very related to our approach
because it proposes to use the Lasso or post-Lasso to estimate the
optimal instrument. This allows to approximate the conditional
expectations with a number of functionals of a low dimensional
vector of instrument that is exponential in the sample size.  They
consider a second stage which is the heteroscedastic robust {\em IV}
estimator and obtain confidence sets for the low dimensional
structural equation.

In this section we investigate the properties of a high-dimensional
version of the two-stage least squares for high-dimensional
structural equations with endogenous regressors. In that case, we
have $L>K$ and $K$ is potentially much larger than $n$. For
simplicity, we will not consider estimating the optimal instruments
in the first stage.  One would also have to properly define
optimality in high-dimension.  This is a difficult question that
should be addressed in the future.   Moreover, we are aiming at
robust confidence sets that have a good approximate finite
sample validity. For simplicity of exposure we assume that there is
only one endogenous regressor $(x_{k_{\rm end}i})_{i=1}^n$ in
\eqref{estruct}. We write the reduced form equation in the form
\begin{equation}\label{eq:reducedform}
x_{k_{\rm end}i}=\sum_{l=1}^Lz_{li}\zeta_{l}+v_{i},\quad i=1,\dots,n,
\end{equation}
where $\zeta_{l}$ are unknown coefficients of the linear combination
of instruments and
\begin{equation}
\label{eq:reducedform1} \E[z_{li}v_{i}]=0
\end{equation}
for $i=1,\dots,n$, $l=1,\dots,L$. We call
$\sum_{l=1}^Lz_{li}\zeta_{l}$ the linear projection instrument.

For simplicity of the proofs we use the normalization of Section
\ref{sNHT}. What we present in this section can easily be extended
to a situation with more than one endogenous regressor and with the
normalization of Section \ref{s3}.

The first stage consists in estimating the unknown coefficients
$\zeta_{l}$. If $L\ge K>n$ and if the reduced form model
(\ref{eq:reducedform}) has some sparsity (or approximate sparsity),
it is natural to use a high-dimensional procedure, such as the
Lasso, the Dantzig selector, the Square-root Lasso or the
post-Lasso, to produce estimators $\widehat{\zeta}_{l}$ of the
coefficients. Because we are after confidence sets for the
coefficients of the structural equation, we will assume that
\eqref{eq:reducedform} is sparse\footnote{Belloni, Chen,
Chernozhukov et al. (2010) also study the setting where the reduced
form has an approximately sparse vectors of coefficients.  We
could obtain rates of convergence in that case or asymptotically
valid confidence sets. Because we are aiming at confidence sets that
have some finite sample validity we do not study this case in this
article.}.  When the reduced form equation is not sparse or
approximately sparse, for example when all the instruments are
equally weak but relevant, this method will fail. This is what we
observe in the simulations in Section \ref{sSim}. In that situation
we recommend to use the one stage {\em STIV} estimator.

It is easy to check, from the first order condition that the {\em
STIV} estimator is, up to the normalization, equivalent to the
Square-root Lasso when all the regressors are exogenous. We can
deduce from Theorem \ref{t1HT}, that under one of the 5 scenarios,
on an event $E_{\alpha}$ of probability at least $1-\alpha$ (or
approximately at least $1-\alpha$ for scenario 4 and the two-stage
procedure with scenario 5), for any $c_{RF}$ in $(0,1)$,  any solution
$(\widehat{\zeta}^{(c_{RF})},\widehat{\sigma}_{1}^{(c_{RF})})$
to the {\em STIV} estimator applied to the reduced form equation is such that
\begin{equation*}
\left|\left({\bf D_X}^{-1} (\widehat{\zeta}^{(c_{RF})}-\zeta)\right)\right|_1 \le
\frac{2\widehat{\sigma}_{1}^{(c_{RF})}r}{\kappa_{1,J(\zeta)}^{(c_{RF})}}\left(1
-\frac{r^2}{\kappa_{1,J(\zeta)}^{(c_{RF})}}\right)_+^{-1}.
\end{equation*}
This yields, if one is given a maximum sparsity $s_{RF}$, the computable
upper bound
\begin{equation}\label{eq:up1stST}
\left|\left({\bf D_X}^{-1} (\widehat{\zeta}-\zeta)\right)\right|_1 \le
\frac{2\widehat{\sigma}_{RF}r}{\kappa_{1}^{(c_{RF})}(s_{RF})}\left(1
-\frac{r^2}{\kappa_{1}^{(c_{RF})}(s_{RF})}\right)_+^{-1}\triangleq C_1(r,c_{RF},s_{RF}).
\end{equation}

The second stage now consists in using the {\em STIV} estimator from
Section \ref{sNHT} with the enlarged {\em IV}-constraint
\begin{equation}\label{IVC2S}
\widehat{\mathcal{I}}_{2S}\triangleq\left\{(\beta,\sigma):\
\beta\in\R^K,\ \sigma>0,\ \left|\frac1n {\bf D}_{{\bf
Z},2S}\bold{Z}_{2S}^T(\bold{Y}-\bold{X}\beta)\right|_{\infty}\le \sigma
r,\ \widehat{Q}(\beta)\le \sigma^2\right\}
\end{equation}
where ${\bf D}_{{\bf Z},2S}$ is a $K\times K$ diagonal matrix such
that $({\bf D}_{{\bf Z},2S})_{k_{\rm end}k_{\rm
end}}=(\max_{i=1,\hdots,n}|z_i^T\widehat{\zeta}|+2C_1(r,c_{RF},s_{RF}))^{-1}$
and for every $k$ in $\{1,\hdots,K\}\setminus\{k_{\rm end}\}$,
$({\bf D}_{{\bf Z},2S})_{kk}=x_{k*}^{-1}$, and the matrix
$\bold{Z}_{2S}$ is the stacked matrix of the endogenous regressors
and the estimated linear projection instrument
$(z_i^T\hat{\lambda})_{i=1}^n$. We enlarge the {\em IV}-constraint
to account for the error made in the first stage estimation of the
linear projection instrument. In this section, for simplicity, we
restrict our attention to the case where $I$ is a singleton
corresponding to the index of the regressor being unity in
$\widehat{\mathcal{I}}_{2S}$. For this reason we drop the index
exponent $I$ everywhere. We obtain the following result where we
denote by $\kappa^{(c)2S}$ the sensitivities where we replace in the
definition of $\Psi_n^{(I)}$,  $\bold{Z}$ by  $\bold{Z}_{2S}$ and ${\bf
D_Z^{(I)}}$ by ${\bf D}_{{\bf Z},2S}$.

\begin{theorem}\label{t12S}
For every $\beta$ in $\mathcal{B}_s$, under one of the 5 scenarios
of Section \ref{s51} together with its respective assumption and
choice of $r$, with probability at least $1-\alpha$ (approximately
at least $1-\alpha$ for scenario 4 and the two-stage procedure with
scenario 5), for any $c$ in $(0,1)$, for any solution $(\widehat\beta^{(c)},\widehat\sigma^{(c)})$ of the
minimization problem (\ref{IVS}), replacing $\widehat{\mathcal{I}}$
by $\widehat{\mathcal{I}}_{2S}$, we have for all $p\in[1,\infty]$
and $J_0\subset\{1,\hdots,K\}$,
\begin{equation}\label{eq:t12S:1}
\left|\left({\bf D_X}^{-1} (\widehat{\beta}^{(c)}-\beta)\right)_{J_0}\right|_p \le
\frac{2\widehat{\sigma}^{(c)}r}{\kappa_{p,J_0,J(\beta)}^{(c)2S}}\left(1-\frac{r}{\kappa_{k_{\rm end},J(\beta)}^{(c)2S*}}
-\frac{r^2}{\kappa_{1,\{1,\hdots,K\}\setminus\{k_{\rm end}\},J(\beta)}^{(c)2S}}\right)_+^{-1}\,
\end{equation}
for all $k=1,\dots,K,$
\begin{equation}\label{eq:t12S:3}
|\widehat{\beta}_k^{(c)}-\beta_k| \le
\frac{2\widehat{\sigma}^{(c)}r}{x_{k*}\,\, \kappa_{k,J(\beta)}^{(c)2S*}
}\left(1-\frac{r}{\kappa_{k_{\rm end},J(\beta)}^{(c)2S*}}-\frac{r^2}
{\kappa_{1,\{1,\hdots,K\}\setminus\{k_{\rm end}\},J(\beta)}^{(c)2S}}\right)_+^{-1}\
\end{equation}
and
\begin{equation}\label{eq:t12S:4}
\widehat{\sigma}^{(c)}\le\sqrt{\widehat{Q}(\beta)}\left(1+\frac{r}{c\kappa_{1,J(\beta),J(\beta)}^{(c)2S}}\right)
\left(1-\frac{r}{c\kappa_{1,J(\beta),J(\beta)}^{(c)2S}}\right)_+^{-1}.
\end{equation}
\end{theorem}
In turn, this yields the following {\it bona-fide} confidence
sets.
\begin{corollary}\label{cor3}
For every $\beta$ in $\mathcal{B}_s$, under one of the 5 scenarios
of Section \ref{s51} together with its respective assumption and
choice of $r$, with probability at least $1-\alpha$ (approximately
at least $1-\alpha$ for scenario 4 and the two-stage procedure with
scenario 5), for any $c$ in $(0,1)$, for any solution $(\widehat\beta^{(c)},\widehat\sigma^{(c)})$ of the
minimization problem (\ref{IVS}), replacing $\widehat{\mathcal{I}}$
by $\widehat{\mathcal{I}}_{2S}$, we have
\begin{equation}\label{eq:CIlp2S}
\left|\left({\bf D_X}^{-1} (\widehat{\beta}^{(c)}-\beta)\right)_{J_0}\right|_p \le
\frac{2\widehat{\sigma}^{(c)}r}{\overline{\kappa}_{p,J_0}^{(c)2S}(s)}\left(1-\frac{r}{\kappa_{k_{\rm end}}^{(c)2S*}(s)}
-\frac{r^2}{\kappa_{1}^{(c)2S}(s)}\right)_+^{-1},
\quad \forall \ p\in[1,\infty],\ \forall J_0\subset\{1,\hdots,K\},
\end{equation}
where $\overline{\kappa}_{p,J_0}^{(c)2S}(s)$, is any lower bounds on
$\kappa_{p,J_0, J(\beta)}^{(c)2S}$ based on the sparsity certificates
that is convenient to calculate (see, {\it e.g.}, \eqref{eq:certif2}
and \eqref{eq:certifgen}), and for all $k=1,\dots,K$,
\begin{equation}\label{eq:CI2S}
|\widehat{\beta}_k^{(c)}-\beta_k| \le
\frac{2\widehat{\sigma}^{(c)}r}{x_{k*} \kappa^{(c)*}_{k}(s)
}\left(1-\frac{r}{\kappa_{k_{\rm end}}^{(c)2S*}(s)}
-\frac{r^2}{\kappa_{1}^{(c)2S}(s)}\right)_+^{-1}.
\end{equation}
\end{corollary}
All lower bounds on the sensitivities appearing in \eqref{eq:CIlp2S} and
\eqref{eq:CI2S}
are easy to calculate.  Indeed, we have
minorized $\kappa_{1,\{1,\hdots,K\}\setminus\{k_{\rm
end}\}}^{(c)2S}(s)$ by $\kappa_{1}^{(c)2S}(s)$ using Proposition \ref{p4}
(ii). It is also possible to proceed using the plug-in strategy
where we replace $J(\beta)$ by an estimate $\widehat{J}$ and obtain
confidence sets of level at least $1-\gamma$ where $\gamma>\alpha$
but $\gamma$ is very close to $\alpha$ under the assumption that the
non-zero coefficients are sufficiently separated from zero.

Rates of convergence and model selection results, similar to those
of Section \ref{s54}, could also be obtained.

%%%%%%%%%%%%%%%%%%%%%%%%%%%%%%%

\section{Models with Possibly Non-valid Instruments}\label{s6}  In this section, we study the problem of
checking for instrument exogeneity when there is overidentification.
This is a classic problem in econometrics (see, {\it e.g.}, Sargan
(1958) and Basmann (1960) for the linear IV model, and Hansen (1982)
for GMM).  See also Andrews (1999) and Liao~(2010) where one can
find more references.  We propose a two-stage method based on the
{\em STIV} estimator.   The main purpose of the suggested method is
to construct confidence sets for non-validity indicators, and
to detect non-valid ({\it i.e.}, endogenous) instruments in the
high-dimensional framework.   We restrict our attention to the case
where \eqref{estruct} is point identified, {\it i.e.}
$\mathcal{I}dent=\{\beta^*\}$. The model can be written in the form:
\begin{align}
&y_i=x_i^T\beta^*+u_i\label{estruct2},\\
&\E\left[z_iu_i\right]=0\label{einstr2},\\
&\E\left[\overline{z}_iu_i\right]=\theta^*,\label{enonvalidinstr2}
\end{align}
where $x_i,\ z_i,$ and ${\overline{z}}_i$ are vectors of dimensions
$K$, $L$ and $L_1$, respectively.  For simplicity of the discussion,
we assume in this section that \eqref{estruct2}-\eqref{einstr2} is
point identified.  for this reason we use the notation $\beta^*$
instead of $\beta$.  The instruments are decomposed in two parts,
$z_i$ and $\overline{z}_i$, where
$\overline{z}_i^T=(\overline{z}_{1i},\dots,\overline{z}_{L_1i})$ is
a vector of possibly non-valid instruments.  A component of the
unknown vector $\theta^*\in \R^{L_1}$ is equal to zero when the
corresponding instrument is indeed valid. The component $\theta_l^*$
of $\theta^*$ will be called the {\it non-validity indicator} of the
instrument $\overline{z}_{li}$. Our study covers the models with
dimensions $K$, $L$ and $L_1$ that can be much larger than the
sample size.

%The use of Lasso-type methods for estimation of $\beta^*$ and
%$\theta^*$ has been recently discussed by Liao (2010) who focused on
%asymptotic study of  when the sample size $n$ tends to $\infty$ and
%the dimensions of the unknown parameters stay fixed. Our purpose is
%quite different since we explore {\em IV} models with the dimensions
%$K$, $L$ and $L_1$
%that can be much larger than the sample size. %In this case,
%%classical arguments based on local Taylor expansion of the contrast
%%function and asymptotic normality are no longer applicable.

As above, we assume independence and allow for heteroscedasticity.
The difference from the previous sections is only in introducing
equation (\ref{enonvalidinstr2}).  We
observe realizations of independent random vectors
$(x_i^T,y_i,z_i^T,\overline{z}_{i}^T)$, $i=1,\dots,n$.
The components $\overline{z}_{li}$ satisfy $\E[\overline{z}_{li}
u_i]=\theta_l^*$ for all $l=1,\dots,L_1$, $i=1,\dots,n$.  We denote
by $\overline{\bold{Z}}$ the matrix of dimension $n\times L_1$ with
rows $\overline{z}_i^T$, $i=1,\hdots,n$. Set
$$
\overline{z}_{*}= \max_{l=1,\dots,L_1}\left(
\frac1{n}\sum_{i=1}^n\overline{z}_{li}^2\right)^{1/2}.
$$
In this section, we assume that we have a pilot estimator
$\widehat{\beta}$ and a statistic $\widehat{b}$ such that, with
probability close to 1,
\begin{equation}\label{eq:nv1}
\left|{\bf D_X}^{-1} (\widehat{\beta}-\beta^*)\right|_1
\le\widehat{b}.
\end{equation}
For example, $\widehat{\beta}$ can be the {\em STIV} estimator based
only on the vectors of valid instruments $z_1,\dots,z_n$.  In this
case, an explicit expression for $\widehat{b}$ can be obtained from
\eqref{eq:CIlp}\footnote{\eqref{eq:tapproxsparse} is not explicit.},
\eqref{eq:CIlpPI}, \eqref{eq:t1RWI:1} or \eqref{eq:CIHT}, depending
on the situation.

We define the
{\em STIV-NV} estimator $(\widehat{\theta}^{(c)}, \widehat{\sigma}_1^{(c)})$ as
any solution of the problem
\begin{equation}\label{def:STIV_est_nonvalid}
\min_{(\theta, \sigma_1)\in
\widehat{\mathcal{I}}_1}\big(|\theta|_1+c\sigma_1 \big),
\end{equation}
where $0<c<1$,
\begin{equation*}\label{IVC_NV}
\widehat{\mathcal{I}}_1\triangleq\left\{(\theta,\sigma_1):\
\theta\in\R^{L_1},\ \sigma_1>0,\ \left|\frac1n
\overline{\bold{Z}}^T(\bold{Y}-\bold{X}\widehat{\beta})-\theta\right|_{\infty}\le
\sigma_1 r_1+ \widehat{b}\overline{z}_*,\ F(\theta,\widehat{\beta}) \le
\sigma_1 + \widehat{b}\overline{z}_* \right\}
\end{equation*}
for some $r_1>0$ (to be specified below), where for all
$\theta=(\theta_1,\dots,\theta_{L_1})\in\R^{L_1}$, $\beta\in\R^{K}$,
$$
F(\theta,\beta)\triangleq\max_{l=1,\dots,L_1}\sqrt{
\widehat{Q}_l(\theta_l,\beta)}
$$
with
$$\widehat{Q}_l(\theta_l,\beta)\triangleq
\frac1n\sum_{i=1}^n\left(\overline{z}_{li}(y_i-x_i^T\beta)-\theta_l\right)^2.$$
It is not hard to see that the optimization problem
(\ref{def:STIV_est_nonvalid}) can be re-written as a conic program.

The following theorem provides a basis for constructing confidence
sets for the non-validity indicators.

\begin{theorem}\label{th:nonvalid}
Under one of the 5 scenarios of Section \ref{s51}, replacing
$z_{li}$ by $\overline{z}_{li}$ and $L$ by $L_1$, together with its
respective assumption and choice of $r$ with $\alpha$ replaced by
$\alpha_1$, when $\widehat{\beta}$ is an estimator satisfying
(\ref{eq:nv1}) with probability at least $1-\alpha_2$ for some
$0<\alpha_2<1$, then, with probability at least
$1-\alpha_1-\alpha_2$ (approximately at least $1-\alpha_1-\alpha_2$
for scenario 4 and the two-stage procedure with scenario 5), for any $c$ in $(0,1)$,
for any solution $(\widehat{\theta}^{(c)}, \widehat{\sigma}_1^{(c)})$ of the
minimization problem (\ref{def:STIV_est_nonvalid}), we have
\begin{equation}\label{eq:th:nonvalid:1}
|\widehat{\theta}^{(c)}-\theta^*|_\infty \le
\frac{2\Big[\widehat{\sigma}_1^{(c)}r_1 +
(1+r_1(1-c)^{-1})\widehat{b}\overline{z}_*\Big]}{(1-2r_1(1-c)^{-1}
|J(\theta^*)| )_+} \triangleq V(\widehat{\sigma}_1,c,\widehat{b},
|J(\theta^*)|)\,,
\end{equation}
and
\begin{equation}\label{eq:th:nonvalid:2}
|\widehat{\theta}^{(c)}-\theta^*|_1 \le
\frac{2\left[2|J(\theta^*)|\Big(\widehat{\sigma}_1^{(c)}r_1 +
(1+r_1)\widehat{b}\overline{z}_*\Big) + c\widehat{b}\overline{z}_*\right]}{(1-c-2r_1
|J(\theta^*)| )_+}\,.
\end{equation}
\end{theorem}
This theorem should be naturally applied when $r_1$ is small, {\em
i.e.}, $n\gg \log(L_1)$. In addition, we need a small $\widehat b$,
which is guaranteed by the results of Section \ref{s5} under the
condition $n\gg \log(L)$ if the pilot estimator $\widehat{\beta}$ is
the {\em STIV} estimator. Note also that the bounds
(\ref{eq:th:nonvalid:1}) and (\ref{eq:th:nonvalid:2}) are meaningful
if their denominators are positive, which is roughly equivalent to
the following bound on the sparsity of $\theta^*$: $|J(\theta^*)|=
O(1/r_1)= O(\sqrt{n/\log(L_1)})$.

Bounds for all the norms $|\widehat{\theta}^{(c)}-\theta^*|_p$, $\forall \
1<p<\infty$, follow immediately from (\ref{eq:th:nonvalid:1}) and
(\ref{eq:th:nonvalid:2}) by the standard interpolation argument. We
note that, in Theorem~\ref{th:nonvalid}, $\widehat{\beta}^{(c)}$ can be
any estimator satisfying (\ref{eq:nv1}), not necessarily the {\em
STIV} estimator.

We now consider that $c$ in $(0,1)$ is fixed.
To turn (\ref{eq:th:nonvalid:1}) and (\ref{eq:th:nonvalid:2}) into
valid confidence bounds, we can replace there $|J(\theta^*)|$ by
$|J(\widehat\theta^{(c)})|$, as follows from Theorem~\ref{th:nonvalid1}
(ii) below. In addition, Theorem~\ref{th:nonvalid1} establishes the
rate of convergence of the {\em STIV-NV} estimator and justifies the
selection of non-valid instruments by thresholding. To state the
theorem, we will need an extra assumption that the random variable
$F(\theta^*,\beta^*)$ is bounded in probability by a constant
$\sigma_{1*}>0$:

\begin{assumption}\label{ass_nonvalid2}
There exist constants $\sigma_{1*}>0$ and $0<\varepsilon<1$ such
that, with probability at least $1-\varepsilon$,
\begin{equation}\label{eq:ass_nonvalid2}
\max_{l=1,\hdots,L_1}\frac1{n}\sum_{i=1}^n \left(\overline{z}_{li}
u_i-\theta_l^*\right)^2 \le \sigma_{1*}^2.
\end{equation}
\end{assumption}

As in (\ref{eq:thresh}) we define a thresholded estimator
\begin{equation}\label{eq:thresh_nv}
\widetilde{\theta}_l^{(c)}\triangleq\left\{\begin{array}{ll}
                  \widehat{\theta}_l^{(c)} & {\rm if\ }
                  |\widehat{\theta}_l^{(c)}|>
                  \omega^{(c)}, \\
                  0 & {\rm otherwise,}
                \end{array}\right.
\end{equation}
where $\omega^{(c)}>0$ is some threshold. For $b_*>0$, $s_1>0$, define
$$
\overline{\sigma}_* =
\left(1-\frac{4r_1s_1}{c(1-c-2r_1s_1)_+}\right)_+^{-1}
\left[\sigma_{1*} +\frac{2b_*\overline{z}_*(1+2(1+r_1)s_1/c)}
{(1-c-2r_1s_1)_+}\right]\,.
$$

\begin{theorem}\label{th:nonvalid1}
Let the assumptions of Theorem \ref{th:nonvalid} and Assumption
\ref{ass_nonvalid2} be satisfied. Then the following holds.
\begin{itemize}
\item[(i)]
Let $\widehat{\beta}$ be an estimator satisfying
\begin{equation}\label{eq:th:nonvalid1:1}
\left|{\bf D_X}^{-1} (\widehat{\beta}-\beta^*)\right|_1 \le b_*
\end{equation}
with probability at least $1-\alpha_2$ for some $0<\alpha_2<1$ and
some constant $b_*$. Assume that $|J(\theta^*)|\le s_1$. Then, with
probability at least $1-\alpha_1-\alpha_2-\varepsilon$, for any
solution $\widehat{\theta}^{(c)}$ of the minimization problem
(\ref{def:STIV_est_nonvalid}) we have
\begin{equation}\label{eq:th:nonvalid1:2}
|\widehat{\theta}^{(c)}-\theta^*|_\infty \le
V({\overline\sigma}_{*},c, b_*,s_1).
%\frac{2\Big[{\overline\sigma}_{*}r_1 +
%(1+r_1(1-c)^{-1})\widehat{b}\overline{z}_*\Big]}{(1-2r_1(1-c)^{-1}
%|J(\theta^*)| )_+}\,.
\end{equation}
\item[(ii)] Let $(\widehat{\beta},\widehat{\sigma})$ be the STIV
estimator, and let the assumptions of all the items of Theorem
\ref{t1b} be satisfied (with $p=1$ in item {\rm (ii)}). Assume that
$|J(\theta^*)|\le s_1$, $|J(\beta^*)|\le s$, and $|\theta_l^*| >
V({\overline\sigma}_{*},c, b_*, s_1)$ for all $l\in J(\theta^*)$,
where
\begin{equation}\label{eq:th:nonvalid1:bstar}
b_*= \frac{2\sigma_*r s\tau^{(c,I)*}(s)}{c_1}\,.
\end{equation}
Then, with probability at least $1-\alpha_1-\varepsilon-\gamma$, for
any solution $\widehat{\theta}^{(c)}$ of the minimization problem
(\ref{def:STIV_est_nonvalid}) we have
\begin{equation}\label{eq:th:nonvalid1:jhat}
J(\theta^*) \subseteq J(\widehat{\theta}^{(c)}).
\end{equation}
\item[(iii)] Let the assumptions of item {\rm (ii)} and
Assumption~\ref{ass_select} hold.  Assume that $|\theta_l^*| >
2V({\overline\sigma}_{*},c, b_*, s_1)$ for all $l\in J(\theta^*)$.
Let $\tilde \theta$ be the thresholded estimator defined in
(\ref{eq:thresh_nv}) where $\widehat{\theta}^{(c)}$ is any solution of
the minimization problem (\ref{def:STIV_est_nonvalid}), and the
threshold is defined by $\omega^{(c)}= V(\widehat{\sigma}_1,c, \widehat
b, s_1)$ with
$$
\widehat b = \frac{2 \widehat{\sigma}
rs}{\kappa_{1}(s)}\left(1-\frac{r}{\kappa_{1}(s)
}\right)_+^{-1}\,.
$$
Then, with probability at least $1-\alpha_1-\varepsilon-\gamma$, we
have
\begin{equation}\label{eq:th:nonvalid1:3}
\overrightarrow{{\rm sign}
(\widetilde{\theta}^{(c)})}=\overrightarrow{{\rm sign}(\theta^*)}.
\end{equation}
As a consequence, $ J(\widetilde{\theta}^{(c)})=J(\theta^*)$.
\end{itemize}
\end{theorem}

In practice, the parameter $s$ %in Theorem~\ref{th:nonvalid1} (iii)
may not be known and it can be replaced by $|J(\widehat{\theta})|$;
this is a reasonable upper bound on $|J(\theta^*)|$ as suggested by
Theorem~\ref{th:nonvalid1} (ii). It is interesting to analyze the
dependence of the rate of convergence in (\ref{eq:th:nonvalid1:2})
on $r,r_1,s$, and $s_1$. As discussed above, a meaningful framework
is to consider small $r$,$r_1$ and the sparsities $s$,$s_1$ such
that $rs$, $r_1s_1$ are comfortably smaller than 1. In this case,
the value $b_*$ given in (\ref{eq:th:nonvalid1:3}) is of the order
$O(rs)$ and the rate of convergence in (\ref{eq:th:nonvalid1:2}) is
of the order $O(r_1) + O(rs)$. We see that the rate does not depend
on the sparsity $s_1$ of $\theta^*$ %provided $s_1\ll r_1^{-1}$
but it does depend on the sparsity $s$ of $\beta^*$. It is
interesting to explore whether this rate is optimal, {\em i.e.},
whether it can be improved by estimators different from the {\em
STIV-NV} estimator.

%%%%%%%%%%%%%%%%%%%%%%

\section{Practical Implementation}\label{s8}
\subsection{Computational Aspects}\label{s900}
For simplicity of exposition, in this section we only consider the case where
$I=\{i_0\}$ and $i_0$ is the index of unity.  Extension to several cones is easy. 
Finding a solution 
$\left(\widehat{\beta}^{(c,\{i_0\})},\widehat{\sigma}^{(c,\{i_0\})}\right)$ of the
minimization problem \eqref{IVS} reduces to the following conic
program: find $\beta\in\R^K$ and $t>0$ ($\sigma=t/\sqrt{n}$), which
achieve the minimum
\begin{equation}\label{LinProg}
\min_{(\beta,t,v,w)\in {\mathcal V}}
\left(\sum_{k=1}^Kw_k+c\frac{t}{\sqrt{n}}\right)
\end{equation}
where ${\mathcal V}$ is the set of $(\beta,t,v,w)$, with
satisfying:
\begin{eqnarray*}
&& v=\bold{Y}-\bold{X}\beta, \qquad -rt\bold{1}\le
\frac{1}{\sqrt{n}}{\bf D}_{{\bf
Z}}\bold{Z}^T\left(\bold{Y}-\bold{X}\beta\right)\le rt\bold{1},\\
&&  -w\le {\bf D}_{{\bf X}}^{-1}\beta\le w, \qquad w\ge\bold{0},
\qquad (t,v)\in C.
\end{eqnarray*}
Here and below $\bold{0}$ and $\bold{1}$ are vectors of zeros and
ones respectively, the inequality between vectors is understood in
the componentwise sense, and $C$ is a cone:
$C\triangleq\{(t,v)\in\R\times\R^n:\ t\ge|v|_2\}$. Conic programming
is a standard tool in optimization and many open source toolboxes
are available to implement it (see, {\em e.g.}, Sturm (1999)).

The expression in curly brackets in the lower bound
\eqref{eq:certif1} is equal to the value of the following
optimization
program: %find $v\in\R$, which achieves the minimum
\begin{equation}\label{LPcertif}
\min_{\epsilon=\pm1}\min_{(w,\Delta,v)\in {\mathcal V}_{k,j}} v
\end{equation}
where ${\mathcal V}_{k,j}$ is the set of $(w,\Delta,v)$ with
$w\in\R^K$, $\Delta\in\R^K$, $v\in \R$ satisfying:
\begin{eqnarray*}
&& v\ge0, \qquad -v\bold{1}\le\Psi_n^{(I)}\Delta\le v\bold{1}, \qquad w\ge
\bold{0}, \qquad -w_{I^c}\le\Delta_{I^c}\le w_{I^c} \ \ {\rm for} \
I=\{j,k\},\\
&& w_{I}=\bold{0}, \qquad \Delta_k=1, \qquad \epsilon\Delta_j\ge0,
\qquad \sum_{i=1}^{K}w_i+1\le
       \epsilon(a+g)\Delta_j
\end{eqnarray*}
where $g$ is the constant such that
$$g= \left\{\begin{array}{ll}
                    0&  {\rm if\ } k= j\\
                  -1 & {\rm otherwise.}
\end{array}\right.$$
Note that, here, $\epsilon$ is the sign of $\Delta_j$, and
(\ref{LPcertif}) is the minimum of two terms, each of which is the
value of a linear program. Analogously, the expression in curly
brackets in \eqref{eq:certif2} can be computed by solving
$2^{|J_0|}$ linear programs. The reduction is done in the same way
as in (\ref{LPcertif}) with the only difference that instead of
$\epsilon$ we introduce a vector $(\epsilon_k)_{k\in J_0}$ of signs
of the coordinates $\Delta_k$ for indices $k\in J_0$.

The coordinate-wise sensitivities
\begin{equation*}
\kappa_{k,J}^{(c,\{i_0\})*}= \inf_{\ \Delta_k= 1,\,|\Delta_{J^c}|_1\le
\frac{1+c}{1-c}|\Delta_J|_1}\left|\Psi_n^{(I)}\Delta \right|_{\infty}
\end{equation*}
can be efficiently computed for given $J$ when the cardinality $|J|$
is small. Indeed, it is enough to find the minimum of the values of
$2^{|J|}$ linear programs:
\begin{equation}\label{LPkappa}
\min_{(\epsilon_j)_{j\in J}\in\{-1,1\}^{|J|}}\min_{(w,\Delta,v)\in
{\mathcal U}_{k,J}} v
\end{equation}
where ${\mathcal U}_{k,J}$ is the set of $(w,\Delta,v)$ with
$w\in\R^K$, $\Delta\in\R^K$, $v\in \R$ satisfying:
\begin{eqnarray*}
&& v\ge0, \qquad -v\bold{1}\le\Psi_n^{(I)}\Delta\le v\bold{1}, \qquad w\ge
\bold{0}, \qquad -w_{I^c}\le\Delta_{I^c}\le w_{I^c} \ \ {\rm for} \
I=J\cup\{k\},\\
&& w_{I}=\bold{0}, \qquad \Delta_k=1, \qquad \epsilon_j\Delta_j\ge0,
\ \ {\rm for} \ j\in J, \qquad \sum_{i=1}^{K}w_i\le
       \frac{1+c}{1-c}\sum_{j\in J}\epsilon_j\Delta_j +g.
\end{eqnarray*}
Here $(\epsilon_j)_{j\in J}$ is the vector of signs of the
coordinates $\Delta_j$ with $j\in J$ and $g$ is the constant defined
by
$$
g= \left\{\begin{array}{ll}
                  0 & {\rm if\ }
                  k\in J,\\
                  -1 & {\rm otherwise.}
                \end{array}\right.
$$

\subsection{Simulations}\label{sSim}
In this section, we consider the performance of the {\em STIV}
estimator of Section \ref{sNHT} on simulated data.  The model is as
follows:
\begin{align*}
&y_i=\sum_{k=1}^{K}x_{ki}\beta_k^*+u_i,\\
&x_{1i}=\sum_{l=1}^{L-K+1}z_{li}\zeta_l+v_i, \\
&x_{l'i}=z_{li} \quad \mbox{for} \ l'=l-L+K \ \  \mbox{and} \ \
l\in\{L-K+2,\hdots,L\},
\end{align*}
where $(y_i,x_i^T,z_i^T,u_i,v_i)$ are i.i.d., $(u_i,v_i)$ have the
joint normal distribution
$$\mathcal{N}\left(0,\left(\begin{array}{cc}
\sigma_{{\rm struct}}^2 & \rho\sigma_{{\rm struct}}\sigma_{{\rm end}} \\
\rho\sigma_{{\rm struct}}\sigma_{{\rm end}} & \sigma_{{\rm end}}^2
\end{array}\right)\right),$$
$z_i^T$ is a vector of independent standard normal random variables,
and $z_i^T$ is independent of $(u_i,v_i)$. Clearly, in this model
$\E[z_iu_i]=0$.  We take $n=49$, $L=50$, $K=25$, $\sigma_{{\rm
struct}}=\sigma_{{\rm end}}=\rho=0.3$,
$\beta^*=(1,1,1,1,1,0,\hdots,0)^T$ and $\zeta_l=0.15$ for
$l=1,\hdots,L-K+1$. We have 50 instruments and only 49 observations,
so we are in a framework of application of high-dimensional
techniques. We set $c=0.1$ and set $r$ according to \eqref{er3} with
$\alpha=0.05$.  The three columns on the left of Table \ref{Tab1}
present simulation results for the {\em STIV} estimator. It is
straightforward to see that only the first five variables (the true
support of $\beta^*$) are eligible to be considered as relevant.
This set will be denoted by $\widehat{J}$. The second and third
columns in Table \ref{Tab1} present the true coordinate-wise
sensitivities $\kappa_{k,\widehat{J}}^*$ as well as their lower
bounds $\kappa_{k}^*(5)$ obtained via the sparsity certificate with
$s=5$. These lower bounds are easy to compute, and we see that they
yield reasonable approximations from below of the true
sensitivities. The estimate $\widehat{\sigma}$ is 0.247 which is
quite close to $\sigma_{{\rm struct}}$. Next, based on
(\ref{eq:t1:3}), the fact that $J_{{\rm exo}}^c=\{1\}$, and the
bounds on the sensitivities in Proposition \ref{p4} and in
(\ref{eq:certif1}) -- (\ref{eq:certif2}), we have the following
formulas for the confidence intervals
\begin{equation}\label{eq:simul_conf1}
|\widehat{\beta}_k-\beta^*_k| \le
\frac{2\widehat{\sigma}r}{x_{k*}\,\, \kappa^*_{k,\widehat{J}}
}\left(1-\frac{r}{\kappa_{1,\widehat{J}}^*}-\frac{r^2}
{\kappa_{1,\widehat{J}}}\right)_+^{-1} \,,%\triangleq
%\widehat{u}_{\widehat{J}}\,,
\end{equation}
\begin{equation}\label{eq:simul_conf2}
|\widehat{\beta}_k-\beta^*_k| \le
\frac{2\widehat{\sigma}r}{x_{k*}\,\, \kappa^*_{k}(s)
}\left(1-\frac{r}{\kappa_{1}^*(s)}-\frac{r^2}
{\kappa_{1}(s)}\right)_+^{-1}%\triangleq  \widehat{u}(s)
\, \quad {\rm
with} \ s=5.
\end{equation}
Here, $\kappa^*_{k,\widehat{J}}$ and $\kappa^*_{k}(s)$ are computed
directly via the programs (\ref{LPkappa}) and (\ref{LPcertif})
respectively. The value $\kappa_{1}(s)$ is then obtained from
\eqref{eq:kappa1}, and for $\kappa_{1,\widehat{J}}$ we use a lower
bound analogous to \eqref{eq:kappa1}:
%$\kappa_{1,\widehat{J}}\ge
%\underline\kappa_{1,\widehat{J}}$, where
$$
%\underline\kappa_{1,\widehat{J}}=
\kappa_{1,\widehat{J}}\ge
\frac{1-c}{2|\widehat{J}|}\min_{k=1,\hdots,K}
\kappa_{k,\widehat{J}}^*.
$$
%{\bf ???? In all the examples that we considered, the terms with
%$r^2$ in the brackets in (\ref{eq:simul_conf1}) and
%(\ref{eq:simul_conf2}) were negligible.}
We get
$\kappa_{1,\widehat{J}}^*=0.0096$ and $\kappa^*_{1}(5)=0.0072$. In
particular, we have $r/{\kappa}_{1,\widehat{J}}^*=4.40>1$, so that
(\ref{eq:simul_conf1}) and (\ref{eq:simul_conf2}) do not provide
confidence intervals in this numerical example.

%and $\kappa^*_{k}(s)$ from From \eqref{eq:kappa} we obtain 0.107 and
%0.080 as lower bounds on $\kappa_{\infty,J(\beta^*)}$ based on the
%coordinate-wise sensitivities, with $J(\beta^*)$ estimated by
%$\widehat{J}$ and by the sparsity certificate respectively. These
%yield 0.0096 and 0.0072 as lower bounds on $\kappa_{1,J(\beta^*)}$
%using \eqref{k2}. Unfortunately, in this case inequality
%\eqref{eq:t1:3} does not provide meaningful confidence intervals
%since the right-hand side of \eqref{eq:t1:3} turns out to be greater
%than 1.
%
%because $r/\underline{\kappa}_{1,\widehat{J}}^*=4.40>1$, where we
%denote by $\underline{\kappa}_{1,\widehat{J}}^*$ the lower bound,
%and of equation \eqref{eq:t1:3}, we are not able to compute the
%confidence intervals or to apply the thresholding strategy to
%estimate of the support.

\begin{table}%[H]
  \centering
  \caption{Results for the {\em STIV} estimator without and with estimated instruments, $n=49$}\label{Tab1} {\footnotesize
\begin{tabular}{|c||c|c|c||c|c|c|} \hline
& $\hat{\beta} \ \ \ \ (1)$ & $\kappa_{k,\widehat{J}}^*\ \ \ (1)$ & $\kappa_{k}^*(5)\ \ \ (1)$  & $\hat{\beta}\ \ \ \ (2)$ & $\kappa_{k,\widehat{J}}^*\ \ \ (2)$ & $\kappa_{k}^*(5)\ \ \ (2)$ \\
\hline
\cellcolor[gray]{.7} $\beta_1^*$ &  \cellcolor[gray]{.7}    1.03    & \cellcolor[gray]{.7}  0.107   &      \cellcolor[gray]{.7}  0.103       &  \cellcolor[gray]{.7}              1.03                                     & \cellcolor[gray]{.7}  0.085   & \cellcolor[gray]{.7}  0.068   \\
\cellcolor[gray]{.7} $\beta_2^*$ &  \cellcolor[gray]{.7}    0.98    &  \cellcolor[gray]{.7} 0.308   &      \cellcolor[gray]{.7}  0.157       &  \cellcolor[gray]{.7}              0.98                                     & \cellcolor[gray]{.7}  0.367   & \cellcolor[gray]{.7}  0.075   \\
\cellcolor[gray]{.7} ${\beta}_3^*$ &  \cellcolor[gray]{.7}      0.96    &  \cellcolor[gray]{.7} 0.129   &  \cellcolor[gray]{.7}  0.103       &  \cellcolor[gray]{.7}              0.96                                     & \cellcolor[gray]{.7}  0.126   &  \cellcolor[gray]{.7} 0.071   \\
\cellcolor[gray]{.7}  ${\beta}_4^*$ &  \cellcolor[gray]{.7}     0.95    & \cellcolor[gray]{.7}  0.150   &  \cellcolor[gray]{.7}  0.109       &  \cellcolor[gray]{.7}              0.95                                     & \cellcolor[gray]{.7}  0.115   &  \cellcolor[gray]{.7} 0.057   \\
\cellcolor[gray]{.7}  ${\beta}_5^*$ &  \cellcolor[gray]{.7}     0.90    & \cellcolor[gray]{.7}  0.253   &  \cellcolor[gray]{.7} 0.175        &  \cellcolor[gray]{.7}              0.90                                     & \cellcolor[gray]{.7}  0.177   &  \cellcolor[gray]{.7} 0.086   \\
${\beta}_6^*$ &         0.00    &   0.166   &   0.095   &   0.00    &   0.126   &   0.065   \\
${\beta}_7^*$ &         0.00    &   0.155   &   0.080   &   0.00    &   0.148   &   0.060   \\
${\beta}_8^*$ &         0.00    &   0.154   &   0.110   &   0.00    &   0.122   &   0.056   \\
\vdots &\vdots&\vdots&\vdots&\vdots&\vdots&\vdots\\
%${\beta}_9^*$ &         -0.01   &   0.198   &   0.126   &   -0.01   &   0.176   &   0.095   \\
%${\beta}_{10}^*$ &      0.00    &   0.344   &   0.213   &   0.00    &   0.276   &   0.103   \\
%${\beta}_{11}^*$ &      0.02    &   0.224   &   0.151   &   0.02    &   0.158   &   0.112   \\
%${\beta}_{12}^*$ &      0.03    &   0.275   &   0.152   &   0.03    &   0.236   &   0.089   \\
%${\beta}_{13}^*$ &      0.00    &   0.205   &   0.129   &   0.00    &   0.150   &   0.091   \\
%${\beta}_{14}^*$ &      0.00    &   0.351   &   0.240   &   0.00    &   0.288   &   0.156   \\
%${\beta}_{15}^*$ &      0.09    &   0.168   &   0.110   &   0.09    &   0.155   &   0.083   \\
%${\beta}_{16}^*$ &      0.00    &   0.209   &   0.119   &   0.00    &   0.166   &   0.078   \\
%${\beta}_{17}^*$ &      -0.04   &   0.176   &   0.086   &   -0.04   &   0.120   &   0.044   \\
%${\beta}_{18}^*$ &      0.00    &   0.202   &   0.150   &   0.00    &   0.170   &   0.107   \\
%${\beta}_{19}^*$ &      0.04    &   0.173   &   0.124   &   0.04    &   0.148   &   0.089   \\
%${\beta}_{20}^*$ &      -0.05   &   0.284   &   0.213   &   -0.05   &   0.233   &   0.157   \\
%${\beta}_{21}^*$ &      0.00    &   0.191   &   0.114   &   0.00    &   0.150   &   0.065   \\
%${\beta}_{22}^*$ &      0.00    &   0.296   &   0.240   &   0.00    &   0.284   &   0.193   \\
${\beta}_{23}^*$ &      0.02    &   0.287   &   0.170   &   0.02    &   0.231   &   0.128   \\
${\beta}_{24}^*$ &      0.00    &   0.243   &   0.137   &   0.00    &   0.195   &   0.105   \\
${\beta}_{25}^*$ &      0.00    &   0.141   &   0.109   &   0.00    &   0.106   &   0.067   \\
\hline
\end{tabular}\\
We use dots because the values that do not appear are similar.\\
(1): With all the 50 instruments,\\
(2): With 25 instruments including an estimate of the linear
projection instrument.}
\end{table}
 The columns on the right in Table \ref{Tab1} present
the results where we use the same data, estimate the linear
projection instrument by the Square-root Lasso and then take only
$K$ instruments: $z_{il}$, $l=L-K+2,\hdots,L$, and
$\widehat{x}_{i1}=\sum_{l=1}^{L}\widehat{\zeta}_l z_{il}$, where
$\widehat{\zeta}_l$ are the Square-root Lasso estimators of
${\zeta}_l$, $l=1,\hdots,L$.
 The
Square-root Lasso with parameter $c_{\sqrt{\rm Lasso}}=1.1$
recommended in Belloni, Chernozhukov and Wang (2010) \footnote{One should not confuse
the constant $c_{\sqrt{\rm Lasso}}$ denoted by $c$ in Belloni,
Chernozhukov and Wang (2010) with the constant
$c=c_{STIV}$ in the definition of the {\em STIV} estimator;
$c_{\sqrt{\rm Lasso}}$ is an equivalent of $\sqrt{n}/c_{STIV}$, up
to constants.} yields all coefficients equal to zero when keeping
only the first three digits. This is disappointing since we get an
instrument equal to zero. It should be noted that estimation in this
setting is a hard problem since the dimension $L$ is larger than the
sample size, the number of non-zero coefficients ${\zeta}_l$ is
large ($L-K+1=26$), and their values are relatively small (equal to
$0,15$). To improve the estimation, we adjusted the parameter
$c_{\sqrt{\rm Lasso}}$ empirically, based on the value of the estimates. %and of
%$\sqrt{\widehat{Q}(\widehat{\beta})}$ {\bf je ne crois pas pour le deuxieme}.
Ultimately, we have chosen
$c_{\sqrt{\rm Lasso}}=0.3$. This choice is not covered by the theory
of Belloni, Chernozhukov and Wang (2010) because there $c_{\sqrt{\rm
Lasso}}$ should be greater than~1. However, it leads to
$\sqrt{\widehat{Q}(\widehat{\beta})}=0.309$, which is very close to
$\sigma_{{\rm end}}$. The corresponding estimates
$\widehat{\zeta}_l$ are given in Table \ref{Tab2}. We see that they
are not very close to the true ${\zeta}_l$; some of the relevant
coefficients are erroneously set to 0 and several superfluous
variables are included, sometimes with significant coefficients,
such as $\widehat{\zeta}_{32}$.
\begin{table}%[H]
  \centering
  \caption{Estimates of the coefficients of the linear projection
  instrument}\label{Tab2}
{\footnotesize \begin{tabular}{|c|c|c|c|c|c|c|c|c|c|c|c|c|c|}
  \hline
$\widehat{\zeta}_1$    &   $\widehat{\zeta}_2$    &   $\widehat{\zeta}_3$    &   $\widehat{\zeta}_4$      &   $\widehat{\zeta}_6$    &   $\widehat{\zeta}_8$    &   $\widehat{\zeta}_9$    &   $\widehat{\zeta}_{10}$ & $\widehat{\zeta}_{14}$ &   $\widehat{\zeta}_{15}$ &   $\widehat{\zeta}_{16}$ &   $\widehat{\zeta}_{17}$ &   $\widehat{\zeta}_{18}$ & $\widehat{\zeta}_{20}$ \\
0.084                   &       0.130                &           0.190          &       0.142              &       0.115               &        0.083           &        0.104              &            0.126          &            0.176    &         0.030             &         0.023             &         0.157             &         0.135          &         0.082          \\
\hline
 $\widehat{\zeta}_{21}$  &   $\widehat{\zeta}_{23}$ &   $\widehat{\zeta}_{24}$ &   $\widehat{\zeta}_{25}$ &   $\widehat{\zeta}_{26}$ &$\widehat{\zeta}_{27}$ &    $\widehat{\zeta}_{32}$ &   $\widehat{\zeta}_{33}$ &   $\widehat{\zeta}_{34}$  &  $\widehat{\zeta}_{44}$   &    $\widehat{\zeta}_{47}$ &      $\widehat{\zeta}_{49}$ &   $\widehat{\zeta}_{50}$& \\
       0.100                      &         0.125             &         0.038             &        0.025              &            0.026      &      -0.058           &            0.108            &         0.005             &         -0.053         &        -0.006            &       -0.009                       &   -0.063               & 0.033 & \\
\hline
\end{tabular}\\
%&   $\widehat{\zeta}_{22}$
%&         0.000
We only show the non-zero coefficients.}
\end{table}
%We obtain 0.085 and 0.044 as lower bounds on
%$\kappa_{\infty,J(\beta^*)}$ based on the estimated coordinate-wise
%sensitivities and the sparsity certificate with $s=5$ respectively.
We get $\kappa_{1,\widehat{J}}^*=0.0076$ and
$\kappa^*_{1}(5)=0.0040$. Again, $r/{\kappa}_{1,\widehat{J}}^*>1$,
so that we cannot use (\ref{eq:simul_conf1}) and
(\ref{eq:simul_conf2}) to get the confidence intervals. Note that
this approach based on the estimated linear projection instrument
gives sensitivities, which are lower than with the full
set of instruments.  %This could seem counter-intuitive but recall
%that a sensitivity is a minimum over a set which allows combinations
%of the instruments.  Also, because sensitivities involve the maximum
%of scalar products of the rows of $\Psi_n$ (instruments) with
%$\Delta$, the more we have instruments the higher is the
%sensitivity.  On the other hand, $\sqrt{\log(L)/\log(K)}=1,102$,
%which is not enough to compensate the overall deterioration of the
%sensitivities.
This is mainly due to the fact that the estimation of the linear
projection instrument is quite imprecise. Indeed, we add an
instrument $\widehat{x}_{i1}$, which is not so good, and at the same
time we drop a large number of other instruments, which may be not
so bad. The overall effect on the sensitivities turns out to be
negative. Recall that since the sensitivities involve the maximum of
the scalar products of the rows of $\Psi_n^{(I)}$ with $\Delta$, the more
we have rows ({\it i.e.}, instruments) the higher is the
sensitivity.
%some error, which can be
%large (infinite with the choice of $c=1.1$ in the Square-root
%Lasso).
The same deterioration of the sensitivities occurred in other
simulated data sets. In conclusion, the approach based on estimation
of the linear projection instrument was not helpful to realize the
above confidence intervals in this small sample situation. However,
we will see that it achieves the task when the sample size gets
large.

Although in this numerical example we were not able to use
(\ref{eq:simul_conf1}) and (\ref{eq:simul_conf2}) for the confidence
intervals, we got evidence that the performance of the {\em STIV}
estimator is quite satisfactory. Table \ref{Tab3} shows a
Monte-Carlo study where we keep the same values of the parameters of
the model, of the sample size $n=49$, and of the parameter $A$
defining the set $\widehat{\mathcal{I}}$, simulate 1000 data sets,
and compute 1000 estimates.
\begin{table}%[H]
  \centering
\caption{Monte-Carlo study, 1000 replications}\label{Tab3}
{\footnotesize
\begin{tabular}{|c|c|c|c||c|c|c|c|} \hline
& $5^{th}$ percentile & Median  & $95^{th}$ percentile & & $5^{th}$ percentile & Median  & $95^{th}$ percentile\\
\hline
\cellcolor[gray]{.7} $\beta_1^*$     &  \cellcolor[gray]{.7} 0.872   &  \cellcolor[gray]{.7} 0.986   & \cellcolor[gray]{.7}  1.093   &   $\beta_8^*$     &   -0.057  &   0.000   &   0.055    \\
\cellcolor[gray]{.7} $\beta_2^*$     &  \cellcolor[gray]{.7} 0.877   &  \cellcolor[gray]{.7} 0.970   & \cellcolor[gray]{.7}  1.048   &   $\beta_9^*$     &   -0.052  &   0.000   &   0.059    \\
\cellcolor[gray]{.7} $\beta_3^*$     &  \cellcolor[gray]{.7} 0.879   &  \cellcolor[gray]{.7} 0.970   & \cellcolor[gray]{.7}  1.049   &    \vdots &  \vdots  & \vdots  &   \vdots     \\
\cellcolor[gray]{.7} $\beta_4^*$     &  \cellcolor[gray]{.7} 0.886   &  \cellcolor[gray]{.7} 0.971   & \cellcolor[gray]{.7}  1.051   &    $\beta_{23}^*$  &   -0.051  &   0.000   &   0.051   \\
\cellcolor[gray]{.7} $\beta_5^*$     &  \cellcolor[gray]{.7} 0.877   &  \cellcolor[gray]{.7} 0.968   & \cellcolor[gray]{.7}  1.049   &    $\beta_{24}^*$  &   -0.057  &   0.000   &   0.051   \\
$\beta_6^*$ &   -0.048  &   0.000   &   0.055   &                                                                                  $\beta_{25}^*$  &   -0.053  &   0.000   &   0.049   \\
$\beta_7^*$     &   -0.059  &   0.000   &   0.063   &                                                                              $\hat{\sigma}$  &   0.181   &   0.233   &   0.291   \\
%$\beta_{19}^*$  &   -0.059  &   0.000   &   0.056   &
% $\beta_{20}^*$  &   -0.049  &   0.000   &   0.055  &
%
%$\beta_8^*$     &   -0.057  &   0.000   &   0.055   &   $\beta_{21}^*$  &   -0.050  &   0.000   &   0.052   \\
%$\beta_9^*$     &   -0.052  &   0.000   &   0.059   &   $\beta_{22}^*$  &   -0.052  &   0.000   &   0.050   \\
%$\beta_{10}^*$  &   -0.050  &   0.000   &   0.050   &   $\beta_{23}^*$  &   -0.051  &   0.000   &   0.051   \\
%$\beta_{11}^*$  &   -0.053  &   0.000   &   0.052   &   $\beta_{24}^*$  &   -0.057  &   0.000   &   0.051   \\
%$\beta_{12}^*$  &   -0.053  &   0.000   &   0.055   &   $\beta_{25}^*$  &   -0.053  &   0.000   &   0.049   \\
%$\beta_{13}^*$  &   -0.059  &   0.000   &   0.050   &   $\hat{\sigma}$  &   0.181   &   0.233   &   0.291   \\
\hline
\end{tabular}}
\end{table}
The empirical performance of the {\em STIV} estimator is extremely
good, even for the endogenous variable. The Monte-Carlo estimation
of the variability of $\widehat{\beta}_1$ is very similar to that of
the exogenous variables. With $c=0.1$ the estimate
$\widehat{\sigma}$ is larger than $\sigma_{{\rm struct}}$ in 95\% of
the simulations. This suggests that there remains some margin to
penalize less for the ``variance" in (\ref{IVS}), {\em i.e.}, to
decrease $c$ and thus to obtain higher sensitivities.

Next, we study the empirical behavior of the non-pivotal {\em STIV}
estimator.  We consider the same model and the same values of all
the parameters, and we choose $\sigma_*=2\cdot 0.233$ where $0.233$
is the median of $\widehat{\sigma}$ from Table \ref{Tab3}. Indeed
$\mathbb{P}\left(\mathbb{E}_n[U^2]\le\sigma_*^2\right)$ should be
close to 1 (see Assumption \ref{ass1b}). The results are given in
Table \ref{Tab4}.
\begin{table}%[H]
\centering \caption{Monte-Carlo study of the non-pivotal estimator,
1000 replications}\label{Tab4} {\footnotesize
\begin{tabular}{|c|c|c|c||c|c|c|c|} \hline
& $5^{th}$ percentile & Median  & $95^{th}$ percentile & & $5^{th}$ percentile & Median  & $95^{th}$ percentile\\
\hline
\cellcolor[gray]{.7}$\beta_1^*$     &  \cellcolor[gray]{.7} 0.714   &  \cellcolor[gray]{.7} 0.914   &  \cellcolor[gray]{.7} 1.110   &  $\beta_8^*$     &   -0.003  &   0.000   &   0.016     \\
\cellcolor[gray]{.7}$\beta_2^*$     &  \cellcolor[gray]{.7} 0.803   &  \cellcolor[gray]{.7} 0.909   &  \cellcolor[gray]{.7} 1.010   &  $\beta_9^*$     &   0.000   &   0.000   &   0.024     \\
\cellcolor[gray]{.7}$\beta_3^*$     &  \cellcolor[gray]{.7} 0.789   &  \cellcolor[gray]{.7} 0.904   &  \cellcolor[gray]{.7} 1.019   &  $\beta_{10}^*$  &   0.000   &   0.000   &   0.018     \\
\cellcolor[gray]{.7}$\beta_4^*$      & \cellcolor[gray]{.7}  0.793   & \cellcolor[gray]{.7}  0.904   & \cellcolor[gray]{.7}  1.023   &    \vdots   &  \vdots    & \vdots   & \vdots    \\
\cellcolor[gray]{.7}$\beta_5^*$     &  \cellcolor[gray]{.7} 0.796   &  \cellcolor[gray]{.7} 0.907   &  \cellcolor[gray]{.7} 1.017   &   $\beta_{23}^*$  &   0.000   &   0.000   &   0.021   \\
$\beta_6^*$                         &   0.000                       &   0.000                       &   0.020                        &  $\beta_{24}^*$  &   0.000   &   0.000   &   0.016   \\
$\beta_7^*$                         &   0.000                       &   0.000                       &   0.021                        &  $\beta_{25}^*$  &   0.000   &   0.000   &   0.005   \\
% $\beta_{19}^*$  &   -0.002  &   0.000   &   0.019   \\
% $\beta_{20}^*$  &   0.000   &   0.000   &   0.011   \\
% $\beta_{17}^*$  &   0.000   &   0.000   &   0.015   \\
%$\beta_{18}^*$  &   0.000   &   0.000   &   0.013   \\
%$\beta_8^*$     &   -0.003  &   0.000   &   0.016   &   $\beta_{21}^*$  &   0.000   &   0.000   &   0.016   \\
%$\beta_9^*$     &   0.000   &   0.000   &   0.024   &   $\beta_{22}^*$  &   0.000   &   0.000   &   0.015   \\
%$\beta_{10}^*$  &   0.000   &   0.000   &   0.018   &   $\beta_{23}^*$  &   0.000   &   0.000   &   0.021   \\
%$\beta_{11}^*$  &   0.000   &   0.000   &   0.025   &   $\beta_{24}^*$  &   0.000   &   0.000   &   0.016   \\
%$\beta_{12}^*$  &   0.000   &   0.000   &   0.015   &   $\beta_{25}^*$  &   0.000   &   0.000   &   0.005   \\
%$\beta_{13}^*$  &   0.000   &   0.000   &   0.016   &       &       &       &       \\
\hline
\end{tabular}}
\end{table}
The non-pivotal procedure seems to better estimate as zeros the zero
coefficients. This is because we minimize the $\ell_1$ norm of the
coefficients without an additional $c\sigma$ term.  On the other
hand, the non-zero coefficients are better estimated using the
pivotal estimator. The non-pivotal procedure yields some shrinkage
to zero (especially for large $\sigma_*$).  Using the pivotal
procedure in the first place allows us to have a good initial guess
of $\sigma_*$.

Let us now increase $n$ to see whether we can obtain interval
estimates and take advantage of thresholding for variable selection.
We consider the same model as above and the same values of the
parameters of the method but we replace $n=49$ by $n=8000$. Then we
are no longer in a situation where we must use specific
high-dimensional techniques. However, it is still a challenging task
to select among 25 candidate variables, one of them being
endogenous. Indeed, classical selection procedures like the BIC
would require to solve $2^{25}$ least squares problems. Our methods
are much less numerically intensive. They are based on linear and
conic programming, and their computational cost is polynomial in the
dimension.
 We study both the setting with
all the 50 instruments and the setting where we estimate the linear
projection instrument.

Consider first the case where we use all the instruments. Set for
brevity
$$
\widehat{w}\triangleq\left(1-\frac{r}{\kappa_{1,\widehat{J}}^*}
-\frac{r^2} {\kappa_{1,\widehat{J}}}\right)_+^{-1}\,,\quad
w(5)\triangleq\left(1-\frac{r}{\kappa_{1}^*(5)}-\frac{r^2}
{\kappa_{1}(5)}\right)_+^{-1}\,.
$$
These are the quantities appearing in (\ref{eq:simul_conf1}) and
(\ref{eq:simul_conf2}). As above, we take $\widehat{J}$ equal to the
set of the first five coordinates; $w(5)$ corresponds to the
sparsity certificate approach with $s=5$. Computing the exact
coordinate-wise sensitivities we obtain the bound $ \widehat w
\le1.6277.$ The sparsity certificate approach with $s=5$ yields
$w(5)\le1.6306.$
%$$
%\left(1-\frac{r}{\kappa_{J_{\rm end},J(\beta^*)}}-\frac{r^2}
%{\kappa_{J_{\rm end}^c,J(\beta^*)}}\right)_+^{-1}\le1.6306.
%$$
%We also obtain $\underline{\kappa}_{1,\widehat{J}}=0.0121$,
%$\underline{\kappa}_{\infty,\widehat{J}}=0.1344$ and
%$\underline{\kappa}_{1,SC}=0.0121$,
%$\underline{\kappa}_{\infty,SC}=0.1340$, where we underline the
%sensitivities because the numbers are lower bounds.
We obtain $\widehat{\sigma}=0.2970$ and the estimates in Table
\ref{Tab5}. The values $\hat{\beta}_{l,\widehat{J}}$ and
$\hat{\beta}_{u,\widehat{J}}$ are the lower and upper confidence
limits respectively obtained from (\ref{eq:simul_conf1});
$\hat{\beta}_{l,SC}$ and $\hat{\beta}_{u,SC}$ are the lower and
upper confidence limits obtained from
(\ref{eq:simul_conf2}) (sparsity certificate approach with $s=5$). %The
%results are summarized in Table 5.
The thresholds $\omega_{k,\widehat{J}}$ and $\omega_{k}(5)$ are
computed from the formulas
$$
\omega_{k,\widehat{J}}= \frac{2\cdot
1.6277\hat{\sigma}r}{x_{k*}\kappa_{k,\widehat{J}}^*}, \qquad
\omega_{k}(5)= \frac{2\cdot
1.6306\hat{\sigma}r}{x_{k*}\kappa_{k}^*(5)}.
$$

\begin{table}%[H]
  \centering
  \caption{Confidence intervals and selection of variables, $n=8000$}\label{Tab5}
{\footnotesize \begin{tabular}{|c|c|c|c|c|c|c|c|c|c|}
  \hline
&$\hat{\beta}_{l,SC}$               & $\hat{\beta}_{l,\widehat{J}}$ & $\hat{\beta}$ & $\hat{\beta}_{u,\widehat{J}}$ & $\hat{\beta}_{u,SC}$ & $\kappa_{k,\widehat{J}}^*$ & $\kappa_{k}^*(5)$ & $\omega_{k,\widehat{J}}$ & $\omega_{k,SC}$\\
  \hline
\cellcolor[gray]{.7} $\beta_1^*$    & \cellcolor[gray]{.7} 0.131   &  \cellcolor[gray]{.7} 0.135   &  \cellcolor[gray]{.7} 1.048   & \cellcolor[gray]{.7}  1.960   & \cellcolor[gray]{.7}  1.965   &  \cellcolor[gray]{.7} 0.134   & \cellcolor[gray]{.7}  0.134   & \cellcolor[gray]{.7}  0.912   &  \cellcolor[gray]{.7} 0.917   \\
\cellcolor[gray]{.7} $\beta_2^*$    & \cellcolor[gray]{.7} 0.795   &  \cellcolor[gray]{.7} 0.804   &  \cellcolor[gray]{.7} 0.995   & \cellcolor[gray]{.7}  1.185   & \cellcolor[gray]{.7}  1.195   &  \cellcolor[gray]{.7} 0.897   & \cellcolor[gray]{.7}  0.855   & \cellcolor[gray]{.7}  0.191   &  \cellcolor[gray]{.7} 0.200   \\
\cellcolor[gray]{.7} ${\beta}_3^*$  & \cellcolor[gray]{.7} 0.824   &  \cellcolor[gray]{.7} 0.829   &  \cellcolor[gray]{.7} 1.004   & \cellcolor[gray]{.7}  1.179   & \cellcolor[gray]{.7}  1.185   &  \cellcolor[gray]{.7} 0.796   & \cellcolor[gray]{.7}  0.775   & \cellcolor[gray]{.7}  0.175   &  \cellcolor[gray]{.7} 0.180   \\
\cellcolor[gray]{.7} ${\beta}_4^*$  & \cellcolor[gray]{.7} 0.817   &  \cellcolor[gray]{.7} 0.822   &  \cellcolor[gray]{.7} 0.998   & \cellcolor[gray]{.7}  1.173   & \cellcolor[gray]{.7}  1.178   &  \cellcolor[gray]{.7} 0.858   & \cellcolor[gray]{.7}  0.833   & \cellcolor[gray]{.7}  0.175   &  \cellcolor[gray]{.7} 0.181   \\
\cellcolor[gray]{.7} ${\beta}_5^*$  & \cellcolor[gray]{.7} 0.833   &  \cellcolor[gray]{.7} 0.834   &  \cellcolor[gray]{.7} 1.001   & \cellcolor[gray]{.7}  1.168   & \cellcolor[gray]{.7}  1.168   &  \cellcolor[gray]{.7} 0.793   & \cellcolor[gray]{.7}  0.790   & \cellcolor[gray]{.7}  0.167   &  \cellcolor[gray]{.7} 0.168   \\
${\beta}_6^*$                       &  -0.163  &   -0.160  &   0.003   &   0.166   &   0.169   &   0.807   &   0.791   &   0.163   &   0.166   \\
${\beta}_7^*$                       &  -0.173  &   -0.168  &   0.002   &   0.172   &   0.177   &   0.846   &   0.823   &   0.170   &   0.175   \\
${\beta}_8^*$                       &  -0.173  &   -0.170  &   0.001   &   0.173   &   0.175   &   0.789   &   0.779   &   0.172   &   0.174   \\
\vdots                              &\vdots&\vdots&\vdots&\vdots&\vdots&\vdots&\vdots&\vdots&\vdots\\
${\beta}_{23}^*$                    & -0.190  &   -0.188  &   0.003   &   0.194   &   0.197   &   0.802   &   0.793   &   0.191   &   0.193   \\
${\beta}_{24}^*$                    & -0.171  &   -0.166  &   0.001   &   0.168   &   0.172   &   0.842   &   0.821   &   0.167   &   0.172   \\
${\beta}_{25}^*$                    & -0.172  &   -0.169  &   -0.005  &   0.158   &   0.162   &   0.828   &   0.811   &   0.163   &   0.167   \\
\hline
\end{tabular}}
\end{table}
\noindent Table 5 shows that in this example thresholding works
well: The true support of $\beta^*$ is recovered exactly by
selecting the variables, for which the estimated coefficient is
larger than the threshold. Note that the threshold for the
endogenous variable is very close to the estimate of the first
coefficient $\widehat{\beta}_1$ since the confidence intervals are
wider for the
endogenous variable. %and its true coefficient should be larger in
%absolute value than for exogenous variables in order to properly
%select the endogenous variable by thresholding.

We now consider the case where we use only 25 instruments; the 24
exogenous variables serve as their own instruments and the
Square-root Lasso estimator of the linear projection instrument is
used for the endogenous variable. This time, we apply the
Square-root Lasso with the recommended choice $c_{\sqrt{\rm
Lasso}}=1.1$. We get $\sqrt{\widehat{Q}(\widehat{\beta})}=0.3012$.
The estimates of $\widehat{\zeta}_l$ are given in Table
\ref{Tab6}.%{\bf They are all downwards biased. Why????}
\begin{table}%[H]
  \centering
  \caption{Estimates of the coefficients in the linear projection
  instrument}\label{Tab6}
{\footnotesize \begin{tabular}{|c|c|c|c|c|c|c|c|c|c|c|c|c|c|}
  \hline
$\widehat{\zeta}_1$    &   $\widehat{\zeta}_2$    &   $\widehat{\zeta}_3$    &   $\widehat{\zeta}_4$      &   $\widehat{\zeta}_5$    &   $\widehat{\zeta}_6$    &   $\widehat{\zeta}_7$    &   $\widehat{\zeta}_{8}$ & $\widehat{\zeta}_{9}$ &   $\widehat{\zeta}_{10}$ &   $\widehat{\zeta}_{11}$ &   $\widehat{\zeta}_{12}$ &   $\widehat{\zeta}_{13}$ & $\widehat{\zeta}_{14}$ \\
0.142                   &       0.145            &                      0.134 &                     0.136 &             0.137        &  0.135                   &               0.139       &       0.139           &        0.134           &          0.140            &          0.146           &            0.140        &          0.134           &       0.136             \\
\hline
$\widehat{\zeta}_{15}$ &   $\widehat{\zeta}_{16}$ &   $\widehat{\zeta}_{17}$ &   $\widehat{\zeta}_{18}$ &   $\widehat{\zeta}_{19}$ &   $\widehat{\zeta}_{20}$ &    $\widehat{\zeta}_{21}$ &    $\widehat{\zeta}_{22}$ &   $\widehat{\zeta}_{23}$ &   $\widehat{\zeta}_{24}$  &  $\widehat{\zeta}_{25}$   &    $\widehat{\zeta}_{26}$ &       &   \\
0.137                  &               0.138       &        0.141            &         0.128           &            0.142          &           0.137            &              0.133        &           0.135          &            0.135           &      0.142       &                    0.137           &              0.138 & & \\
\hline
\end{tabular}\\
We only show the non-zero coefficients (keeping only three digits).}
\end{table}
%Here we have $r/\underline{\kappa}_{1,\widehat{J}}=0.0622$ and
%$r/\underline{\kappa}_{1,SC}=0.0651$ ($r$ has a different value than
%above because we only work with $K$ instruments). This allows us to
%use \eqref{eq:t1:3}.
Next, we use (\ref{eq:simul_conf1}) and (\ref{eq:simul_conf2}) to
obtain the confidence intervals. Computing the exact coordinate-wise
sensitivities we get the bound $\widehat w\le 1.0941. $ The sparsity
certificate approach with $s=5$ yields $ w(5)\le 1.0990. $
%We also obtain $\underline{\kappa}_{1,\widehat{J}}=0.0500$,
%$\underline{\kappa}_{\infty,\widehat{J}}=0.5557$ and
%$\underline{\kappa}_{1,SC}=0.0478$,
%$\underline{\kappa}_{\infty,SC}=0.5309$.
We also get $\widehat{\sigma}=0.2970$.
\begin{table}%[H]
  \centering
  \caption{Confidence intervals and selection of variables, $n=8000$}\label{Tab7}
{\footnotesize \begin{tabular}{|c|c|c|c|c|c|c|c|c|c|}
  \hline
&$\hat{\beta}_{l,SC}$               & $\hat{\beta}_{l,\widehat{J}}$ & $\hat{\beta}$ & $\hat{\beta}_{u,\widehat{J}}$ & $\hat{\beta}_{u,SC}$ & $\kappa_{k,\widehat{J}}^*$ & $\kappa_{k}^*(5)$ & $\omega_{k,\widehat{J}}$ & $\omega_{k,SC}$\\
  \hline
\cellcolor[gray]{.7} $\beta_1^*$    & \cellcolor[gray]{.7} 0.901   & \cellcolor[gray]{.7}  0.909   & \cellcolor[gray]{.7}  1.048   & \cellcolor[gray]{.7}  1.187   & \cellcolor[gray]{.7}  1.194   & \cellcolor[gray]{.7}  0.556   & \cellcolor[gray]{.7}  0.531   & \cellcolor[gray]{.7}  0.139   & \cellcolor[gray]{.7}  0.146   \\
\cellcolor[gray]{.7} $\beta_2^*$    & \cellcolor[gray]{.7} 0.872   & \cellcolor[gray]{.7}  0.883   & \cellcolor[gray]{.7}  0.995   & \cellcolor[gray]{.7}  1.106   & \cellcolor[gray]{.7}  1.118   & \cellcolor[gray]{.7}  0.968   & \cellcolor[gray]{.7}  0.882   & \cellcolor[gray]{.7}  0.111   & \cellcolor[gray]{.7}  0.123   \\
\cellcolor[gray]{.7} ${\beta}_3^*$  & \cellcolor[gray]{.7} 0.896   & \cellcolor[gray]{.7}  0.905   & \cellcolor[gray]{.7}  1.004   & \cellcolor[gray]{.7}  1.103   & \cellcolor[gray]{.7}  1.112   & \cellcolor[gray]{.7}  0.888   & \cellcolor[gray]{.7}  0.821   & \cellcolor[gray]{.7}  0.099   & \cellcolor[gray]{.7}  0.108   \\
\cellcolor[gray]{.7} ${\beta}_4^*$  & \cellcolor[gray]{.7} 0.885   & \cellcolor[gray]{.7}  0.893   & \cellcolor[gray]{.7}  0.998   & \cellcolor[gray]{.7}  1.102   & \cellcolor[gray]{.7}  1.110   & \cellcolor[gray]{.7}  0.907   & \cellcolor[gray]{.7}  0.848   & \cellcolor[gray]{.7}  0.105   & \cellcolor[gray]{.7}  0.113   \\
\cellcolor[gray]{.7} ${\beta}_5^*$  & \cellcolor[gray]{.7} 0.899   & \cellcolor[gray]{.7}  0.902   & \cellcolor[gray]{.7}  1.001   & \cellcolor[gray]{.7}  1.100   & \cellcolor[gray]{.7}  1.103   & \cellcolor[gray]{.7}  0.843   & \cellcolor[gray]{.7}  0.823   & \cellcolor[gray]{.7}  0.099   & \cellcolor[gray]{.7}  0.102   \\
${\beta}_6^*$                       & -0.098  &   -0.092  &   0.003   &   0.099   &   0.104   &   0.868   &   0.822   &   0.095   &   0.101   \\
${\beta}_7^*$                       & -0.103  &   -0.098  &   0.002   &   0.102   &   0.107   &   0.907   &   0.869   &   0.100   &   0.105   \\
${\beta}_8^*$                       & -0.099  &   -0.095  &   0.001   &   0.098   &   0.102   &   0.886   &   0.853   &   0.096   &   0.101   \\
\vdots                              &\vdots&\vdots&\vdots&\vdots&\vdots&\vdots&\vdots&\vdots&\vdots\\
${\beta}_{23}^*$                    & -0.115  &   -0.109  &   0.003   &   0.115   &   0.121   &   0.862   &   0.825   &   0.112   &   0.118   \\
${\beta}_{24}^*$                    & -0.104  &   -0.099  &   0.001   &   0.101   &   0.106   &   0.888   &   0.848   &   0.100   &   0.105   \\
${\beta}_{25}^*$                    & -0.109  &   -0.104  &   -0.005  &   0.093   &   0.098   &   0.870   &   0.830   &   0.098   &   0.103   \\
\hline
\end{tabular}}
\end{table}
The thresholds $\omega_{k,\widehat{J}}$ and $\omega_{k}(5)$ are
obtained from the formulas
$$
\omega_{k,\widehat{J}}=\frac{2\cdot
1.0941\hat{\sigma}r}{x_{k*}\kappa_{k,\widehat{J}}^*}, \quad
\omega_{k}(5)= \frac{2\cdot
1.0990\hat{\sigma}r}{x_{k*}\kappa_{k}^*(5)}.
$$
The results are summarized in Table \ref{Tab7}. Note that the
confidence intervals and the thresholds are sharper than in the
approach including all the instruments. The particularly good news
is that the confidence interval for the coefficient of the
endogenous variable
becomes much tighter. %if we use the estimated linear projection
%instrument.

In conclusion, when the sample size is large, the coordinate-wise
sensitivities based on the sparsity certificate work remarkably well
for estimation, confidence intervals, and variable selection. We
also get a significant improvement from using the two-stage
procedure with estimated linear projection instrument.

\section{Appendix}\label{s9}

\subsection{Lower Bounds on $\kappa_{p,J}^{(c,I)}$ for
Square Matrices $\Psi_n^{(I)}$}\label{s91} In this section we consider the
case where all regressors are exogenous and we take ${\bf
D_Z^{(I)}}={\bf D_X}$. We therefore drop the exponent $I$. The
following propositions establish lower bounds on
$\kappa_{p,J}^{(c)}$ when $\Psi_n$ is a square $K\times K$ matrix.
For any $J\subseteq \{1,\dots,K\}$ we define the following
restricted eigenvalue (RE) constants
$$
\kappa_{{\rm RE},J}^{(c)} \triangleq \inf_{\Delta\in\R^K\setminus\{0\}:\
\Delta\in C_{J}^{(c)}} \frac{|\Delta^T\Psi_n\Delta|}{|\Delta_{J}|_2^2},
\quad \quad
\kappa^{'(c)}_{{\rm RE},J} \triangleq \inf_{\Delta\in\R^K\setminus\{0\}:\
\Delta\in C_{J}^{(c)}}
\frac{|J|\,|\Delta^T\Psi_n\Delta|}{|\Delta_{J}|_1^2}.
$$
%Clearly, $\kappa_{{\rm RE},J} \le \kappa'_{{\rm RE},J}$.
\begin{proposition}\label{p3} For any $J\subseteq \{1,\dots,K\}$ and $c$ in $(0,1)$
we have
$$\kappa_{1,J}^{(c)}
\ge\frac{(1-c)^2}{4|J|}\kappa^{'(c)}_{{\rm RE},J} \ge
\frac{(1-c)^2}{4|J|}\kappa_{{\rm RE},J}^{(c)}.$$
\end{proposition}
\noindent{\bf Proof.} For such that
$|\Delta_{J^c}|_1\le\frac{1+c}{1-c}|\Delta_{J}|_1$ we have
$|\Delta|_1\le\frac{2}{1-c}|\Delta_{J}|_1$.  Thus,
$$
\frac{|\Delta^T\Psi_n\Delta|}{|\Delta_J|_1^2} \le
\frac{|\Delta|_1|\Psi_n\Delta|_\infty}{|\Delta_J|_1^2} \le
\frac4{(1-c)^2}\frac{|\Psi_n\Delta|_\infty}{|\Delta|_1}\,.
$$
This proves the first inequality of the proposition.  The second
inequality is obvious. %since $|\Delta|_1^2\le|J||\Delta_{J}|_2^2$.
\hfill $\square$

\begin{proposition}\label{p1} Let $J\subseteq \{1,\dots,K\}$ and $c$ in $(0,1)$
be such that
\begin{equation}\label{gram}
\inf_{\Delta\in\R^K\setminus\{0\}:\ \Delta\in C_{J}^{(c)}}\frac{|{\bf
X}{\bf D_X}\Delta|_2}{\sqrt{n}|\Delta_{J}|_2} \ge\widetilde{\kappa}^{(c)}
\end{equation}
for some $\widetilde{\kappa}>0$, and let there exist $0<\delta<1$
such that
\begin{equation}\label{prox}
\left|\frac1n\left( {\bf X}{\bf D_X}-{\bf Z}{\bf D}_{{\bf
Z}}\right)^T {\bf X}{\bf D_X}\right|_{\infty}\le
\frac{\delta(1-c)^2(\widetilde{\kappa}^{(c)})^2}{4|J|}.
\end{equation}
Then
$$\kappa_{1,J}^{(c)}\ge\frac{(1-\delta)(1-c)^2(\widetilde{\kappa}^{(c)})^2}{4|J|}\,.$$
\end{proposition}
\noindent{\bf Proof.}  We have
\begin{align*}
|\Psi_n\Delta|_{\infty}|\Delta|_1&\ge|\Delta^T\Psi_n\Delta|\\
&\ge\left|\Delta^T\frac1n{\bf D_X}{\bf X}^T{\bf X}{\bf D_X}\Delta\right|-
\left|\Delta^T\frac1n({\bf X}{\bf D_X}-{\bf Z}{\bf D}_{{\bf Z}})^T
{\bf X}{\bf D_X}\Delta\right|
\end{align*}
where
\begin{align*}
\left|\Delta^T\frac1n\left(
{\bf X}{\bf D_X}-{\bf Z}{\bf D}_{{\bf Z}}\right)^T
{\bf X}{\bf D_X}\Delta\right|&\le\left|\frac1n\left(
{\bf X}{\bf D_X}-{\bf Z}{\bf D}_{{\bf Z}}\right)^T
{\bf X}{\bf D_X}\right|_{\infty}|\Delta|_1^2\\
&\le \frac{\alpha(1-c)^2\widetilde{\kappa}^2}{4|J|}|\Delta|_1^2.
\end{align*}
Combining these inequalities and using that
$|\Delta|_1^2\le\frac{4}{(1-c)^2}|J||\Delta_{J}|_2^2$ for all
$\Delta\in C_J^{(c)}$ (cf. proof of Proposition \ref{p1}) we get the
result.\hfill $\square$\vspace{0.3cm}

Relation \eqref{prox} accounts for the closeness between the
normalized instruments and the normalized regressors.  In the case
where there is only one endogenous regressor in the structural
equation and all exogenous regressors are used as their own
instrument, $ {\bf X}{\bf D_X}$ and ${\bf Z}{\bf D}_{{\bf Z}}$ only
differ through one column.  They are equal if there is no
endogeneity.  In that case we are left with the usual condition
\eqref{gram}. It is the restricted eigenvalue condition of Bickel,
Ritov and Tsybakov~(2009) for the Gram matrix of $X$-variables, up
to the normalization by ${\bf D_X}$.

We now obtain bounds for sensitivities $\kappa_{p,J}^{(c)}$ with $1< p\le
2$ and $c$ in $(0,1)$.  For any $s\le K$, we consider a uniform version of the
restricted eigenvalue constant: $ \kappa_{{\rm RE}}(s) \triangleq
\min_{|J|\le s} \kappa_{{\rm RE},J}^{(c)}$.

\begin{proposition}\label{p3_1} For any $s\le K/2$ and $1<
p\le 2$, we have
$$\kappa_{p,J}^{(c)}
\ge C(p) s^{-1/p}\kappa_{{\rm RE}}^{(c)}(2s), \quad \forall \ J: \ |J|\le
s,
$$ where $C(p)=2^{-1/p-1/2}(1-c)
\left(1 + \frac{1+c}{1-c}\left(p-1\right)^{-1/p}\right)^{-1}$.
\end{proposition}
\noindent{\bf Proof.} For $\Delta\in\ R^K$ and a set $J\subset
\{1,\dots,K\}$, let $J_1=J_1(\Delta,J)$ be the subset of indices in
$\{1,\dots,K\}$ corresponding to the $s$ largest in absolute value
components of $\Delta$ outside of $J$. Define $J_+=J\cup J_1$. If
$|J|\le s$ we have $|J_+|\le 2s$. It is easy to see that the $k$th
largest absolute value of elements of $\Delta_{J^c}$ satisfies
$|\Delta_{J^c}|_{(k)}\le |\Delta_{J^c}|_1/k$. Thus,
$$
|\Delta_{J^c_+}|_p^p \le |\Delta_{J^c}|_1^p\sum_{k\ge
s+1}\frac1{k^p} \le \frac{|\Delta_{J^c}|_1^p}{(p-1)s^{p-1}}\, .
$$
For $\Delta\in C_J^{(c)}$, this implies
$$
|\Delta_{J^c_+}|_p \le \frac{|\Delta_{J^c}|_1}{(p-1)^{1/p}s^{1-1/p}}
\le \frac{c_0|\Delta_{J}|_1}{(p-1)^{1/p}s^{1-1/p}} \le
\frac{c_0|\Delta_{J}|_p}{(p-1)^{1/p}}\, ,
$$
where $c_0=\frac{1+c}{1-c}$. Therefore, for $\Delta\in C_J^{(c)}$,
\begin{equation}\label{eq:prop_RE}
|\Delta|_p \le (1+c_0(p-1)^{-1/p})|\Delta_{J_+}|_p\le
(1+c_0(p-1)^{-1/p})(2s)^{1/p-1/2}|\Delta_{J_+}|_2.
\end{equation}
Using (\ref{eq:prop_RE}) and the fact that
$|\Delta|_1\le\frac{2}{1-c}|\Delta_{J}|_1 \le
\frac{2\sqrt{s}}{1-c}|\Delta_{J}|_2$ for $\Delta\in C_J^{(c)}$, we get
\begin{eqnarray*}
\frac{|\Delta^T\Psi_n\Delta|}{|\Delta_{J_+}|_2^2} &\le &
\frac{|\Delta|_1|\Psi_n\Delta|_\infty}{|\Delta_{J_+}|_2^2}
\\
&\le & \frac{2\sqrt{s}|\Psi_n\Delta|_\infty}{(1-c)|\Delta_{J_+}|_2}
\\
&\le & \frac{s^{1/p}|\Psi_n\Delta|_\infty}{C(p)|\Delta|_p}\,.
\end{eqnarray*}
Since $|J_+|\le 2s$, this proves the proposition.  \hfill
$\square$\vspace{0.3cm}

The lower bounds in Propositions \ref{p3} and \ref{p3_1} require to
control from below $|\Delta^T\Psi_n\Delta|$ (where $\Psi_n$ is a
non-symmetric possibly non-positive definite matrix) by a quadratic
form with many zero eigenvalues for vectors in a cone of dominant
coordinates. This is potentially a strong restriction on the
instruments that we can use. In other words, the sensitivity
characteristics $\kappa_{p,J}^{(c)}$ can be much larger than the above
bounds. The propositions of this section imply that, even in the
case of symmetric matrices, these characteristics are more general
and potentially lead to better results than the restricted
eigenvalues $\kappa_{\rm RE}^{(c)}(\cdot)$ appearing in the usual RE
condition of Bickel, Ritov and Tsybakov~(2009).

\subsection{Moderate Deviations for Self-normalized Sums}\label{s90}
Throughout this section $X_1,\dots, X_n$ are independent
random variables such that, for
every $i$, $\mathbb{E}[X_i]=0$.
The following result is due to Efron~(1969).
\begin{theorem}\label{th:efron}
If $X_i$ for $i=1,\hdots,n$ are symmetric, then for every $r$ positive,
$$\mathbb{P}\left(\frac{\left|\frac{1}{n}\sum_{i=1}^nX_i\right|}
{\sqrt{\frac{1}{n}\sum_{i=1}^nX_i^2}}\ge r\right)\le
2\exp\left(-\frac{nr^2}{2}\right).$$
\end{theorem}
This upper bound is refined in Pinelis~(1994) for i.i.d. random variables.
\begin{theorem}\label{th:pinelis}
If $X_i$ for $i=1,\hdots,n$ are symmetric and identically distributed,
then for every $r$ in $[0,1)$,
$$\mathbb{P}\left(\frac{\left|\frac{1}{n}\sum_{i=1}^nX_i\right|}
{\sqrt{\frac{1}{n}\sum_{i=1}^nX_i^2}}\ge r\right)\le
\frac{4e^3}{9}\Phi(-\sqrt{n}r).$$
\end{theorem}
The following result is from Jing, Shao and Wang~(2003).
\begin{theorem}\label{th:jing}
Assume that $0<\mathbb{E}[|X_i|^{2+\delta}]<\infty$ for some
$0<\delta\le1$ and
set
$$B_n^2=\sum_{i=1}^n\mathbb{E}[X_i^2],\ L_{n,\delta}=\sum_{i=1}^n
\mathbb{E}\left[|X_i|^{2+\delta}\right],\
d_{n,\delta}=B_n/L_{n,\delta}^{1/(2+\delta)}.$$  Then
$$\forall 0\le r\le \frac{d_{n,\delta}}{\sqrt{n}},\
\mathbb{P}\left(\frac{\left|\frac{1}{n}\sum_{i=1}^nX_i\right|}
{\sqrt{\frac{1}{n}\sum_{i=1}^nX_i^2}}\ge r\right)\le
2\Phi(-\sqrt{n}r)\left(1+A_0\left(\frac{1+\sqrt{n}r}{d_{n,\delta}}
\right)^{2+\delta}\right)$$ where
$A_0>0$ is an absolute constant.
\end{theorem}
Despite of its great interest to understand the
large deviations behavior of self-normalized sums,
the bound has limited practical use because $A_0$ is not an
explicit constant.\\
The following result is a corollary of Theorem 1
in Bertail, Gauth\'erat and Harari-Kermadec (2009).
\begin{theorem}\label{th:bghk}
Assume that $X_i$ for $i=1,\hdots,n$ are identically distributed and
$0<\mathbb{E}[|X_i|^{4}]<\infty$.  Then
\begin{equation}\label{emdBGHK}
\forall r\ge0,\
\mathbb{P}\left(\frac{\left|\frac{1}{n}\sum_{i=1}^nX_i\right|}
{\sqrt{\frac{1}{n}\sum_{i=1}^nX_i^2}}\ge r\right)\le
(2e+1)\exp\left(-\frac{nr^2}{2+\gamma_4r^2}\right)
\end{equation}
where $\gamma_4=\frac{\mathbb{E}[X_i^4]}{\mathbb{E}[X_i^2]^2}$, while
$$\forall r\ge\sqrt{n},\
\mathbb{P}\left(\frac{\left|\frac{1}{n}\sum_{i=1}^nX_i\right|}
{\sqrt{\frac{1}{n}\sum_{i=1}^nX_i^2}}\ge r\right)=0.$$
\end{theorem}
\noindent{\bf Proof}
Bertail, Gauth\'erat and Harari-Kermadec (2009) obtain the
upper bound for $r\ge\sqrt{n}$ and that for $0\le r<\sqrt{n}$
$$\mathbb{P}\left(\frac{\left|\frac{1}{n}\sum_{i=1}^nX_i\right|}
{\sqrt{\frac{1}{n}\sum_{i=1}^nX_i^2}}\ge r\right)\le
\inf_{a>1}\left\{2e\exp\left(-\frac{nr^2}{2(1+a)}\right)
+\exp\left(-\frac{n}{2\gamma_4}\left(1-\frac{1}{a}\right)^2\right)\right\}.$$
Because
$$\frac{1}{1+a}=\frac{1}{a}\frac{1}{1+\frac{1}{a}}
\ge\frac{1}{a}\left(1-\frac{1}{a}\right)$$
we obtain
$$-\frac{r^2}{1+a}\le-\frac{r^2}{a}\left(1-\frac{1}{a}\right).$$
This yields \eqref{emdBGHK} by choosing $a$ to equate the
two exponential terms.\hfill $\square$

\subsection{Proofs}\label{s93}
\noindent{\bf Proof of Proposition \ref{p4}.} Parts \eqref{p410},
\eqref{p411} and \eqref{p44i} of the proposition are straightforward.
The upper bound
in \eqref{k2} follows immediately from (\ref{eq:technical}).
%\eqref{k1} simply follows from the fact that
%for $1\le p\le p'\le \infty$, $\forall\Delta\in\R^K,\
%|\Delta|_{p'}\le|\Delta|_p$.\\
%Since $\Delta$ belongs to $C_{J}$,
%\begin{equation}\label{econe}
%|\Delta|_1=\left|\Delta_{J}\right|_1+\left|\Delta_{J^c}\right|_1
%\le \frac{2}{1-c}|\Delta_{J}|_1.
%\end{equation}
%The H\"older inequality also yields that
%$$\left|\Delta_{J}\right|_1\le|J|^{1-1/p}
%\left|\Delta_J\right|_p\le|J|^{1-1/p}|\Delta|_p.$$
%Thus,
%\begin{equation}\label{Hold}
%|\Delta|_1\le\frac{2}{1-c}|J|^{1-1/p}|\Delta|_p
%\end{equation}
%which gives .\\
Next, obviously,
$|\Delta|_p\le|\Delta|_1^{1/p}|\Delta|_{\infty}^{1-1/p}$ and we get
that, for $\Delta\ne0$,
$$\frac{|\Psi_n^{(I)}\Delta|_{\infty}}{|\Delta|_p}\ge
\frac{|\Psi_n^{(I)}\Delta|_{\infty}}{|\Delta|_{\infty}}
\left(\frac{|\Delta|_{\infty}}{|\Delta|_1}\right)^{1/p}.$$
Furthermore, \eqref{eq:technical} implies
$|\Delta|_1\le\frac{2}{1-c}|J||\Delta|_{\infty}$ for $\Delta \in
C_J$. Combining this with the above inequality we obtain the lower
bound in \eqref{k2}. The sequence of inequalities in \eqref{k3}
follow from the fact that
$$|J_0|{-1+1/p}|\Delta_{J_0}|_1\le |\Delta_{J_0}|_p \le |J_0|{1/p}
|\Delta_{J_0}|_1.$$
\hfill $\square$\vspace{0.3cm}

\noindent{\bf Proof of Proposition \ref{p5}.}
For all $1\le k\le K$ and $1\le l\le L$,
$$\left|\left(\Psi_n^{(I)}\Delta\right)_l-(\Psi_n^{(I)})_{lk}\Delta_k\right|\le
|\Delta|_1\max_{k'\ne k}|(\Psi_n^{(I)})_{lk'}| ,
$$
which yields
$$\left|(\Psi_n^{(I)})_{lk}\right|\left|\Delta_k\right|\le
|\Delta|_1\max_{k'\ne k}|(\Psi_n^{(I)})_{lk'}|
+\left|\left(\Psi_n^{(I)}\Delta\right)_l\right|.
$$
The two inequalities of the assumption yield
$$\left|(\Psi_n^{(I)})_{l(k)k}\right|\left|\Delta_k\right|\le
|\Delta|_1\frac{(1-\eta_2)(1-c)}{2|J|}|(\Psi_n^{(I)})_{l(k)k}|
+\frac{1-c}{\eta_1}\left|\left(\Psi_n^{(I)}\Delta\right)_{l(k)}\right|
\left|(\Psi_n^{(I)})_{l(k)k}\right|.$$ This inequality, together with the
fact that $\left|\left(\Psi_n^{(I)}\Delta\right)_{l(k)}\right|\le
\left|\Psi_n^{(I)}\Delta\right|_{\infty}$, we obtain
\begin{equation}\label{eloc}
\left|\Delta_j\right|\le
|\Delta|_1\frac{(1-\eta_2)(1-c)}{2|J|}
+\frac{1-c}{\eta_1}\left|\Psi_n^{(I)}\Delta\right|_{\infty}
\end{equation}
Summing the inequalities
over $j$ in $J$, yields
$$\left|\Delta_{J}\right|_1\le
\frac{(1-\eta_2)(1-c)}{2}|\Delta|_1+\frac{|J|(1-c)}{\eta_1}
\left|\Psi_n^{(I)}\Delta\right|_{\infty}.$$ This and the first inequality
in \eqref{eq:technical} imply that
$$\frac{1-c}{2}|\Delta|_1\le \frac{(1-\eta_2)(1-c)}{2}|\Delta|_1+\frac{|J|(1-c)}{\eta_1}
\left|\Psi_n^{(I)}\Delta\right|_{\infty}$$ which yields
$$\frac{\eta_2(1-c)}{2}|\Delta|_1\le \frac{|J|(1-c)}{\eta_1}
\left|\Psi_n^{(I)}\Delta\right|_{\infty}$$ and
\begin{equation}\label{eq:Notobvious}
\frac{\eta_1\eta_2}{2|J|}|\Delta|_1\le
\left|\Psi_n^{(I)}\Delta\right|_{\infty}.
\end{equation}
We conclude, using the
definition of the $l_1$-sensitivity, that
\begin{equation}\label{epropk1}
\kappa_{1,J}=\frac{\eta_1\eta_2}{2|J|}.
\end{equation}
Next, plugging \eqref{eq:Notobvious} into \eqref{eloc}, we deduce
\begin{align*}
\left|\Delta_j\right|&\le
\left(\frac{1-\eta_2}{\eta_1\eta_2}
+\frac{1}{\eta_1}\right)(1-c)\left|\Psi_n^{(I)}\Delta\right|_{\infty}\\
&\le\frac{1-c}{\eta_1\eta_2}\left|\Psi_n^{(I)}\Delta\right|_{\infty},
\end{align*}
which implies
$$\kappa_{\infty,J}\ge\frac{\eta_1\eta_2}{1-c}.$$
This and the lower bound in \eqref{k2} yield the result. \hfill
$\square$\vspace{0.3cm}

\noindent{\bf Proof of Theorem \ref{t1}.} Fix $\beta$ in
$\mathcal{B}_s$, denote by $u_i=y_i-x_i^T\beta$ and define the event
$$\mathcal{G}=\left\{\forall l=1,\hdots,L,\
\left({\bf D_Z^{(I)}}\right)_{ll}\left|\frac1n\sum_{i=1}^nz_{li}
u_i\right|\le \max_{l\in I}\left({\bf D_Z^{(I)}}\right)_{ll}
\sqrt{\widehat{Q}_l(\beta)}r\right\}.$$
Note that $\widehat{Q}_l(\beta)=\E_n[Z_l^2U^2]$.
We have
$$\mathcal{G}^c=\left(\bigcup_{l\in I^c}\left\{\frac{1}{n}
\left|\frac{\sum_{i=1}^nz_{li}u_i}{z_{l*}\sqrt{\E_n[U^2]}}
\right|\ge
r \right\} \right)\bigcup\left(\bigcup_{l\in I}\left\{\frac{1}{n}
\left|\frac{\sum_{i=1}^nz_{li}u_i}{\sqrt{\E_n[Z_l^2U^2]}}
\right|\ge
r \right\}\right).$$
Under scenario 1, $\mathbb{P}(\mathcal{G}^c)=\alpha$.
For the other scenarios, we use
\begin{equation}\label{eq:union}
\mathcal{G}^c\subset\bigcup_{l=1,\hdots,L}\left\{
\left|\frac{\sum_{i=1}^nz_{li}u_i}{\sqrt{\sum_{i=1}^n(z_{li}u_i)^2}}\right|\ge
\sqrt{n}r\right\}.
\end{equation}
The union bound yields
\begin{equation*}
\xP(\mathcal{G}^c)\le\sum_{l=1}^L\xP\left(
\left|\frac{\sum_{i=1}^nz_{li}u_i}{\sqrt{\sum_{i=1}^n(z_{li}u_i)^2}}\right|\ge
\sqrt{n}r\right).
\end{equation*}
We conclude, using the moderate deviations result from Section
\ref{s90} that the event $\mathcal{G}$ holds with probability at
least $1-\alpha$ (approximately at least $1-\alpha$ for scenario 4
and the two-stage procedure with scenario 5) under the respective choices of $r$.
Because the event containing $\mathcal{G}^c$ on the right hand side of \eqref{eq:union}
does not depend on $I$, we obtain statements which are uniform in $I$.
Because $\mathcal{G}^c$ does not depend on $c$ we obtain statements which are uniform in $c$
in $(0,1)$ under all 5 scenarios.
Set $\Delta\triangleq
{\bf D_X}^{-1}(\widehat{\beta}^{(c,I)} -\beta)$.  On the event $\mathcal{G}$
we have:
\begin{align}
\left|\Psi_n^{(I)}\Delta \right|_{\infty} &\le \left|\frac{1}{n}{\bf D_Z^{(I)}}
\bold{Z}^T(\bold{Y}-\bold{X}\widehat{\beta}^{(c,I)})\right|_{\infty}
+\left|\frac{1}{n}{\bf D_Z^{(I)}}\bold{Z}^T(\bold{Y}-\bold{X}\beta)
\right|_{\infty}\label{stop1}\\
&\le r\widehat{\sigma}^{(c,I)}+\left|\frac1n{\bf D_Z^{(I)}}{\bf Z}^T
{\bf U}\right|_{\infty}\label{eq:Notobvious2}\\
&\le
r\left(\widehat{\sigma}^{(c,I)}+\max_{l\in I}({\bf D_{Z}})_{ll}\sqrt{
\widehat{Q}_l(\beta)}\right).\nonumber\label{stop2}
\end{align}
The inequality \eqref{eq:Notobvious2} holds because
$(\widehat{\beta}^{(c,I)},\widehat{\sigma}^{(c,I)})$ belongs to the
set $\widehat{\mathcal{I}}^{(I)}$ by definition.  Notice that, on
the event $\mathcal{G}$, the pair
$\left(\beta,\max_{l=1,\hdots,L}({\bf
D_Z^{(I)}})_{ll}\sqrt{\widehat{Q}_l(\beta)}\right)$ belongs to the
set $\widehat{\mathcal{I}}^{(I)}$. On the other hand,
$(\widehat{\beta}^{(c,I)},\widehat{\sigma}^{(c,I)})$ minimizes the
criterion $\left|{\bf D_X}^{-1}\beta\right|_1 +c\sigma$ on the same
set $\widehat{\mathcal{I}}^{(I)}$. Thus, on the event $\mathcal{G}$,
\begin{equation}\label{eq:main}
\left|{\bf D_X}^{-1}\widehat{\beta}^{(c,I)}\right|_1 +c\widehat{\sigma}^{(c,I)}\le
|{\bf D_X}^{-1}\beta|_1+c\max_{l\in I}({\bf D_{Z}})_{ll}\sqrt{\widehat{Q}_l(\beta)}.
\end{equation}
This implies, again on the event $\mathcal{G}$,
\begin{align}
\left|\Delta_{J(\beta)^c}\right|_1&=\sum_{k\in J(\beta)^c}
\left|\mathbb{E}_n[X_k^2]^{1/2}\widehat{\beta}_k^{(c,I)}\right|\label{in0}\\
&\le\sum_{k\in J(\beta)}\left(
\left|\mathbb{E}_n[X_k^2]^{1/2}\beta_k\right|
-\left|\mathbb{E}_n[X_k^2]^{1/2}\widehat{\beta}_k^{(c,I)}\right|\right)\nonumber\\
&\quad +c\left(\max_{l\in I}({\bf D_Z^{(I)}})_{ll}\sqrt{\widehat{Q}_l(\beta)}-\max_{l\in I}({\bf D_Z^{(I)}})_{ll}\sqrt{\widehat{Q}_l(\widehat\beta^{(c,I)})}\right).\nonumber
\end{align}
For every $l$ in $I$, $\gamma\mapsto({\bf
D_Z^{(I)}})_{ll}\sqrt{\widehat{Q}_l(\gamma)}$ is almost surely
differentiable.  Thus almost surely its subgradient is the gradient
and thus the singleton
$$\partial \left(({\bf D_Z^{(I)}})_{ll}\sqrt{\widehat{Q}_l(\cdot)}\right)(\beta)=
\nabla\left(({\bf D_Z^{(I)}})_{ll}\sqrt{\widehat{Q}_l(\cdot)}\right)(\beta)=\frac{\mathbb{E}_n[U X^TZ_l^2]} {\sqrt{\E_n[(UZ_l)^2]}}.$$
Because the function is also convex, almost surely,
\begin{align*}({\bf D_Z^{(I)}})_{ll}\sqrt{\widehat{Q}_l(\beta)}-({\bf D_Z^{(I)}})_{ll}\sqrt{\widehat{Q}_l(\widehat{\beta}^{(c,I)})}&\le \nabla\left(({\bf D_Z^{(I)}})_{ll}\sqrt{\widehat{Q}_l(\cdot)}\right)(\beta)^T(\beta-\widehat{\beta}^{(c,I)})\\
&\le({\bf D_Z^{(I)}})_{ll}\left(\frac{\mathbb{E}_n[U X^TZ_l^2]} {\sqrt{\E_n[(UZ_l)^2]}}\right)^T{\bf
D_X}{\bf
D_X}^{-1}(\beta-\widehat{\beta}^{(c,I)})\\
&\le-({\bf D_Z^{(I)}})_{ll}\frac{\mathbb{E}_n[U X^TZ_l^2]} {\sqrt{\E_n[(UZ_l)^2]}}^T{\bf
D_X}\Delta.
\end{align*}
By the Cauchy-Schwartz inequality,
$$\left|-({\bf D_Z^{(I)}})_{ll}\frac{\mathbb{E}_n[U X^TZ_l^2]} {\sqrt{\E_n[(UZ_l)^2]}}^T{\bf
D_X}\right|\le 1.$$ We will now use
the same notation for $\mathcal{G}$ and $\mathcal{G}$ intersected
with the probability 1 event where each mapping is differentiable
and work on the later event. By the Dubovitsky-Milutin theorem (see,
{\em e.g.}, Alekseev, Tikhomirov and Fomin (1987), Chapter 2)
$$\left(\max_{l\in I}({\bf D_Z^{(I)}})_{ll}\sqrt{\widehat{Q}_l(\beta)}
-\max_{l\in I}({\bf D_Z^{(I)}})_{ll}\sqrt{\widehat{Q}_l(\widehat\beta^{(c,I)})}\right)
\le|\Delta|_1.$$
Thus (\ref{in0}) implies
\begin{equation}\label{eq:conec}
\left|\Delta_{J(\beta)^c}\right|_1\le
\frac{1+c}{1-c}\left|\Delta_{J(\beta)}\right|_1.
\end{equation}
Thus, $\Delta\in C_{J(\beta)}^{(c,I)}$ on the event $\mathcal{G}$. Using
\eqref{stop1} and arguing as in (\ref{in0}) we find
\begin{align}
\left|\Psi_n^{(I)}\Delta \right|_{\infty}&\le r\left(2\widehat{\sigma}^{(c,I)}+
\max_{l\in I}({\bf D_Z^{(I)}})_{ll}\sqrt{\widehat{Q}_l(\beta)}-\widehat{\sigma}^{(c,I)}\right)\label{in1}\\
&\le
r\left(2\widehat{\sigma}^{(c,I)}+\max_{l\in I}({\bf D_Z^{(I)}})_{ll}\sqrt{\widehat{Q}_l(\beta)}
-\max_{l\in I}({\bf D_Z^{(I)}})_{ll}\sqrt{\widehat{Q}_l
(\hat\beta^{(c,I)})}\right)\nonumber\\
&\le r\left(2\widehat{\sigma}^{(c,I)}+ |\Delta|_1\right).
\nonumber
\end{align}
This again uses the convexity for every $l$ in $I$, of
$\gamma\mapsto({\bf D_Z^{(I)}})_{ll}\sqrt{\widehat{Q}_l(\gamma)}$.
Using the definition of the sensitivities we get that, on the event
${\mathcal G}$,
$$
\left|\Psi_n^{(I)}\Delta \right|_{\infty}\le r\left(2\widehat{\sigma}^{(c,I)}
+\frac{\left|\Psi_n^{(I)}\Delta
\right|_{\infty}}{\kappa_{1,J(\beta)}^{(c,I)}}\right),$$
which implies
\begin{equation}\label{eq:ubc}
\left|\Psi_n^{(I)}\Delta \right|_{\infty}\le2\widehat{\sigma}^{(c,I)}r
\left(1-\frac{r}{\kappa_{1,J(\beta)}^{(c,I)}}\right)_+^{-1}.
\end{equation}
This inequality and the definition of the sensitivities yield
\eqref{eq:t1:1} and \eqref{eq:t1:3}.

To prove \eqref{eq:t1:4}, it suffices to note that, by
\eqref{eq:main} and by the definition of
$\kappa_{1,J(\beta),J(\beta)}^{(c,I)}$,
\begin{align*}
c\widehat{\sigma}^{(c,I)}&\le|\Delta_{J(\beta)}|_1+c
\max_{l\in I}({\bf D_Z^{(I)}})_{ll}\sqrt{\widehat{Q}_l(\beta)}\\
&\le \quad
\frac{|\Psi_n^{(I)}\Delta|_{\infty}}{\kappa_{1,J(\beta),J(\beta)}^{(c,I)}}
+c\max_{l\in I}({\bf D_Z^{(I)}})_{ll}\sqrt{\widehat{Q}_l(\beta)},
\end{align*}
and to combine this inequality with \eqref{stop1}.\hfill
$\square$\vspace{0.3cm}
%%%%%%%%%%%%%%%%%%%%%%%%%%%%

\noindent{\bf Proof of Theorem \ref{tapproxsparse}.} Take $\beta$ in
$\mathcal{B}_s$.  Fix an arbitrary subset $J$ of $\{1,\hdots,K\}$.
Acting as in \eqref{in0} with $J$ instead of $J(\beta)$, we get:
\begin{align*}
\sum_{k\in J^c}
\left|\mathbb{E}_n[X_k^2]^{1/2}\widehat{\beta}_k^{(c,I)}\right|+
\sum_{k\in J^c}
\left|\mathbb{E}_n[X_k^2]^{1/2}\beta_k\right|
&\le\sum_{k\in J}\left(
\left|\mathbb{E}_n[X_k^2]^{1/2}\beta_k\right|
-\left|\mathbb{E}_n[X_k^2]^{1/2}\widehat{\beta}_k^{(c,I)}\right|\right)
+2\sum_{k\in J^c}
\left|\mathbb{E}_n[X_k^2]^{1/2}\beta_k\right|\\
&\quad+c\left(\sqrt{\widehat{Q}(\beta)}-\sqrt{\widehat{Q}(\widehat
\beta^{(c,I)})}\right)\\
&\le\left|\Delta_{J}\right|_1+2\left|\left({\bf D_X}^{-1}
\beta\right)_{J^c}\right|_1+c|\Delta|_1.
\end{align*}
This yields
\begin{equation}\label{eminL1}
|\Delta_{J^c}|_1\le|\Delta_J|_1+2\left|\left({\bf D_X}^{-1}\beta
\right)_{J^c}\right|_1+c|\Delta|_1.
\end{equation}
Assume now that we are on the event $\mathcal{G}$. Consider the two
possible cases. First, if $2\left|\left({\bf
D_X}^{-1}\beta\right)_{J^c}\right|_1\le|\Delta_J|_1$, then
$\Delta\in\widetilde{C}_{J}^{(c,I)}$.  From this, using the definition of
the sensitivity $\widetilde{\kappa}_{p,J_0,J}^{(c,I)}$, we get that
$|\Delta_{J_0}|_p$ is bounded from above by the first term of the
maximum in (\ref{eq:tapproxsparse}).  Second, if $2\left|\left({\bf
D_X}^{-1}\beta\right)_{J^c}\right|_1>|\Delta_J|_1$, then for any
$p\in[1,\infty]$ we have a simple bound
$$ |\Delta_{J_0}|_p \le |\Delta|_1
=|\Delta_{J^c}|_1+|\Delta_J|_1\le \frac{6}{1-c}\left|\left({\bf
D_X}^{-1}\beta\right)_{J^c}\right|_1.
$$
In conclusion, $|\Delta|_p$ is smaller than the maximum of the two
bounds.\hfill $\square$\vspace{0.3cm}

%%%%%%%%%%%%%%%%%%%%%%%%%%%%

\noindent{\bf Proof of Theorem \ref{t1b}.} Part (i) of the theorem
is a consequence of \eqref{eq:t1:4} and Assumptions \ref{ass1b} and
\ref{ass1c}. Parts (ii) and (iii) follow immediately from
\eqref{eq:t1:1}, \eqref{eq:t1:3}, and Assumptions \ref{ass1b} and
\ref{ass1c}. Part (iv) is straightforward in view of
\eqref{eq:th1b:3}.\hfill $\square$\vspace{0.3cm}

%%%%%%%%%%%%%%%%%%%%%%%%

\noindent{\bf Proof of Theorem \ref{th:threshold}.} Fix $\beta$ in
$\mathcal{B}_s$.  Let ${\mathcal G}_j$ be the events of
probabilities at least $1-\gamma_j$ respectively appearing in
Assumptions \ref{ass1b}, \ref{ass1X}, \ref{ass_select}. Assume that
all these events hold, as well as the event ${\mathcal G}$.  Then,
using Theorem \ref{t1b} (i),
$$
\omega_{k}^{(c,I)}(s) \le \frac{2\sigma_* r }{c_{k}^{(c)*}(s)v_k
}\left(1+\frac{r}{c c_{1,J(\beta)}^{(c,I)}}\right)
\left(1-\frac{r}{c
c_{1,J(\beta)}^{(c,I)}}\right)_+^{-1}\left(1-\frac{r}{c_{1}^{(c,I)}(s)}\right)_+^{-1}
\triangleq \omega_k^{(c)*}.
$$
By assumption, $|\beta_k|> 2\omega_k^{(c)*}$ for $k\in J(\beta)$. Note
that the following two cases can occur. First, if $k\in J(\beta)^c$
(so that $\beta_k=0$) then, using \eqref{eq:t1:3} and Assumptions
\ref{ass1b} and \ref{ass_select}, we obtain
$|\widehat{\beta}_k^{(c,I)}|\le\omega_k^{(c,I)}$, which implies
$\widetilde{\beta}_k^{(c,I)}=0$. Second, if $k\in J(\beta)$, then using
again \eqref{eq:t1:3} we get
$||\beta_k|-|\widehat{\beta}_k^{(c,I)}||\le|\beta_k-
\widehat{\beta}_k^{(c,I)}|\le\omega_k^{(c,I)}\le \omega_k^{(c)*}$. Since $|\beta_k|>
2\omega_k^{(c)*}$ for $k\in J(\beta)$, we obtain that
$|\widehat{\beta}_k^{(c,I)}|>\omega_k^{(c,I)} $, so that
$\widetilde{\beta}_k^{(c,I)}=\widehat{\beta}_k^{(c,I)}$ and the signs of $\beta_k$
and $\widehat{\beta}_k^{(c,I)}$ coincide. This yields the result. \hfill
$\square$\vspace{0.3cm}

\noindent{\bf Proof of Theorem \ref{t1HT}.} Fix $\beta$ in
$\mathcal{B}_s$ and define the event
$$\mathcal{G}=\left\{\left|\frac1n{\bf D_Z}{\bf Z}^T
{\bf U}\right|_{\infty}\le \sqrt{\widehat{Q}(\beta)}r\right\} .$$
Note that $\widehat{Q}(\beta)=\E_n[U^2]$.  Under scenario 1, the
event $\mathcal{G}$ holds with probability at least
$1-\alpha$. For the other scenarios, the union bound yields
\begin{align}\label{eq:union}
\xP(\mathcal{G}^c)&\le \sum_{l=1}^L\xP\left(\frac{1}{n}
\left|\frac{\sum_{i=1}^nz_{li}u_i}{z_{l*}\sqrt{\E_n[U^2]}}
\right|\ge
r \,\right)\\
&\le \sum_{l=1}^L\xP\left(
\left|\frac{\sum_{i=1}^nz_{li}u_i}{\sqrt{\sum_{i=1}^n(z_{li}u_i)^2}}\right|\ge
\sqrt{n}r\right).\nonumber
\end{align}
We conclude, using the moderate deviations result from Section
\ref{s90} that the event $\mathcal{G}$ holds with probability at
least $1-\alpha$ (approximately at least $1-\alpha$ for scenario 4
and the two-stage procedure with scenario 5).

Take $c$ in $(0,1)$, set $\Delta\triangleq {\bf D_X}^{-1}(\widehat{\beta}^{(c)}-\beta)$. On
the event $\mathcal{G}$ we have:
\begin{align}\label{stop1HT}
\left|\Psi_n\Delta \right|_{\infty} &\le \left|\frac{1}{n}{\bf D_Z}\bold{Z}^T(\bold{Y}-\bold{X}\widehat{\beta}^{(c)})\right|_{\infty}
+\left|\frac{1}{n}{\bf D_Z}\bold{Z}^T(\bold{Y}-\bold{X}\beta)
\right|_{\infty}\\
&\le r\widehat{\sigma}^{(c)}+\left|\frac1n{\bf D_Z}{\bf Z}^T
{\bf U}\right|_{\infty}\nonumber\\
&\le
r\left(\widehat{\sigma}^{(c)}+\sqrt{\widehat{Q}(\beta)}\right).\nonumber\label{stop2HT}
\end{align}
Notice that, on the event $\mathcal{G}$, the pair
$\left(\beta,\sqrt{\widehat{Q}(\beta)}\right)$ belongs to the set
$\widehat{\mathcal{I}}$. On the other hand,
$(\widehat{\beta}^{(c)},\widehat{\sigma}^{(c)})$ minimizes the
criterion $\left|{\bf D_X}^{-1}\beta\right|_1 +c\sigma$ on the same
set $\widehat{\mathcal{I}}$. Thus, on the event $\mathcal{G}$,
\begin{equation}\label{eq:mainHT}
\left|{\bf D_X}^{-1}\widehat{\beta}^{(c)}\right|_1 +c\widehat{\sigma}^{(c)}\le
|{\bf D_X}^{-1}\beta|_1+c\sqrt{\widehat{Q}(\beta)}.
\end{equation}
This implies, again on the event $\mathcal{G}$,
\begin{align}
\left|\Delta_{J(\beta)^c}\right|_1&=\sum_{k\in J(\beta)^c}
\left|x_{k*}\widehat{\beta}_k^{(c)}\right|\label{in0HT}\\
&\le\sum_{k\in J(\beta)}\left(
\left|x_{k*}\beta_k\right|
-\left|x_{k*}\widehat{\beta}_k^{(c)}\right|\right)+c\left(\sqrt{\widehat{Q}(\beta)}-\sqrt{\widehat{Q}(\widehat\beta^{(c)})}\right)\nonumber\\
&\le\left|\Delta_{J(\beta)}\right|_1+c\left(\sqrt{\widehat{Q}(\beta)}-\sqrt{\widehat{Q}(\widehat\beta^{(c)})}\right)\nonumber\\
&\le\left|\Delta_{J(\beta)}\right|_1 +c\left|\frac{\E_n[U X^T]{\bf
D_X}\Delta} {\sqrt{\E_n[U^2]}}\right| \quad \Big(\mbox{by convexity
of}\ \beta\mapsto \sqrt{\widehat{Q}(\beta)}\Big)
\nonumber\\
&\le\left|\Delta_{J(\beta)}\right|_1+ c\left|\frac{\E_n[U X^T]{\bf
D_X}} {\sqrt{\E_n[U^2]}}\right|_\infty|\Delta|_1
\nonumber\\
&\le\left|\Delta_{J(\beta)}\right|_1+ c|\Delta|_1\quad \mbox{(by
the Cauchy-Schwarz inequality).}\nonumber
\end{align}
Note that (\ref{in0HT}) can be re-written as a cone condition:
\begin{equation}\label{eq:coneHT}
\left|\Delta_{J(\beta)^c}\right|_1\le
\frac{1+c}{1-c}\left|\Delta_{J(\beta)}\right|_1.
\end{equation}
Thus, $\Delta\in C_{J(\beta)}^{(c)}$ on the event $\mathcal{G}$. Using
\eqref{stop1} and arguing as in (\ref{in0}) we find
\begin{align}
\left|\Psi_n\Delta \right|_{\infty}&\le r\left(2\widehat{\sigma}^{(c)}+
\sqrt{\widehat{Q}(\beta)}-\widehat{\sigma}^{(c)}\right)\label{in1HT}\\
&\le
r\left(2\widehat{\sigma}^{(c)}+\sqrt{\widehat{Q}(\beta)}-\sqrt{\widehat{Q}(\widehat{\beta}^{(c)})}\right)
\quad (\mbox{since\ }\sqrt{\widehat{Q}(\widehat{\beta}^{(c)})}\le\widehat{\sigma}^{(c)})\nonumber\\
&\le r\left(2\widehat{\sigma}^{(c)}+ \left|\frac{\E_n[U X^T]{\bf
D_X}\Delta} {\sqrt{\E_n[U^2]}}\right|
\right)
\nonumber
%\label{in2}%\quad \Big(\mbox{by convexity of}\
%\beta\mapsto \sqrt{\widehat{Q}(\beta)}\Big)
\\
&\le r\left(2\widehat{\sigma}^{(c)}+\max_{j\in J_{\rm exo}^c
}\left|\frac{\E_n[U X_j]} {\sqrt{\E_n[X_j^2U^2]}}\right|
|\Delta_{J_{\rm exo}^c}|_1 + \max_{j\in J_{\rm exo}
}\left|\frac{\E_n[U X_j]} {\sqrt{\E_n[X_j^2U^2]}}\right|
|\Delta_{J_{\rm exo}}|_1\right)\nonumber
\\
&\le r\left(2\widehat{\sigma}^{(c)}+ |\Delta_{J_{\rm exo}^c}|_1 + \max_{j\in
J_{\rm exo} }\left|\frac{\E_n[U X_j]}
{\sqrt{\E_n[X_j^2U^2]}}\right| |\Delta_{J_{\rm exo}}|_1\right) \ \
\mbox{(by the Cauchy-Schwarz inequality).}\nonumber
\end{align}
Since the exogenous variables serve as their own instruments, we
obtain that, on the event ${\mathcal G}$,
$$
\max_{j\in J_{\rm exo} }\left|\frac{\E_n[U X_j]}
{\sqrt{\E_n[X_j^2U^2]}}\right| \le r.
$$
Combining this with (\ref{in1HT}) and using the definition of the
block sensitivity $\kappa_{1,J_0,J(\beta)}^{(c)}$ with $J_0=J_{{\rm
exo}}^c$, $J_0=J_{{\rm exo}}$, we get that, on the event ${\mathcal
G}$,
\begin{align}
\left|\Psi_n\Delta \right|_{\infty}&\le r\left(2\widehat{\sigma}^{(c)}+
\sqrt{\widehat{Q}(\beta)}-\widehat{\sigma}^{(c)}\right)\label{in3HT}\\
&\le r\left(2\widehat{\sigma}^{(c)}+\frac{\left|\Psi_n\Delta
\right|_{\infty}}{\kappa_{1,J_{{\rm exo}}^c,J(\beta)}^{(c)}}
+r\frac{\left|\Psi_n\Delta \right|_{\infty}}{\kappa_{1,J_{{\rm
exo}},J(\beta)}^{(c)}}\right)\,, \nonumber
\end{align}
which implies
\begin{equation}\label{eq:ubHT}
\left|\Psi_n\Delta \right|_{\infty}\le2\widehat{\sigma}^{(c)}r\left(1-\frac{r}{\kappa_{1,J_{{\rm exo}}^c,J(\beta)}^{(c)}}-\frac{r^2}{\kappa_{1,J_{{\rm exo}},J(\beta)}^{(c)}}\right)_+^{-1}.
\end{equation}
This inequality and the definition of the sensitivities yield
\eqref{eq:t1HT:1} and \eqref{eq:t1HT:3}.

To prove \eqref{eq:t1HT:4}, it suffices to note that, by
\eqref{eq:mainHT} and by the definition of
$\kappa_{1,J(\beta),J(\beta)}^{(c)}$,
\begin{align*}
c\widehat{\sigma}^{(c)}&\le|\Delta_{J(\beta)}|_1+c\sqrt{\widehat{Q}(\beta)}\\
&\le \quad
\frac{|\Psi_n\Delta|_{\infty}}{\kappa_{1,J(\beta),J(\beta)}^{(c)}}
+c\sqrt{\widehat{Q}(\beta)},
\end{align*}
and to combine this inequality with \eqref{stop1HT}.\hfill
$\square$\vspace{0.3cm}

\noindent{\bf Proof of Theorem \ref{t12S}.}
Take $\beta$ in $\mathcal{B}_s$.
Note that
\begin{align*}
|z_i^T\widehat\zeta^{(c_{RF})}|&\ge |z_i^T\zeta|-\left|\left({\bf D_Z}\right)^{-1}(\widehat{\zeta}^{(c_{RF})}-\zeta)\right|_1\\
&\ge |z_i^T\widehat\zeta^{(c_{RF})}|-C_1(r,c_{RF},s_{RF})
\end{align*}
thus
$$\max_{i=1,\hdots,n}|z_i^T\widehat\zeta^{(c_{RF})}|+C_1(r,c_{RF},s_{RF})r\ge \max_{i=1,\hdots,n}|z_i^T\zeta|.$$
As well, we have on $E_{\alpha}$,
\begin{align*}
\frac1n|\widehat{\zeta}^T{\bf Z}^T{\bf U}|&\le \left|\left({\bf D_Z}\right)^{-1}(\widehat{\zeta}^{(c_{RF})}-\zeta)\right|_1\sqrt{\widehat{Q}(\beta)}r
+\frac1n|\zeta^{T}{\bf Z}^T{\bf U}| \\
&\le C_1(r,c_{RF},s_{RF})\sqrt{\widehat{Q}(\beta)}r+\frac1n|\zeta^{T}{\bf Z}^T{\bf U}|,
\end{align*}
thus
\begin{align*}
\frac{\frac1n|\widehat{\zeta}^{(c_{RF})T}{\bf Z}^T{\bf U}|}{\max_{i=1,\hdots,n}|z_i^T\widehat\zeta^{(c_{RF})}|+C_1(r,c_{RF},s_{RF})}&\le \frac{C_1(r,c_{RF},s_{RF})\sqrt{\widehat{Q}(\beta)}r}{\max_{i=1,\hdots,n}|z_i^T\widehat\zeta^{(c_{RF})}|+C_1(r,c_{RF},s_{RF})}+\frac{\frac1n|\zeta^{T}{\bf Z}^T{\bf U}|}{\max_{i=1,\hdots,n}|z_i^T\zeta|}\\
&\le \frac{C_1(r, c_{RF} ,s_{RF})\sqrt{\widehat{Q}(\beta)}r}{\max_{i=1,\hdots,n}|z_i^T\widehat\zeta|+C_1(r, c_{RF}, s_{RF})}+\sqrt{\widehat{Q}(\beta)}r\\
&\le \sqrt{\widehat{Q}(\beta)}r\left(1+\frac{C_1(r, c_{RF},  s_{RF})}{\max_{i=1,\hdots,n}|z_i^T\widehat\zeta^{(c_{RF})}|+C_1(r , c_{RF} ,s_{RF})}\right).
\end{align*}
This yields
$$\frac{\frac1n|\widehat{\zeta}^{(c_{RF})T}{\bf Z}^T{\bf U}|}{\max_{i=1,\hdots,n}|z_i^T\widehat\zeta^{(c_{RF})}|+2C_1(r, c_{RF}, s_{RF})}\le\sqrt{\widehat{Q}(\beta)}r.$$
The rest of the proof is the same as for Theorem \ref{t1HT}.\hfill
$\square$

%%%%%%%%%%%%%%%%%%%%%%%

\noindent{\bf Proof of Theorem \ref{th:nonvalid}.} Throughout the
proof, we assume that we are on the event of probability at least
$1-\alpha_2$ where (\ref{eq:nv1}) holds. It follows easily from
(\ref{eq:nv1}) that
\begin{equation}\label{eq:proof:nv1}
\left|\frac1n{\bf \overline{Z}}^T{\bf
X}(\widehat{\beta}-\beta^*)\right|_\infty \le\widehat{b}{\overline
z}_{*}.
\end{equation}
Next, an argument similar to (\ref{eq:union}) and
Theorem~\ref{th:jing} yield that, with probability at least
$1-\alpha_1$,
\begin{equation}\label{eq:proof:nv2}
\left|\frac1n{\bf \overline{Z}}^T{\bf U}-\theta^*\right|_\infty \le
r_1 \max_{l=1,\dots,L_1}\sqrt{\frac1n \sum_{i=1}^n
(\overline{z}_{li}u_i-\theta^*_l)^2} = r_1 F(\theta^*,\beta^*).
\end{equation}
In what follows, we assume that we are on the event of probability
at least $1-\alpha_1-\alpha_2$ where both (\ref{eq:proof:nv1}) and
(\ref{eq:proof:nv2}) are satisfied.

We will use the properties of $F(\theta,\beta)$ stated in the next
lemma that we prove in Section \ref{subsec:proof:lem}.
\begin{lemma}\label{lem:Ftheta} We have
\begin{eqnarray}\label{eq:proof:nv3}
F(\theta^*,\widehat{\beta}) - F(\widehat{\theta},\widehat{\beta})
&\le& |\widehat{\theta}-\theta^*|_1,\\
|F(\theta^*,\widehat{\beta}) - F(\theta^*,\beta^*)| &\le& {\overline
z}_{*}\left|{\bf D_X}^{-1}(\widehat{\beta}-\beta^*)\right|_1
\le\widehat{b}{\overline z}_{*} .\label{eq:proof:nv4}
\end{eqnarray}
\end{lemma}
We proceed now to the proof of Theorem~\ref{th:nonvalid}. First, we
show that the pair $(\theta, \sigma_1) = (\theta^*,
F(\theta^*,\beta^*))$ belongs to the set $\widehat{\mathcal I}_1$.
Indeed, from (\ref{eq:proof:nv1}) and (\ref{eq:proof:nv2}) we get
\begin{eqnarray*}\label{eq:proof:nv}
\left|\frac1n{\bf \overline{Z}}^T({\bf Y}-{\bf
X}\widehat{\beta})-\theta^*\right|_\infty &\le& \left|\frac1n{\bf
\overline{Z}}^T{\bf U}-\theta^*\right|_\infty +\left|\frac1n{\bf
\overline{Z}}^T{\bf
X}(\widehat{\beta}-\beta^*)\right|_\infty\\
&\le&r_1 F(\theta^*,\beta^*) +\widehat{b}{\overline z}_{*} .%\label{eq:proof:nv}
\end{eqnarray*}
Thus, the pair $(\theta, \sigma_1) = (\theta^*,
F(\theta^*,\beta^*))$ satisfies the first constraint in the
definition of $\widehat{\mathcal I}_1$. It satisfies the second
constraint as well, since $F(\theta^*,\widehat{\beta}) \le
F(\theta^*,\beta^*) + \widehat{b}{\overline z}_{*}$ by
(\ref{eq:proof:nv4}).

Take any $c$ in $(0,1)$, as $(\theta^*, F(\theta^*,\beta^*))\in \widehat{\mathcal I}_1$
and $(\widehat{\theta}^{(c)},\widehat{\sigma}_1^{(c)})$ minimizes $|\theta|_1 +
c{\sigma}_1$ over $\widehat{\mathcal I}_1$, we have
\begin{equation}\label{eq:proof:nv5}
|\widehat{\theta}^{(c)}|_1 + c\widehat{\sigma}_1^{(c)} \le |\theta^*|_1 +
cF(\theta^*,\beta^*),
\end{equation}
which implies
\begin{equation}\label{eq:proof:nv6}
|\overline{\Delta}_{J(\theta^*)^c}|_1 \le
|\overline{\Delta}_{J(\theta^*)}|_1 +
c(F(\theta^*,\beta^*)-\widehat{\sigma}_1^{(c)}),
\end{equation}
where $\overline{\Delta}=\widehat{\theta}^{(c)}-\theta^*$. Using the fact
that $F(\widehat{\theta}^{(c)},\widehat{\beta})\le \widehat{\sigma}_1^{(c)} +
\widehat{b}{\overline z}_{*}$, (\ref{eq:proof:nv3}), and
(\ref{eq:proof:nv4}) we obtain
\begin{eqnarray}\label{eq:proof:nv6a}
F(\theta^*,\beta^*)-\widehat{\sigma}_1^{(c)} &\le&
F(\theta^*,\beta^*)-F(\widehat{\theta}^{(c)} ,\widehat{\beta})+\widehat{b}{\overline z}_{*}\\
&\le& |\widehat{\theta}^{(c)} -\theta^*|_1 + 2\widehat{b}{\overline
z}_{*}.\nonumber
\end{eqnarray}
This inequality and (\ref{eq:proof:nv6}) yield
$$
|\overline{\Delta}_{J(\theta^*)^c}|_1 \le
|\overline{\Delta}_{J(\theta^*)}|_1 + c|\widehat{\theta}^{(c)}-\theta^*|_1
+ 2c\widehat{b}{\overline z}_{*},
$$
or equivalently,
\begin{equation}\label{eq:proof:nv7}
|\overline{\Delta}_{J(\theta^*)^c}|_1 \le \frac{1+c}{1-c}\,
|\overline{\Delta}_{J(\theta^*)}|_1 +
\frac{2c}{1-c}\widehat{b}{\overline z}_{*}.
\end{equation}
Next, using (\ref{eq:proof:nv1}), (\ref{eq:proof:nv2}) and the
second constraint in the definition of
$(\widehat{\theta}^{(c)},\widehat{\sigma}_1^{(c)})$, we find
\begin{eqnarray*}%\label{eq:proof:nv}
|\widehat{\theta}^{(c)}-\theta^*|_\infty &\le& \left|\frac1n{\bf
\overline{Z}}^T({\bf Y}-{\bf
X}\widehat{\beta})-\widehat{\theta}^{(c)}\,\right|_\infty\\
\nonumber &&  + \ \left|\frac1n{\bf \overline{Z}}^T{\bf
U}-\theta^*\right|_\infty +\left|\frac1n{\bf \overline{Z}}^T{\bf
X}(\widehat{\beta}-\beta^*)\right|_\infty\\
&\le& r_1(\widehat{\sigma}_1^{(c)} +  F(\theta^*,\beta^*))
+2\widehat{b}{\overline z}_{*}.\nonumber
\end{eqnarray*}
This and (\ref{eq:proof:nv6a}) yield
\begin{eqnarray}\label{eq:proof:nv8}
|\widehat{\theta}^{(c)}-\theta^*|_\infty &\le& r_1(2\widehat{\sigma}_1^{(c)} +
|\widehat{\theta}^{(c)}-\theta^*|_1) +2(1+r_1)\widehat{b}{\overline
z}_{*}.
\end{eqnarray}
On the other hand, (\ref{eq:proof:nv7}) implies
\begin{eqnarray}\label{eq:proof:nv9}
|\widehat{\theta}^{(c)}-\theta^*|_1 &\le& \frac{2}{1-c}\,
|\overline{\Delta}_{J(\theta^*)}|_1 +
\frac{2c}{1-c}\widehat{b}{\overline z}_{*}\\
\nonumber &\le&
\frac{2|J(\theta^*)|}{1-c}\,|\widehat{\theta}^{(c)}-\theta^*|_\infty +
\frac{2c}{1-c}\widehat{b}{\overline z}_{*}.
\end{eqnarray}
Inequalities (\ref{eq:th:nonvalid:1}) and (\ref{eq:th:nonvalid:2})
follow from solving (\ref{eq:proof:nv8}) and (\ref{eq:proof:nv9})
with respect to $|\widehat{\theta}^{(c)}-\theta^*|_\infty$ and
$|\widehat{\theta}^{(c)}-\theta^*|_1$ respectively. \hfill
$\square$\vspace{0.3cm}

%%%%%%%%%%%%%%%%%%%%%%%%%

\noindent{\bf Proof of Theorem \ref{th:nonvalid1}.} We first prove
part (i). We will assume that we are on the event of probability at
least $1-\alpha_1-\alpha_2-\varepsilon$ where (\ref{eq:proof:nv2}),
(\ref{eq:ass_nonvalid2}), and (\ref{eq:th:nonvalid1:1}) are
simultaneously satisfied. From (\ref{eq:proof:nv5}) and the fact
that (\ref{eq:ass_nonvalid2}) can be written as
$F(\theta^*,\beta^*)\le \sigma_{1*}$ we obtain
\begin{eqnarray}\label{eq:proof:nv10}
{\widehat\sigma}_1^{(c)}\le |\widehat{\theta}^{(c)}-\theta^*|_1/c + \sigma_{1*}.
\end{eqnarray}
Note also that the argument in the proof of Theorem
\ref{th:nonvalid} and the results of that theorem remain obviously
valid with $\widehat{b}$ replaced by $b_*$. Thus, we can use
(\ref{eq:th:nonvalid:2}) with $\widehat{b}$ replaced by $b_*$, and
combining it with (\ref{eq:proof:nv10}) we obtain
\begin{eqnarray}\label{eq:proof:nv11}
{\widehat\sigma}_1^{(c)}\le {\overline\sigma}_{*}.
\end{eqnarray}
This and (\ref{eq:th:nonvalid:1}) yield (\ref{eq:th:nonvalid1:2}).
%\begin{eqnarray}\label{eq:proof:nv12}
%|\widehat{\theta}-\theta^*|_\infty &\le& r_1(2\sigma_{1*} +
%(1+2/c)|\widehat{\theta}-\theta^*|_1) +2(1+r_1)b_*{\overline z}_{*}.
%\end{eqnarray}
%Bounding here $|\widehat{\theta}-\theta^*|_1$ by
%(\ref{eq:proof:nv9}) where we replace $\widehat{b}$ by $b_*$, and
%then solving with respect to $|\widehat{\theta}-\theta^*|_\infty$ we
%obtain

We now prove part (ii) of the theorem. In the rest of the proof, we
assume that we are on the event ${\mathcal G}'$ of probability at
least $1-\alpha_1-\varepsilon-\gamma$ where (\ref{eq:proof:nv2}),
(\ref{eq:ass_nonvalid2}), and the events ${\mathcal G}$, ${\mathcal
G}_j$ defined in the proofs of Theorems \ref{t1}, \ref{t1b} are
simultaneously satisfied. Then item (ii) of Theorem \ref{t1b} with
$p=1$ implies (\ref{eq:th:nonvalid1:1}) with $b_*$ defined in
(\ref{eq:th:nonvalid1:bstar}). This and (\ref{eq:th:nonvalid1:2})
easily give part (ii) of the theorem.

To prove part (iii), note that, by Theorem~\ref{t1b}~(i) and
Assumption \ref{ass_select},
\begin{eqnarray}\label{eq:proof:nv12}
\widehat b = \frac{2 \widehat{\sigma}
rs}{\kappa_{1}(s)}\left(1-\frac{r}{\kappa_{1}(s)
}\right)_+^{-1} \le b_* \,
\end{eqnarray}
for $b_*$ defined in (\ref{eq:th:nonvalid1:bstar}). This and
(\ref{eq:proof:nv11}) imply that the threshold $\omega$ satisfies
$\omega \triangleq V(\widehat{\sigma}_1^{(c)},c, \widehat b,
J(\widehat{\theta}))\le V({\overline\sigma}_{*},c, b_*, s_1)\triangleq
\omega^{(c)*}$ on the event~${\mathcal G}'$. On the other hand,
(\ref{eq:th:nonvalid:1}) guarantees that $|\widehat
{\theta}_l^{(c)}-\theta_l^*|\le\omega^{(c)}$ and, by assumption,
$|\theta_l^*|>2\omega^{(c)*}$ for all $l \in J(\theta^*)$. In addition,
by (\ref{eq:t1:1}) and (\ref{eq:th:nonvalid:1}) for all $l \in
J(\theta^*)^c$ we have $|\theta_l^*|<\omega^{(c)}$,  which implies
$\widetilde\theta_l^{(c)}=0$. We finish the proof in the same way as the
proof of Theorem \ref{t1b}. \hfill $\square$\vspace{0.3cm}

%%%%%%%%%%%%%%%%%%%%%%%

\subsection{Proof of Lemma \ref{lem:Ftheta}.}\label{subsec:proof:lem}
\label{subsec:proof:lem:Ftheta}
Set
$f_l(\theta_l)\triangleq\sqrt{\widehat{Q}_l(\theta_l,\widehat{\beta})}$,
and $f(\theta)\triangleq \max_{l=1,\hdots,L_1} f_l(\theta_l) \equiv
F(\theta,\widehat{\beta})$. The mappings $\theta\mapsto
f_l(\theta_l)$ are convex, so that by the Dubovitsky-Milutin, the
subdifferential of their maximum $f$ is contained in the convex hull
of the union of the subdifferentials of the $f_l$:
\begin{equation}\label{eq:dubovitsky}
\partial f \subseteq {\rm Conv}\left(\bigcup_{l=1}^{L_1} \partial f_l
\right).
\end{equation}
Since, obviously, $\partial f_l(\theta_l) \subseteq [-1,1]$, we find
that $\partial f(\theta) \subseteq \{w\in \R^{L_1}: \, |w|_\infty\le
1\}$ for all $\theta\in \R^{L_1}$. Using this property and the
convexity of $f$, we get
$$
f(\theta^*) - f(\widehat{\theta}) \le \langle w,\theta^*-
\widehat{\theta}\rangle \le |\widehat{\theta}-\theta^*|_1, \quad
\forall \ w\in \partial f(\theta^*),
$$
where $\langle \cdot,\cdot \rangle$ denotes the standard inner
product in $\R^{L_1}$. This yields (\ref{eq:proof:nv3}). The proof
of (\ref{eq:proof:nv4}) is based on similar arguments. Instead of
$f_l$, we now introduce the functions $g_l$ defined by
$g_l(\beta)\triangleq\sqrt{\widehat{Q}_l(\theta_l^*,\beta)}$, and
set $g(\beta)\triangleq \max_{l=1,\hdots,L_1} g_l(\beta) \equiv
F(\theta^*,\beta)$. Next, notice that the subdifferential of $g_l$
satisfies $\partial g_l(\beta) \subseteq \{w\in \R^{K}: \, |w_k|\le
a_{lk}, \ k=1,\dots,K\}$ for all $\beta\in \R^{K}$, $l=1,\dots,L_1$,
where
$$
a_{lk}=\frac{\left|\frac1n\sum_{i=1}^n\overline{z}_{li}x_{ki}
\left(\overline{z}_{li}
(y_i-x_i^T\beta)-\theta_l^*\right)\right|}{\sqrt{\frac1n\sum_{i=1}^n\left(\overline{z}_{li}
(y_i-x_i^T\beta)-\theta_l^*\right)^2}}\,.
$$
Consequently, by the Cauchy-Schwarz inequality, ${\bf D_X}\partial
g_l(\beta) \subseteq \{w\in \R^{K}: \, |w|_\infty\le {\overline
z}_{*}\}$ for all $\beta\in \R^{K}$, $l=1,\dots,L_1$. This and
(\ref{eq:dubovitsky}) with $g$, $g_l$ instead of $f$, $f_l$ imply
${\bf D_X}\partial g(\beta) \subseteq \{w\in \R^{K}: \,
|w|_\infty\le {\overline z}_{*}\}$ for all $\beta\in \R^{K}$.  Using
this property and the convexity of $g$, we get
$$
g({\beta})- g(\beta') \le \langle w, ({\beta}-\beta')\rangle \le
|{\bf D_X}w|_\infty \left|{\bf D_X}^{-1}({\beta}-\beta')\right|_1
\le {\overline z}_{*}\left|{\bf D_X}^{-1}({\beta}-\beta')\right|_1,
\quad \forall \ w\in
\partial g({\beta}),
$$
for any ${\beta}, \beta'\in \R^K$. This proves (\ref{eq:proof:nv4}).
\hfill $\square$\vspace{0.3cm}

\noindent{\bf Proof of Theorem \ref{tRWI}.}  The only difference
with the proof of Theorem \ref{t1} is that, because we do not have
the $l_1$-norm in the objective function \eqref{IVSRWI}, we drop the
discussion leading to \eqref{eq:conec}.\hfill $\square$

\end{document}